%% file: RationalApprox.tex
\documentclass[review]{elsarticle}

\input{ex_shared}

\ifpdf
\hypersetup{
	pdftitle={Multivariate Rational Approximation},
	pdfauthor={Austin, Krishnamoorthy, Leyffer, Mrenna, M\"uller, and Schulz}
}
\fi

\begin{document}

\begin{frontmatter}
	

	\title{Multivariate Rational Approximation}
	\input{authorandaff}
		
	\begin{abstract}
		We present two approaches for computing rational approximations to multivariate
		  functions, motivated by their effectiveness as surrogate models for high-energy
		  physics (HEP) applications. Our first approach builds on the Stieltjes process to
		  efficiently and robustly compute the coefficients of the rational approximation.
		  Our second approach is based on an optimization formulation that allows us to include
		  structural constraints on the rational approximation, resulting in a semi-infinite
		  optimization problem that we solve using an outer approximation approach. We present
		  results for synthetic and real-life HEP data, and we compare the approximation
		  quality of our approaches with that of traditional polynomial approximations.
	\end{abstract}
	
	
\end{frontmatter}


%
%
%


\input{intro}
\input{la}
\input{opti}
\input{numerics}

\input{apps}
\input{concl}

\section*{Acknowledgments}

This work was supported by the U.S. Department of Energy, Office of Science, Advanced Scientific Computing Research, under Contract DE-AC02-06CH11357.
Support for this work was provided through the SciDAC program funded by U.S. Department of Energy, Office of Science, Advanced Scientific Computing Research.
This work was also supported by
the U.S. Department of Energy through grant DE-FG02-05ER25694, and
by Fermi Research Alliance, LLC under Contract No. DE-AC02-07CH11359 with the U.S. Department of Energy, Office of Science, Office of High Energy Physics.
This work was supported in part by the U.S. Department of Energy, Office of Science, Office of Advanced Scientific
Computing Research and Office of Nuclear Physics, Scientific Discovery through Advanced Computing (SciDAC) program through the FASTMath Institute under Contract No. DE-AC02-05CH11231 at Lawrence Berkeley National Laboratory.

\bibliography{NLP,robust,minlp,jshort,ndrathep,ra,hep}

\appendix
\input{testProblems}

\input{multistartstrategy}


\vfill
\begin{flushright}
\scriptsize
\framebox{\parbox{\textwidth}{
The submitted manuscript has been created by UChicago Argonne, LLC, Operator of Argonne National Laboratory (“Argonne”). 
Argonne, a U.S. Department of Energy Office of Science laboratory, is operated under Contract No. DE-AC02-06CH11357. 
The U.S. Government retains for itself, and others acting on its behalf, a paid-up nonexclusive, irrevocable worldwide 
license in said article to reproduce, prepare derivative works, distribute copies to the public, and perform publicly 
and display publicly, by or on behalf of the Government.  The Department of Energy will provide public access to these 
results of federally sponsored research in accordance with the DOE Public Access Plan. 
\url{http://energy.gov/downloads/doe-public-access-plan}.
}}
\normalsize
\end{flushright}

\end{document}


\begin{frontmatter}
	\title{Supplementary Materials: Multivariate Rational Approximation}
	
	\input{authorandaff}

\end{frontmatter}

\makeatletter
\def\@cline#1-#2\@nil{%
	\omit
	\@multicnt#1%
	\advance\@multispan\m@ne
	\ifnum\@multicnt=\@ne\@firstofone{&\omit}\fi
	\@multicnt#2%
	\advance\@multicnt-#1%
	\advance\@multispan\@ne
	\leaders\hrule\@height\arrayrulewidth\hfill
	\cr
	\noalign{\nobreak\vskip-\arrayrulewidth}}
\makeatother

\input{suppPolesAndError}
\input{suppTimesAndIterations}

\clearpage

%% file: ex_shared.tex
\usepackage{longtable}
\usepackage{amsmath}
\usepackage{amsfonts}
\usepackage{mathtools}
\usepackage{algorithm}
\usepackage{algorithmic}
\usepackage{subcaption}
\usepackage{multirow}
\usepackage{xr}

\usepackage[algo2e,linesnumbered,ruled,vlined,boxed,noresetcount,algosection]{algorithm2e}
\usepackage{tikz}
\usetikzlibrary{quotes,angles}
\usepackage{graphicx}
\usetikzlibrary{positioning}

\usepackage{pgf}
\usepackage{rotating}
\newcommand{\R}{\mathbb{R}}
\newcommand{\sP}{\mathcal{P}}
\newcommand{\ip}[2]{\left\langle #1, #2 \right\rangle}
\DeclareMathOperator*{\diag}{\mathrm{diag}}
\DeclareMathOperator*{\Ran}{\mathrm{Ran}}
\usepackage[normalem]{ulem}
\usepackage{amssymb, palatino, geometry}
\usepackage{latexsym}
\usepackage{epsfig,graphicx,xcolor}
\usepackage{natbib, epsfig} 
\usepackage{refcheck}           
\usepackage{url}                
\usepackage{version}            
\usepackage{tcolorbox}
\usepackage{longtable}
\graphicspath{{./images/}}
\newcommand{\mini}{\mathop{\mbox{minimize}}}

\newcommand{\st}{\mbox{subject to }}

\newcommand{\dps}{\displaystyle}

\newcommand{\xk}{x^{(k)}}

\usepackage{makecell}
\usepackage{array}
\usepackage{url}
\usepackage{hyperref}
\hypersetup{breaklinks=true}
\newcounter{rowcntr}[table]
\renewcommand{\therowcntr}{\thetable.\arabic{rowcntr}}
\newcolumntype{N}{>{\refstepcounter{rowcntr}\therowcntr}c}
\AtBeginEnvironment{tabular}{\setcounter{rowcntr}{0}}
\usepackage[nameinlink]{cleveref}
\crefname{appendix}{}{}
\Crefname{appendix}{}{}
\crefname{table}{\mbox{Table}}{\mbox{Table}}
\Crefname{table}{\mbox{Table}}{\mbox{Table}}
\crefname{figure}{\mbox{Fig.}}{\mbox{Fig.}}
\Crefname{figure}{\mbox{Fig.}}{\mbox{Fig.}}
\crefname{algocf}{\mbox{Algorithm}}{\mbox{Algorithms}}
\Crefname{algocf}{\mbox{Algorithm}}{\mbox{Algorithms}}
\crefname{equation}{\mbox{Eq.}}{\mbox{Eq.}}
\Crefname{equation}{\mbox{Eq.}}{\mbox{Eq.}}
\crefname{algorithm}{\mbox{Algorithm}}{\mbox{Algorithms}}
\Crefname{algorithm}{\mbox{Algorithm}}{\mbox{Algorithms}}
\crefname{section}{\mbox{section}}{\mbox{sections}}
\Crefname{section}{\mbox{Section}}{\mbox{Sections}}
\crefname{rowcntr}{\mbox{Function}}{\mbox{Functions}}
\Crefname{rowcntr}{\mbox{Function}}{\mbox{Functions}}
\usepackage{lineno,hyperref}
\modulolinenumbers[5]

\norefnames
\nocitenames


%% file: authorandaff.tex
	\author[vt]{Anthony P.~Austin}
	\author[anl]{Mohan Krishnamoorthy\corref{mycorrespondingauthor}}
	\cortext[mycorrespondingauthor]{Corresponding author}
	\ead{mkrishnamoorthy@anl.gov}
	\author[anl]{Sven Leyffer}
	\author[fnl]{Stephen Mrenna}
	\author[lbnl]{Juliane M\"uller}
	\author[uct]{Holger Schulz}
		
	\address[vt]{Department of Mathematics, Virginia Tech, Blacksburg, VA 24061.}
	\address[anl]{Mathematics and Computer Science Division, Argonne National Laboratory, Lemont, IL 60439}
	\address[fnl]{Fermi National Accelerator Laboratory, Batavia, IL 60510}
	\address[lbnl]{Lawrence Berkeley National Laboratory, Berkeley, CA 94720}
	\address[uct]{Department of Physics, University of Cincinnati, Cincinnati, OH 45219}

%% file: intro.tex
\section{Introduction}\label{sec:intro}

Optimization problems arising in complex science and engineering applications
often involve simulations that are computationally expensive to evaluate
(several minutes to hours or more per evaluation). The simulations are usually
nonlinear and black box; in other words, we have no analytic description of the function
$f(\cdot)$ that maps the parameter inputs $x\in D\subset\mathbb{R}^n$ to
simulation outputs.  The computational expense limits the number of evaluations
we can do during the optimization.  A widely used approach to mitigate this
difficulty is to use a fast-to-evaluate \textit{surrogate model}, $s(x)$, as a
proxy for the simulation~\citep{Booker1999}: $f(x) =s(x)+e(x)$, where $e(x)$
denotes the difference between the true function and the surrogate model.  We
fit a surrogate model based on a set of pre-evaluated parameter--function value
pairs and use it during the optimization search, thus reducing the number of
queries to the expensive simulation. Types of surrogate models include Gaussian
process models~\citep{Matheron1963}, radial basis functions~\citep{Powell1992},
multivariate adaptive regression splines~\citep{Friedman1991}, and polynomial
regression models~\citep{Myers1995}.

Polynomial models have several advantages, such as a simple representation and
being easy to build and use; however, they have poor extrapolation behavior and
are severely limited in their ability to cope with singularities.  These
drawbacks can reduce their effectiveness at representing elements of the
physics in many applications.  Because of these drawbacks, one 
turns to models based on rational functions (quotients of polynomials), whose
ability to capture singularities naturally via their poles can make them
considerably more powerful than polynomials \citep{Devore1986,New1964}.
Unfortunately, rational approximations can be numerically fragile to compute
and are prone to having spurious singularities.  Moreover,  how to select the appropriate combination of numerator and denominator
degree is not always clear.

In this article, we investigate the utility of rational approximations as
surrogate models.  We propose two methods for computing
multivariate rational approximations $r(x)=p(x)/q(x)$. The first approach is
based on the univariate methods of \citep{GPT11,PGV12} and provides a robust and
efficient way to compute the coefficients of $p(x)$ and $q(x)$.  Although it
tries to reduce the propensity for the resulting $r(x)$ to contain unwanted
singularities by using ideas from linear algebra to minimize the degree of
$q(x)$, it does  \emph{not guarantee} that $r(x)$ will be pole free in the
parameter domain.

The second approach uses a constrained optimization formulation that includes
structural constraints on $r(x)$ to enforce the absence of poles in $D$. These
constraints are motivated from applications that arise, for example, in high-energy
physics (HEP) simulations.  Although it is computationally more expensive than the
first approach, the second approach allows us to guarantee that the computed
approximation is free of poles in the parameter space, which can be crucial in
the context of computing surrogate models for use in optimization.

\subsection{Previous Work on Rational Approximation}

The literature on rational approximation is too vast to cover comprehensively
here; we refer the reader to standard texts such as \citep[Ch.\ V]{Bra1986},
\citep[Ch.\ 5]{Che1966}, and \citep[Ch.\ 23--27]{Tre2013} for the history and
basic concepts.  In this article, we are concerned with multivariate models.
Multivariate rational approximation has been studied extensively by Cuyt and
co-authors \citep{Cuy1983,CV1985}; in particular, Cuyt and Yang have recently
developed practical error bounds \citep{Cuyt2010}.

Our rational approximations are least-squares models, which can be formulated
in a nonlinear or linearized fashion as we describe below.  Our first approach
is an extension of the algorithms for the linearized problem presented in
\citep{GPT11,PGV12,SCGPS2013}.  The nonlinear problem is an example of a
\emph{separable} nonlinear least-squares problem, and algorithms for it often
exploit this structure.  Examples include the Gauss-Newton algorithm developed
by Golub and Pereyra \cite{GP1973} and the full-Newton algorithm of Borges
\cite{Bor2009}.  The so-called ``AAA'' algorithm of \citep{NST2018} is a
particularly interesting recently developed alternative to traditional methods
for rational least-squares problems; however, as it currently only works in the
univariate context, we do not explore it further here.  Of course,
least-squares approximations are not the only type of rational model.  Recent
work of interest on rational models other than simple least-squares
approximations includes the rational minimax approximation algorithm of
\citep{FNTB2018}.

One of the appeals of a least-squares approach to rational approximation is
that it is naturally robust to noise in the data being fit.  Salazar Celis,
Cuyt, and Verdonk \citep{Salazar2007} have developed an alternative approach to
rational approximation in the face of noisy or uncertain data that computes an
uncertainty interval for each datum and then solves a quadratic optimization
problem to find a rational function that passes through all uncertainty
intervals.

\subsection{High-Energy Physics Motivation}
Our work is motivated by simulations for studying complex physical phenomena,
especially in high-energy physics. Simulations are often used to  guide our
real-world experiments in order to find ``interesting'' physics or  to verify
that the models we derive from our physics understanding  are in agreement with
the experiments \citep{Buckley:2011ms}. However, these simulations (as well as physics experiments)
are generally resource intensive (computationally or other) \citep{Alves:2017she}. A single
simulation may require many hours of compute time on a modern supercomputer,
thus limiting the number of simulation runs that can realistically be done.


This fact severely limits many applications that require extensive parameter
space exploration. Our aim is to replace the costly simulations with rational approximations that are  much
cheaper to evaluate. In particular, we want to
construct and numerically optimize an objective function over a space of model parameters that is defined as the mismatch between
experimental data and simulation predictions.

\subsection{Outline of the Paper}
In \cref{SEC:Setup}, we establish our notation and describe the types of
models that we will generate. In \cref{SEC:LinAlg}, we devise  a method
for constructing rational approximation  models based on linear algebra. This approach is
flexible and easy to implement, but it has the drawback that  singularities may be present.  Although singularities are acceptable in  some contexts, we generally have to prevent  singularities in particular regions of the parameter space because they may cause an unbounded  objective function in our optimization
procedure, which is not acceptable.
In \cref{SEC:opti}, we describe a separate approach based on
semi-infinite optimization that allows us to achieve this goal.
In \cref{SEC:Numerics}, we present some numerical results, and in \cref{SEC:HEP} we describe our high-energy physics application and show the superior performance of our pole-free rational approximations  over polynomial approximations and rational approximations with poles.  
In \cref{SEC:Conclusion}, we summarize our key findings and discuss potential
avenues for further research.

\section{Notation and Setup}
\label{SEC:Setup}

We denote by $n$ the number of parameters in our model, and our generic
variables are $x_1, \ldots, x_n$.  By the \emph{degree} of an $n$-variate
monomial $x_1^{i_1} \cdots x_n^{i_n}$, we mean its \emph{total degree}, in other words,
the sum $i_1 + \cdots + i_n$, as distinguished from its \emph{maximal degree}
$\max(i_1, \ldots, i_n)$.  The degree of an $n$-variate polynomial is the
maximum of the degrees of its constituent monomials.  We write $\sP_d^n$ for
the space of all $n$-variate polynomials of degree at most $d$; this is a real
vector space of dimension $\alpha(d) = \binom{n + d}{d}$.

Let $x^{(0)}, \ldots, x^{({K - 1})}$ be $K$ points in $\R^n$, and let $f_0, \ldots, f_{K - 1}$
be the corresponding real data values.  Our aim is to find an
$n$-variate rational function $r(x) = p(x)/q(x)$ with $p \in \sP_M^n$ and $q
\in \sP_N^n$ such that $r(x^{(k)}) \approx f_k$  for each $k$.  One natural
approach is to choose $p$ and $q$ to solve the discrete least-squares problem
\begin{equation}\label{EQN:NLSQ}
  \mini_{p,q} \quad \sum_{k = 0}^{K - 1} \left(\frac{p(\xk)}{q(\xk)} - f_k\right)^2 \quad
  \st \quad p \in \sP_M^n, q \in \sP_N^n,
\end{equation}
but the nonlinearity in $q$ makes this problem challenging.  It is usually easier to
 work with the linearized problem
\begin{equation}\label{EQN:LSQ}
  \mini_{p,q} \quad \sum_{k = 0}^{K - 1} \left(p\left(\xk\right) - f_k q\left(\xk\right)\right)^2 \quad
  \st \quad p \in \sP_M^n, q \in \sP_N^n,
\end{equation}
which is the formulation we will use in the following.  Note that the solutions to
\cref{EQN:NLSQ} and \cref{EQN:LSQ} do not generally coincide.

As written, \cref{EQN:NLSQ} and \cref{EQN:LSQ} are incompletely specified.
The objective in \cref{EQN:NLSQ} depends only on the ratio of $p$ to $q$, and  additional normalization conditions must be imposed to pin down the
solution.  Likewise, a normalization condition is needed in \cref{EQN:LSQ} to
exclude the trivial solution $p = q = 0$.  We will address these issues in detail  in later sections.

%% file: la.tex
\section{Multivariate Rational Models via Linear Algebra}
\label{SEC:LinAlg}

Our first approach to constructing rational models is based on ideas from
linear algebra following \cite{GPT11,PGV12}; see also \citep{GGT13} and
\cite[Ch.\ 26]{Tre2013}.  We extend the method of these references to the
multivariate case.  One such extension has been proposed in \cite{SCGPS2013};
our method may be viewed as a generalization of that extension to handle
situations in which the data used to construct the model come from arbitrary
sample points in the parameter space instead of from a tensor product grid.

\subsection{Basic Algorithm}

The basic idea is as follows.  Let $L = \max(M, N)$.  Given a basis $\varphi_0,
\ldots, \varphi_{\alpha(L) - 1}$ for $\sP_L^n$, consider the Vandermonde-like
matrices $V_M \in \R^{K \times \alpha(M)}$ and $V_N \in \R^{K \times
\alpha(N)}$ whose $(k, j)$ entries are $\varphi_j(x^{(k)})$.  Express $p$ and
$q$ as
\[
p(x) = \sum_{j = 0}^{\mathclap{\alpha(M) - 1}} a_j \varphi_j(x), \qquad q(x) = \sum_{j = 0}^{\mathclap{\alpha(N) - 1}} b_j \varphi_j(x),
\]
and gather the coefficients $a_j$, $b_j$ into vectors $a \in \R^{\alpha(M)}$
and $b \in \R^{\alpha(N)}$, respectively.  Let $F = \diag(f_0, \ldots f_{K -
1})$.  Then, since the $k$th entries of $V_M a$ and $V_N b$ are $p(\xk)$ and
$q(\xk)$, respectively, the linearized problem \cref{EQN:LSQ} may be
rewritten as the following linear least-squares problem to find coefficients, $a,b$, that
\begin{equation}\label{EQN:LSQ-LA1}
\mini_{a,b} \quad \|V_M a - F V_N b\|_2^2.
\end{equation}
Just as \cref{EQN:LSQ} has the trivial solution $p = q = 0$,
\cref{EQN:LSQ-LA1} has the trivial solution $a = 0$, $b = 0$.  To forbid this
solution, we impose the normalization condition $\|b\|_2 = 1$.  If the choice
of $b$ needed to solve \cref{EQN:LSQ-LA1} is known, then the corresponding
choice of $a$ is given by $a = Z b$, where $Z = V_M^+ F V_N$ and $V_M^+$ is the
Moore-Penrose pseudoinverse of $V_M$.  Substituting this relationship into
\cref{EQN:LSQ-LA1}, we are left with the problem to find the coefficients, $b$,
of the denominator that 
\begin{equation}\label{EQN:LSQ-LA2}
  \mini_b \quad \|(V_M Z - F V_N) b\|_2^2 \quad \st \quad \|b\|_2 = 1,
\end{equation}
and this may be solved by taking $b$ to be the right singular vector
corresponding to the smallest singular value of $W = V_M Z - F V_N$.

If $K = \alpha(M) + \alpha(N) - 1$, then the number of data points matches the
number of degrees of freedom in $p$ and $q$, less 1 for the normalization
condition.  In this case, we expect that the objective in \cref{EQN:LSQ-LA2}
can be driven to zero, yielding a linearized rational interpolant to the data,
sometimes called a \emph{multipoint Pad\'{e} approximation}.  We write $W = (V_M
V_M^+ - I)FV_N$, where $I$ is the identity matrix.  Since $V_M V_M^+ - I$ is
($-1$ times) the orthogonal projector onto $\Ran(V_M)^\perp$, it has rank at
most $K - \alpha(M) = \alpha(N) - 1$.  Since $W$ is of size $K \times
\alpha(N)$, this implies that $W$ is rank deficient---it has at least one zero
singular value---so $b$ can indeed be chosen to satisfy \cref{EQN:LSQ-LA2}
with an objective value of zero, as expected.

\subsection{Discrete Multivariate Orthogonal Polynomials}

While in principle one can use any basis $\varphi_0, \ldots,
\varphi_{\alpha(L) - 1}$ for $\sP_L^n$, some bases are better  suited to
numerical computation than are others.  In particular, it is important that the
basis be chosen so that the Vandermonde-like matrices $V_M$ and $V_N$ are
well conditioned.  We would ideally choose the basis so that $V_M$ and $V_N$
have orthonormal columns; in addition to ensuring that operations involving
these matrices are robust to rounding error, this would make working with the
pseudoinverse of $V_M$ trivial, because we would have $V_M^+ = V_M^*$.  We can
accomplish this by choosing $\varphi_0, \ldots, \varphi_{\alpha(L) - 1}$ so
that they are orthogonal with respect to the discrete inner product%
\footnote{This will be an inner product only if the $\xk$ constitute a set of
linear independence for $\sP_L$.  We assume throughout this article that this
is true.}
\begin{equation}\label{EQN:DiscreteIP}
\ip{h}{g} = \sum_{k = 0}^{K - 1} h(x^{(k)}) g(x^{(k)})
\end{equation}
on $\sP_L^n$ associated with the sample points $\xk$.  The orthogonality
condition $\ip{\varphi_i}{\varphi_j} = \delta_{ij}$, where $\delta_{ij}$ is the
Kronecker delta, is precisely the statement that $V_M$ and $V_N$ have
orthonormal columns.

One way to construct such a basis is via a multivariate version of the familiar
\emph{Stieltjes process}~\citep{carter2000lebesgue} from the theory of
(univariate) orthogonal polynomials.  Discussions may be found elsewhere in the
literature---see, for example, \cite{HL2002,Zac2014}---but to keep this paper
self-contained, we describe the process in the form in which we use it here.

The Stieltjes process may be viewed as a variant of the Gram-Schmidt process
that orthogonalizes the columns of a Vandermonde matrix without performing the
numerically unsavory operation of evaluating high-order monomials, that is,
without explicitly forming the matrix itself.  In a single variable, it works
as follows.  We begin by assigning $\varphi_0(x) = 1/\ip{1}{1}$.  Then, having
constructed $\varphi_0, \ldots, \varphi_{j - 1}$, we construct $\varphi_j$ by
orthogonalizing $x \varphi_{j - 1}(x)$ against $\varphi_0, \ldots, \varphi_{j -
1}$,
\begin{equation}\label{eq:orthogonalize}
\hat{\varphi}_j(x) = x\varphi_{j - 1}(x) - \sum_{i = 0}^{j - 1}\ip{x\varphi_{j - 1}}{\varphi_i} \varphi_i(x),
\end{equation}
and normalizing,
\begin{equation}\label{eq:normalize}
\varphi_j(x) = \frac{\hat{\varphi}_j(x)}{\sqrt{\ip{\hat{\varphi}_j}{\hat{\varphi}_j}}}.    
\end{equation}
Since the operation of multiplication by $x$ is self-adjoint (i.e.,
$\ip{x\varphi}{\psi} = \ip{\varphi}{x\psi}$ for all $\varphi$, $\psi$), the
orthogonality condition can be used to show that only the $i = j - 1$ and $i =
j - 2$ terms in the sum for $\hat{\varphi}_j$ are nonzero, leading to a
three-term recurrence relation for $\varphi_j$.  This recurrence can be used to
evaluate polynomials that are expressed as linear combinations of the
$\varphi_j$ at arbitrary points.

The multivariate case works similarly.  The key difference is that since there
is no canonical ordering of the monomials in several variables---no agreed-upon
order in which to list the columns of a multivariate Vandermonde matrix---we
must first select one and then develop a version of the Stieltjes process
tailored to that ordering.  The ordering we use is as follows.  We say that
$x_1^{i_1} \cdots x_n^{i_n} < x_1^{j_1} \cdots x_n^{j_n}$ if $i_1 + \cdots +
i_n < j_1 + \cdots + j_n$ or if $i_1 + \cdots + i_n = j_1 + \cdots + j_n$ and
$i_k > j_k$, where $k$ is the smallest index such that $i_k \neq j_k$.  For
instance, in $n = 3$ variables $x_1 = x$, $x_2 = y$, and $x_3 = z$, this
ordering lists the monomials of degree $3$ or less in the following sequence:
\[
1, x, y, z, x^2, xy, xz, y^2, yz, z^2, x^3, x^2y, x^2z, xy^2, xyz, xz^2, y^3, y^2z, yz^2, z^3.
\]
This order is related to the popular ``grevlex'' order \cite[Sec.\
2.2]{CLS2010} and has two features that make it convenient.  One is that the
monomials are ordered by degree.  The other is that it yields a simple
inductive process for listing the monomials in sequence.  This is most easily
described by example.  To construct the three-variable sequence above, we begin
with the constant monomial $1$.  We then multiply $1$ by each of the variables
in order to obtain the three linear monomials $x$, $y$, and $z$.  To produce
the quadratic monomials, we first multiply each of the linear monomials by $x$,
retaining the order, to produce $x^2$, $xy$, and $xz$.  We then multiply by $y$ 
the
linear monomials that do not contain $x$, giving $y^2$ and $yz$.
Finally, we multiply the linear monomials that contain neither $x$ nor $y$---in
other words, $z$---by $z$, giving $z^2$.  The cubic monomials are constructed
similarly.  We multiply all of the quadratic monomials by $x$ to obtain $x^3$
through $xz^2$.  Then, we multiply by $y$ the quadratic monomials that do not contain
$x$, giving $y^3$ through $yz^2$.  Finally, we multiply by $x$ the lone
quadratic monomial that contains neither $x$ nor $y$ to produce $z^3$.

The multivariate Stieltjes process that we use is a straightforward outgrowth
of this construction.  We associate each orthogonal polynomial $\varphi_j$ with
its corresponding term in the monomial sequence, beginning with the association
$\varphi_0 \leftrightarrow 1$.  Having constructed $\varphi_0, \ldots,
\varphi_j$, we construct $\varphi_{j + 1}$ by multiplying the appropriate
previously constructed polynomial by the appropriate variable and
orthogonalizing.  For instance, taking the three-variable case as an example
once more, to produce $\varphi_{12}$, which is associated with the monomial
$x^2 z$, we multiply $\varphi_6$, which is associated with the monomial $xz$,
by $x$ and then orthogonalize the result against $\varphi_0, \ldots,
\varphi_{11}$.

\sloppy 
This process can be easily adapted to compute the Vandermonde-like matrix $V_L$
corresponding to $\varphi_0, \ldots, \varphi_{\alpha(L) - 1}$, which is what we
really want, rather than the polynomials themselves.  Such a version of the
process is given in \cref{ALG:MVVandQR}.  In implementing this
procedure in finite-precision arithmetic, all of the standard caveats about the
numerical stability of the Gram-Schmidt processes apply.  In particular, some
form of reorthogonalization is mandatory to ensure that the columns of the
computed $V_L$ are orthogonal to working precision.  In our implementation, we
use the standard technique of performing the orthogonalization twice, which is
usually sufficient \citep{GLRE2005,LBG2013}, \cite[\S6.9]{Par1998}.


\begin{algorithm2e}[htb!]
	\SetNoFillComment \DontPrintSemicolon \setcounter{AlgoLine}{0}
	\caption{Multivariate Stieltjes Process for Vandermonde-like Matrix}\label{ALG:MVVandQR}
	\SetKwInOut{Input}{Input}
	\SetKwInOut{Output}{Output}
	\Input{Points $x^{(0)}, \ldots, x^{(K - 1)} \in \R^n$ that are a set of linear independence for $\sP_L^n$.}
	\Output{Vandermonde-like matrix $V_L$ corresponding to a basis $\varphi_0, \ldots, \varphi_{\alpha(L) - 1}$ for $\sP_L^n$, orthonormal with respect to \cref{EQN:DiscreteIP}, and coefficients $r_{i, j}$ for use in the recurrence of \cref{ALG:MVRecurrence}.}
	$i \gets 1$\;
	\For{$j = 1$ to $n + 1$}{
		$i_j \gets 0$ \tcc*{$i_j$ marks start of sequence last multiplied by $x_j$.}
	}
	$v_0 = \begin{bmatrix}1 & \cdots & 1\end{bmatrix}^T/\sqrt{K}$ \tcc*{Begin with constant polynomial.}
	\For{$d = 1$ to $L$}{
		\For{$j = 1$ to $n$}{
			$i^* \gets i$\;
			\For{$k = i_j$ to $i_{n + 1}$}{
				$\hat{v}_i \gets \diag(x^{(0, j)}, \ldots, x^{(K - 1, j)}) v_k$ \tcc*{Multiply by $x_j$.}
				\For(\tcc*[f]{Orthogonalize (Gram-Schmidt) \cref{eq:orthogonalize}.}){$\ell = 0$ to $i - 1$}{
					$r_{\ell, i} \gets v_\ell^*\hat{v}_i$\;
					$\hat{v}_i \gets \hat{v}_i - r_{\ell, i}v_\ell$\;
				}
				$r_{i, i} \gets \sqrt{\hat{v}_i^*\hat{v}_i}$ \tcc*{Normalize \cref{eq:normalize}}
	 			$v_i \gets \hat{v}_i/r_{i, i}$\;
				$i \gets i + 1$\;
			}
			$i_j \gets i^*$ \tcc*{Update bookkeeping information.}
		}
		$i_{n + 1} \gets i - 1$\;
	}
	$V_L \gets \begin{bmatrix}v_0 & \cdots & v_{\alpha(L) - 1}\end{bmatrix}$\;
\end{algorithm2e}

%


Like the univariate Stieltjes process, the multivariate process produces a
recurrence relation that can be used to evaluate polynomials expressed in the
generated orthogonal basis at arbitrary points.  Unlike the univariate
recurrence, the multivariate recurrence cannot be reduced to three terms, but
it may possess other structure depending on the monomial ordering that is used.
The recurrence generated by \cref{ALG:MVVandQR} is presented in
\cref{ALG:MVRecurrence}.


\begin{algorithm2e}[htb!]
	\DontPrintSemicolon \setcounter{AlgoLine}{0}
	\caption{Recurrence for Evaluating Multivariate Orthogonal Polynomial Series}\label{ALG:MVRecurrence}
	\SetKwInOut{Input}{Input}
	\SetKwInOut{Output}{Output}
	\Input{Evaluation point $x \in \R^n$, expansion coefficients $c_0$, \ldots, $c_{\alpha(L) - 1}$, recurrence coefficients $r_{i, j}$ from  \cref{ALG:MVVandQR}.}
	\Output{$s = c_0 \varphi_0(x) + \cdots + c_{\alpha(L) - 1}\varphi_{\alpha(L) - 1}(x)$, where the $\varphi_i \in \sP_L^n$ are the orthogonal polynomials associated with the Vandermonde matrix constructed by \cref{ALG:MVVandQR}.}
	$i \gets 1$\;
	\For{$j = 1$ to $n + 1$}{
		$i_j \gets 0$ \tcc*{$i_j$ marks start of sequence last multiplied by $x_j$.}
	}
	$y_0 \gets 1/\sqrt{\ip{1}{1}}$ \tcc*{Begin with $y_0 = \varphi_0(x)$ (constant).}
	\For{$d = 1$ to $L$}{
		\For{$j = 1$ to $n$}{
			$i^* \gets i$\;
			\For(\tcc*[f]{Recurrence for $y_i = \varphi_i(x)$.}){$k = i_j$ to $i_{n + 1}$}{
				$\hat{y}_i \gets x_j y_k$\;
				\For{$\ell = 0$ to $i - 1$}{
					$\hat{y}_i \gets \hat{y}_i - r_{\ell, i}y_\ell$\;
				}
				$y_i \gets \hat{y}_i/r_{i, i}$\;
				$i \gets i + 1$\;
			}
			$i_j \gets i^*$ \tcc*{Update bookkeeping information.}
		}
		$i_{n + 1} \gets i - 1$\;
	}
	$s \gets c_0 y_0 + \cdots + c_{\alpha(L) - 1}y_{\alpha(L) - 1}$ \tcc*{Evaluate $s = c_0 \varphi_0(x) + \cdots + c_{\alpha(L) - 1}\varphi_{\alpha(L) - 1}(x)$.}
\end{algorithm2e}


\subsection{Spurious Poles and Degree Reduction}\label{sec:degred}


Rational approximations are powerful because of their ability to capture
singularities in the function being approximated with singularities of their
own; however, if the approximations are computed naively, one often finds that
they possess singularities that bear little resemblance to those of the
function under consideration.  This can happen even when approximating 
well-behaved functions of a single variable, where the unwanted singularities
in the approximation are known as \emph{spurious poles} or \emph{Froissart
doublets}. This is a serious problem:  in our context, an unwanted singularity
in the surrogate model leads to an unbounded objective for our optimization
procedure, which can be problematic.

Unwanted singularities can be broadly classified into two types:  those that
arise in the mathematics and those that arise from noise and numerical
artifacts.  It seems that little can be done about the former; sometimes the solution to the least-squares problem
\cref{EQN:LSQ-LA2} really does have a singularity in an undesirable location
that does not clearly correspond to a singularity of the function being
approximated.  Unwanted singularities of the latter type usually emerge when
the approximation has more degrees of freedom than are necessary to fit the
given data.  One advantage of the construction just described is that it
affords a natural way to handle this situation.  This technique was first
described in \cite{GPT11} for univariate approximation; we extend this idea to
the multivariate case.

Our construction calls for computing $b$ in \cref{EQN:LSQ-LA2} as the right
singular vector of $W = V_M Z - F V_N$ corresponding to the smallest singular
value.  If this singular value is nearly zero, then our rational approximation
will fit the data nearly exactly.  If $W$ has many singular values that are
nearly zero, then there are many possible choices for $b$---and thus many
possible rational approximations---that will have this property.  The key idea
is this:  \emph{If there are many approximations that will work, one should use
the approximation with the lowest-degree denominator.}  In one dimension, reducing the
degree of the denominator by 1 reduces the number of poles of the approximation
by 1.  If the approximation is already fitting the data well, it is highly
likely that the pole that will be eliminated is a spurious one.

Multivariate rational approximations are more complicated than univariate ones:
their singularities may not be isolated, and even if we eliminate unnecessary
degrees of freedom, they will still, in general, have uncountably many singular
points.  As such, it is too much to hope that the degree-reduction approach to
eliminating unwanted singularities will work as well as it does in the
univariate case, especially with noisy input data.  We will see this in some of
our later experiments.  Nevertheless, it can still be highly effective.

The procedure we recommend is summarized in \cref{ALG:ReduceDegree}.
The algorithm attempts to reduce the denominator degree by $1$, checking to see
whether this is possible by examining the smallest singular value of the $W$ matrix
associated with the reduced degree.  If this singular value is smaller than a
chosen threshold $\eta$, it deems the reduction successful and then
tries to reduce the degree by $1$ again.  It continues until the smallest
singular value of $W$ is too large for the reduction to be considered viable.
It then repeats the process to reduce the degree of the numerator by
considering the problem of fitting an approximation to the reciprocal data.

By considering the nullity of $W$, we can reduce the degree in steps
 greater than $1$:  if $W$ has many singular values that lie below the
threshold, we could eliminate many degrees of freedom simultaneously.
Nevertheless, we have found that the stated approach is more robust, especially
in the presence of noise.  For this procedure to succeed,  the
approximation must be expressed in a well-behaved basis such as the discrete
orthogonal polynomial basis described in the preceding section.  With a
badly behaved basis, the singular values of $W$ may not decay as rapidly,
resulting in opportunities for degree reduction (and thus for singularity
reduction) being missed.

How should we choose the threshold $\eta$?  With noiseless input data, a
singular value of $W$ will be negligible if its size relative to the largest
singular value is on the order of the rounding error incurred during the
computation.  In this case, an appropriate choice for $\eta$ is a small power
of $10$ times the machine epsilon; for double-precision arithmetic, values such
as $\eta = 10^{-12}$ or $\eta = 10^{-14}$ work well.  If the input data are
noisy, the threshold should be increased so that any singular value below the
noise level is regarded as negligible.  For instance, if all but the $6$
leading digits of the data are noisy, then setting $\eta = 10^{-5}$ (a factor
of $10$ larger than the relative noise level of $10^{-6}$) may be appropriate.

\begin{algorithm2e}[htb!]
	\SetNoFillComment \DontPrintSemicolon \setcounter{AlgoLine}{0}
	\caption{Degree Reduction}\label{ALG:ReduceDegree}
	\SetKwInOut{Input}{Input}
	\SetKwInOut{Output}{Output}
	\Input{Vandermonde-like matrix $V_L$ computed with \cref{ALG:MVVandQR}, diagonal matrix $F$ of sample values, maximum numerator and denominator degrees $M$ and $N$, threshold $\eta$.}
	\Output{Reduced degrees $M$ and $N$.}
	\tcc{Reduce the denominator degree.}
	\While{true}{
		$Z \gets V_{M - 1}^*FV_N$\;
		$\sigma_{\mathrm{min}}, \sigma{_\mathrm{max}} \gets \text{smallest, largest singular values of $Z$}$\;
		\If{$\sigma_{\mathrm{min}} < \eta \sigma_{\mathrm{max}}$} {
			$M \gets M - 1$\;
		}\Else{
			break\;
		}
	}

	\tcc{Reduce the numerator degree.}
	\While{true}{
		$Z \gets V_{N - 1}^*F^{-1}V_M$\;
		$\sigma_{\mathrm{min}}, \sigma{_\mathrm{max}} \gets \text{smallest, largest singular values of $Z$}$\;
		\If{$\sigma_{\mathrm{min}} < \eta \sigma_{\mathrm{max}}$} {
			$N \gets N - 1$\;
		}\Else{
			break\;
		}
	}
\end{algorithm2e}

\begin{figure}[t!]
\centering
	\begin{subfigure}{0.45\textwidth}
		\centering
		\includegraphics[scale = 0.5]{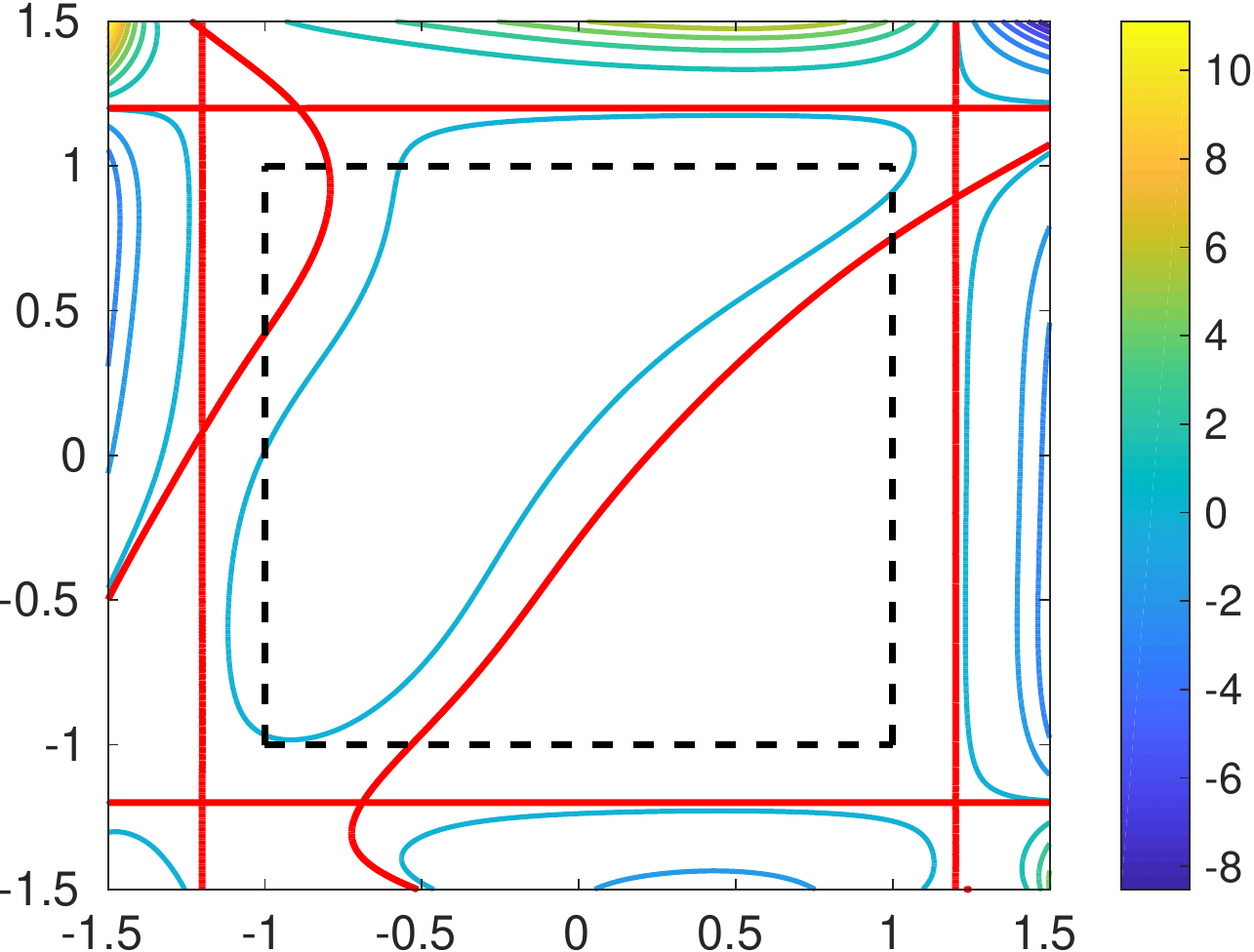}
		\caption{Without degree reduction.}
		\label{SFIG:DegRedExOff}
	\end{subfigure}
	~
	\begin{subfigure}{0.45\textwidth}
		\centering
		\includegraphics[scale = 0.5]{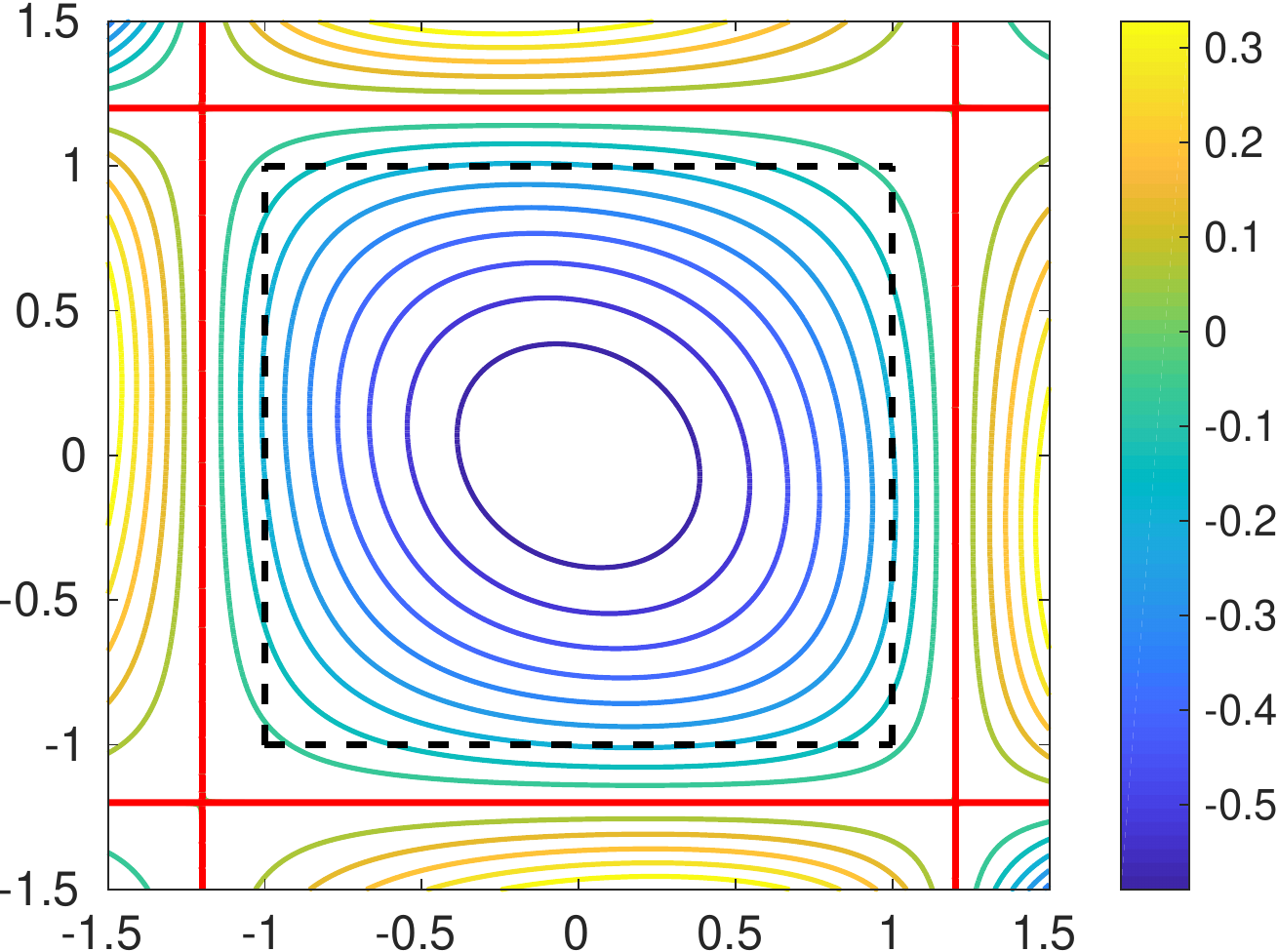}
		\caption{With degree reduction.}
		\label{SFIG:DegRedExOn}
	\end{subfigure}
\caption{Contour plots of the denominators of the rational approximations
computed in the example of \cref{sec:degred}.  Red lines denote zero-level curves (and hence curves of singularities present in the approximation).
The dashed black line outlines the unit square $[-1, 1] \times [-1, 1]$.}
\label{FIG:DegreeReductionExample}
\end{figure}

To illustrate the effectiveness of this general procedure, we consider
the problem of computing a rational approximation to the bivariate function
$f(x, y) = \exp(xy)/\bigl((x^2 - 1.44)(y^2 - 1.44)\bigr)$.  We sample this
function in 1,000 uniformly randomly distributed points in $[-1, 1] \times [-1,
1]$ and attempt to fit a rational approximation with numerator and denominator
degrees $M, N = 20$. \Cref{FIG:DegreeReductionExample} displays a
contour plot of the denominator of the computed approximation.  Without degree
reduction, we obtain the picture in \cref{SFIG:DegRedExOff}.  In addition
to the singularity curves at $x = \pm 1.2$ and $y = \pm 1.2$ that reflect the
true singularities of $f$, the approximation possesses a pair of spurious
singularity curves that wind their way through the middle of the square.
Applying \cref{ALG:ReduceDegree} with $\eta = 10^{-12}$ reduces the
numerator degree to $M = 12$ and the denominator degree to $N = 9$ and produces
a rational approximation with a denominator that generates the contour plot of
\cref{SFIG:DegRedExOn}.  The spurious singularity curves have
disappeared.

%% file: opti.tex
\section{Multidimensional Rational Approximation with Constraints}\label{SEC:opti}


The algorithm just described is simple and powerful; however, even with degree
reduction, it does not guarantee that the computed approximation is free of
singularities in the domain of interest.  In this section, we add constraints
to the rational approximation problem \cref{EQN:LSQ} that enforce this
requirement. We show that these constraints lead to a semi-infinite
optimization problem (see, e.g.,
\cite{stein2013bi,HettichKortanek:1993,ben2009robust}), which we solve using an
outer approximation approach due to Polyak~\citep{LEVITIN19661}. We are
motivated by a class of structural constraints that arise in HEP data analysis,
for which it is known that the underlying function has no poles in a certain
domain $D$ (but may have them outside of $D$), and we exploit this information
by enforcing the same condition for our rational approximation.


Formally, we can write the constraint that  ``$r(x)$ has no poles in $D$'' as the condition that
\[ q(x) \not = 0, \quad \forall \; x \in D .\]
However, this condition is not a convenient constraint to add to \cref{EQN:LSQ} because it describes an open set. Instead, we use the equivalent condition that ``$q(x)$ does not change sign in $D$,'' which can be written without loss of generality as
\begin{equation}\label{eq:SIPcon}
  q(x) \geq \tau>0, \quad \; \forall x \in D,
\end{equation}
where $\tau > 0$ is an arbitrary positive constant.  (We use $\tau=1$ in our experiments.) In most cases, $D$ will be a simple set such as bounded hyper-rectangle  $D=\prod_{i=1}^n[L_i,U_i]$. We then formulate the multivariate constrained rational approximation problem as the following constrained least-squares problem:
\begin{equation} \label{eq:rationalApprox}
  \mini_{p,q} \; \sum_{k=0}^{K-1} \left( p(\xk) - f_k q(\xk) \right)^2 \quad
  \st q(x) \geq \tau, \; \forall x \in D \quad \text{and} \quad p \in \sP_M^n, q \in \sP_N^n.
\end{equation}
This is a linear least-squares problem in the coefficients of the polynomials $p(x)$ and $q(x)$ with a linear semi-infinite constraint; see, for example, the monographs and surveys by \cite{stein2013bi,HettichKortanek:1993,ben2009robust}.

If $q(x)=b^T x + b_0$ and $D$ is affine, then we can use linear programming duality to replace the semi-infinite constraint by a set of equivalent finite-dimensional affine constraints; see, for example,\cite{ben2009robust}.  In general, however, this transformation does not exist unless we also assume that $q(x)$ is convex, which would add a semi-definite constraint in the quadratic case and more complex conic constraints in general. Hence, we will instead consider an outer approximation approach to solving \cref{eq:rationalApprox}.

\subsection{A Practical Algorithm for General Rational Approximation}\label{sec:rasip}

For general denominators, we apply a method due to Polyak. The algorithm maintains a finite set $U$ of points $\xk \in D$  at which the semi-infinite constraint is enforced. It then alternates between solving the finite-dimensional relaxation of \cref{eq:rationalApprox}, which at iteration $l$  is given by
\begin{equation} \label{eq:rationalApproxFD}
  \mini_{p,q} \; \sum_{k=0}^{K-1} \left( p(\xk) - f_k q(\xk) \right)^2 \quad \st q(\xk) \geq \tau, \; \forall k=0,\ldots, K-1+l,
\end{equation}
and an optimization problem to check \cref{eq:SIPcon}.
We let the solution of this problem be $p_l(x), q_l(x)$, and then solve the following minimization problem to global optimality to check whether \cref{eq:SIPcon} holds:
\begin{equation}\label{eq:globalMin} \mini_{x \in D} \; q_l(x).
\end{equation}
Either we  obtain a new point $\hat{x} \in D$ that violates $q_l(\hat{x}) \geq \tau$ or we show that $q_l(x)>0$ for all $x \in D$.
Formally, this procedure is defined in \cref{A:Polyak}.

\begin{algorithm2e}[htb!]
	\DontPrintSemicolon \setcounter{AlgoLine}{0}
	\caption{Alternating Algorithm for Pole-Free Rational Approximation.}\label{A:Polyak}
	\SetKwInOut{Input}{Input}
	\SetKwInOut{Output}{Output}
	\Input{$\{ x^{(0)}, \ldots, x^{(K-1)} \}$}
	\Output{Pole-free rational approximation $p_l(x)/q_l(x)$}
	Set $l \gets 0$, \textit{done} $\gets$ \textit{false}\; 
	\Repeat{done is false}{
		Let $p_l(x), q_l(x)$ be a solution of the relaxation \cref{eq:rationalApproxFD}.\;
		Let $\hat{x}$ be a (global) minimizer of \cref{eq:globalMin}.\;
		\uIf {$q_l(\hat{x}) \geq \tau$}{
			Set \textit{done} $\gets$ \textit{true}\;
		}
		\Else{
			Add a new point: $\{x^{(K+l)} : = \hat{x} \}$ and set $l \gets l+1$\;
		}
	}
\end{algorithm2e}

 
We note that we can stop the algorithm as soon as $q_l(\hat{x}) > 0$, which indicates that $q(x)$ has no poles in $D$. The final $p_l(x)/q_l(x)$ is the best (least-squares) interpolant that has no poles  in $D$. Unfortunately, the algorithm requires the global minimization of the polynomial $q_l(x)$ over $D$. We can either resort to multistarts (multiple local optimizations starting from  different  points), or compute an underestimator of $q_l(x)$ on $D$ using the reformulation-linearization-technique  of \cite{sherali.adams:98}. We discuss a practical way of solving the global optimization problem using multistarts in \cref{app:globaloptprob}.

%% file: numerics.tex
\section{Numerical Experiments}\label{SEC:Numerics}

In this section, we compare the approximation quality and computation times of the rational and polynomial approximation approaches. We also study the effects of using different strategies for sampling the interpolation points from the domain and the effects of constraints on the rational approximation.
\subsection{Experimental Setup}
Our numerical experiments are conducted on a server with 64 Intel Xeon Gold CPU running at 2.30 GHz. There are two threads per core, but each approximation is run on a single thread. The operating system is Linux Ubuntu 16.04. Additionally, the server is equipped with 1.5TB DDR4 2666 MHz of memory. The code is written in Python v3.7.2 where the optimization functions and constraints are compiled with the Numba JIT compiler v0.42.

The experiments are conducted on fast-to-compute analytic test problems whose functional forms are summarized in \cref{T:TPs}. The use of these analytic test problems enables us to  assess the performance of our algorithms efficiently. We show detailed results for five typical test functions that span the range of the functions of interest in this section, and we summarize the remaining results for the other functions, which are included in the electronic supplement \cref{app:supppoleserrors,app:supptimesiters}. 
In the following we show the results for \cref{fn:f4} whose approximation using Taylor series expansions is a polynomial function; \cref{fn:f8}, which is a rational function; \cref{fn:f17}, which is used to describe a resonant particle of mass $M$ and width $\Gamma$ as a function of the particle's energy $E$ in high-energy physics~\citep{PhysRev.49.519,Bohm:2004zi}; and \cref{fn:f18,,fn:f19}, whose approximation using Taylor series expansions is a rational function. Note that the domain of \cref{fn:f18} is close to the true pole.

We sample the interpolation points $\{x^{(0)},\dots,x^{(K-1)}\}$ using sparse grids (SGs)~\citep{Barthelmann2000} and Latin hypercube sampling (LHS)~\citep{MckayThreeMethods}, and we propose a new hybrid strategy called decoupled Latin hypercube design (d-LHD) where the interpolation points are sampled on the faces and inside the domain. A plot of the interpolation points sampled by sing the three different strategies is shown in \cref{fig:samplingstrategies}.
The approximation results change for interpolation points that are sampled by using the LHS and d-LHD strategies  because they have randomness. To account for these changes, each experiment is repeated five times  using  different random number seeds, and we  report the mean and other statistics of the performance metrics for these strategies.
We also experimented with uniform randomly sampled points, but we found the results to be inferior and therefore do not include them here.

Each functional value $f_k$ is obtained by evaluating \textit{f} at $x^{(k)}$.  
The number of interpolation points \textit{K} is set as twice the sum of the number of degrees of freedom of the polynomials in each approximation given by $\alpha(M)$ + $\alpha(N)$ for numerator of degree $M$ and denominator of degree $N$.
We consider both noise-free and noisy data in the experiments. For noisy data, each functional value $f_k$ is multiplied by a fraction $\epsilon$ of the random value $\phi^{(k)}$ sampled from a standard normal distribution $\mathcal{N}(0,1)$ as follows:
\begin{equation}\label{eq:noise}
f_k = f_k \left(1+\epsilon\phi^{(k)}\right),\quad\forall k=0,\dots,K-1.
\end{equation} 
The approximation $r(x)$ is computed in four  ways: (1) $p(x)$ is the polynomial approximation that is computed by a \textit{NumPy} implementation of finding the linear least-squares solution using singular value decomposition within the driver routine \textit{DGELSD}~\citep{laug}, (2) $r_1(x)$ is the rational approximation using \cref{ALG:MVVandQR} without degree reduction, (3) $r_2(x)$ is the rational approximation using \cref{ALG:MVVandQR} with the degree reduction described in \cref{ALG:ReduceDegree}, and (4) $r_3(x)$ is the rational approximation using \cref{A:Polyak}.

To assess the quality of our approximations,  we use a second set of testing points $\{x^{(K)},\dots,x^{(L-1)}\}$ on the faces and inside of the domain, and we compute their function values $\{f_K,\dots,f_{L-1}\}$. No noise is added to the testing data. 
We use  the $l_2$-norm error as a test metric  to compare the quality of the approximation $r(x)$:
\begin{equation}\label{eq:testError}
\Delta_r = \left|\left|r - f\right|\right|_{D,2} = \left(\int_{D} \left(r(x)-f(x)\right)^2 dx\right)^{1/2} \approx \left\{\sum\limits_{k = K}^{L-1} \left[r\left(x^{(k)}\right) - f_k\right] ^2\right\}^{1/2}.
\end{equation}
We consider a solution to be better if it has smaller $\Delta_r$. 
We assume that the degrees of the numerator and the denominator polynomials of the approximations each are 5 . This choice allows us to approximate test functions in which the polynomials are up to degree 4. Choosing the optimal degree of polynomials is a question that is beyond the scope of this paper.

\subsection{Effects of Interpolation Point Selection Method}

\begin{figure}[htbp]
	\centering
	\begin{subfigure}[b]{0.3\textwidth}
		\includegraphics[width=1\textwidth]{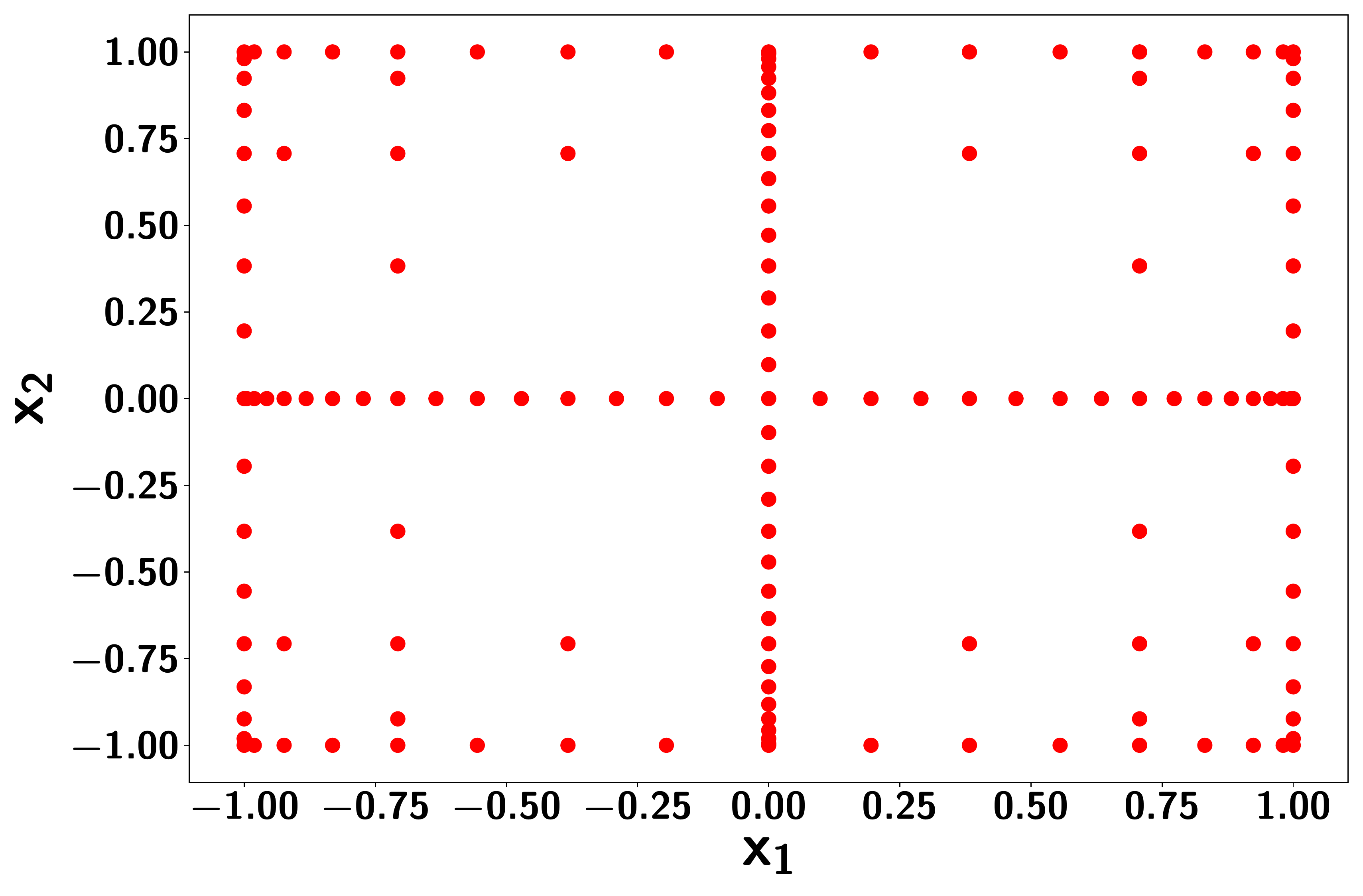}
		\caption{2D Sparse gid of level 5 with 145 points.}%
		\label{fig:sparsegrid}
	\end{subfigure}
	~
	\begin{subfigure}[b]{0.3\textwidth}  
		\includegraphics[width=1\textwidth]{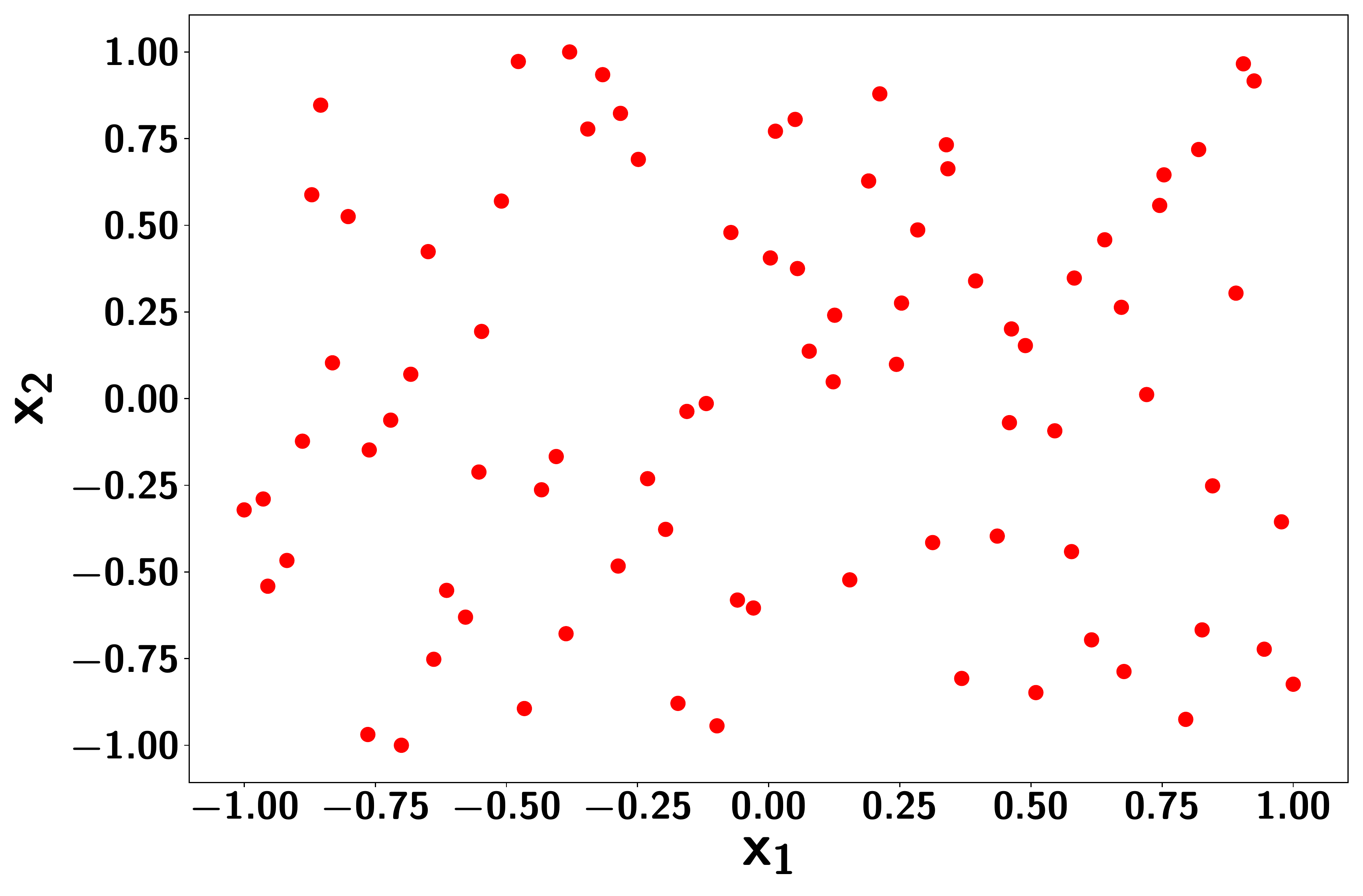}
		\caption{2D Latin hypercube sample with 84 points.}
		\label{fig:lhssample}
	\end{subfigure}
	~
	\begin{subfigure}[b]{0.3\textwidth}  
		\includegraphics[width=1\textwidth]{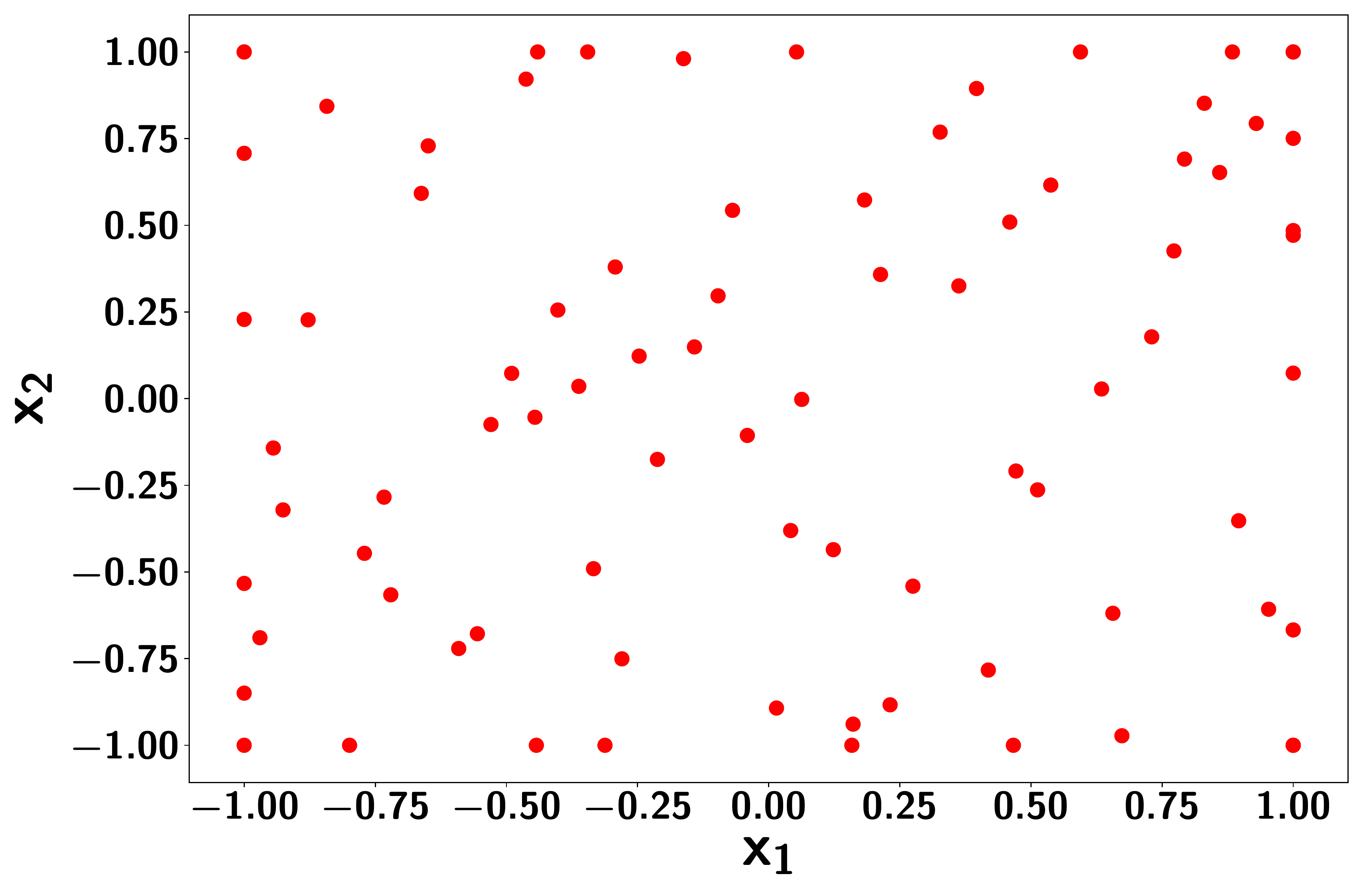}
		\caption{2D decoupled Latin hypercube Design with 84 points.}
		\label{fig:d-lhdsample}
	\end{subfigure}
	\caption{Location of interpolation points using different sampling strategies for a rational approximation of $M=5$,  $N=5$ and $\alpha(M) + \alpha(N) = 42$.}
	\label{fig:samplingstrategies}
\end{figure}

\begin{figure}[htb]
  	\begin{center}
  		\includegraphics[width=0.7\textwidth]{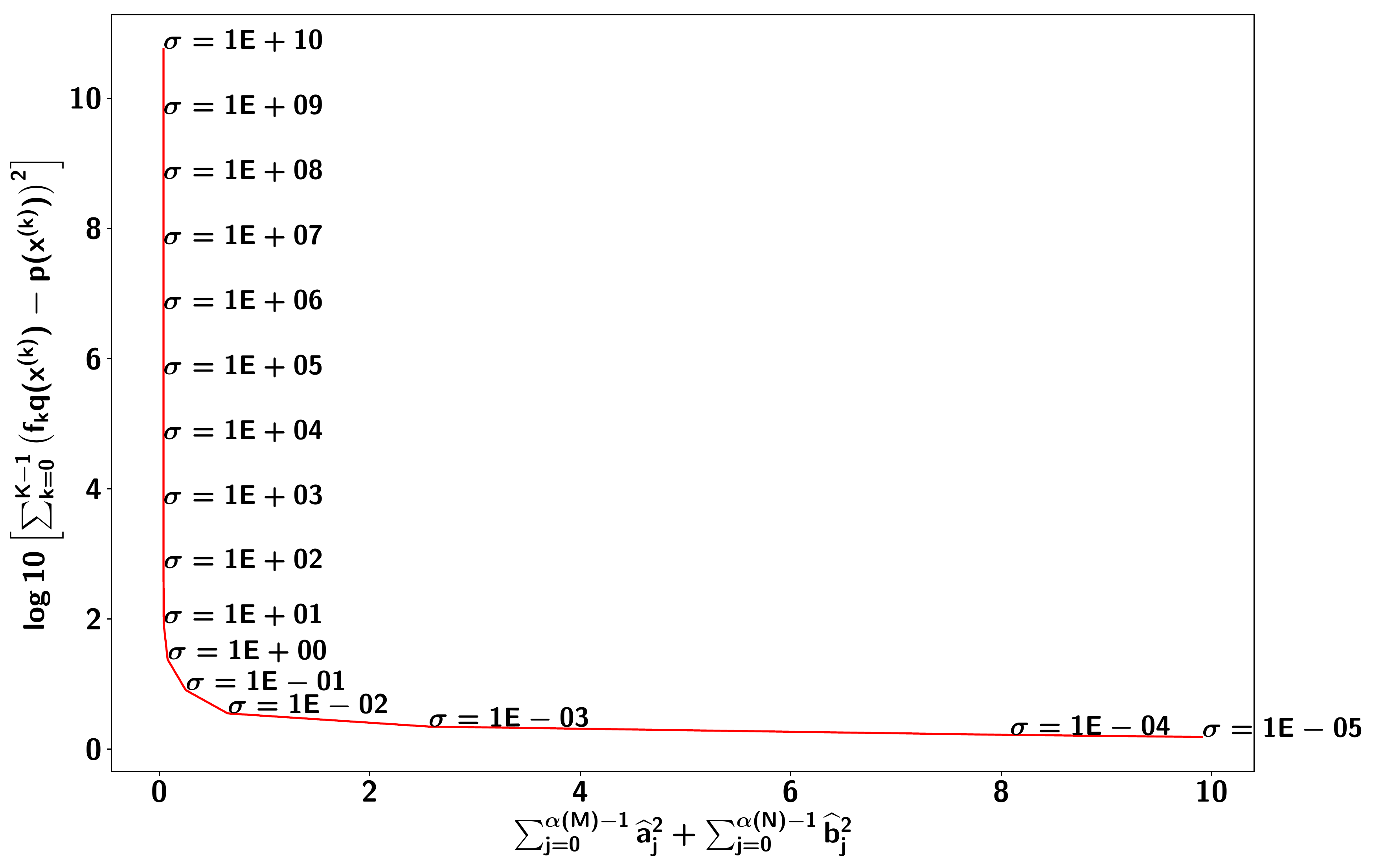}
  		\caption{L-curve method to choose $\sigma$. \Cref{fn:f18} interpolation data is sampled by using SG. The approximations are performed by using \cref{A:Polyak} with \cref{eq:rationalApproxFD2} instead of \cref{eq:rationalApproxFD} for different values of $\sigma$. The degrees of the numerator and denominator polynomials are $M=5$ and $N=5$, respectively. The corner of the L is found at $\sigma=10^{-1}$.}
  		\label{fig:choosingsigma}
  	\end{center}
  \end{figure}

In this section we discuss the choice  of the interpolation points using SGs~\citep{Barthelmann2000}, LHS~\citep{MckayThreeMethods}  and d-LHD strategies. SGs, first proposed by Smolyak, are sparse tensor product spaces.
We also experimented with a uniform random set of points but observed uncompetitive results. With SGs, the grid points are obtained by combining, up to a certain level, the tensor product grid corresponding to the total degree multi-index set.
Here, the SG level is chosen such that the number of points in the grid is at least twice the total degrees of freedom of the polynomials in each approximation.
\cref{fig:sparsegrid} shows a 2D SG. We observe that many points of the SG are collinear, violating the linear independence assumption from \cref{SEC:LinAlg}. Hence, for the chosen SG levels, the Hessian matrix of the fitting problem in \cref{eq:rationalApproxFD} is singular. Because of the null space, there are multiple minimizers that result in unbounded values for $p$ and $q$.
We overcome this issue by adding a regularization term with weight $\sigma>0$ to \cref{eq:rationalApproxFD} and bounding the eigenvalues $>\sigma$. The updated fitting problem in iteration $l$ of \cref{A:Polyak} is 
  \begin{equation} \label{eq:rationalApproxFD2}
  \begin{split}
	&\mini_{p,q} \; \sum_{k=0}^{K-1} \left( f_k q(\xk) - p(\xk) \right)^2 + \sigma\left(\sum_{j=0}^{\alpha(M)-1} \widehat{a}_j^2 + \sum_{j=0}^{\alpha(N)-1} \widehat{b}_j^2\right)\\ &\st q(\xk) \geq 1, \; \forall k=0,\ldots, K-1+l,
 \end{split}
  \end{equation}
  where $\widehat{a}$ and $\widehat{b}$ are the coefficients of the monomial basis expansion of $p(x)$ and $q(x)$, respectively.
  To choose $\sigma$, we ran \cref{A:Polyak} with \cref{eq:rationalApproxFD2} to approximate the data sampled with SG. We found $\sigma$ to be in the vicinity of $10^{-1}$ for all test functions using the L-curve method. An example plot of the L-curve for \Cref{fn:f18} interpolation data is shown in \cref{fig:choosingsigma}. 

	 \begin{figure}[htb]
		\centering
		\begin{minipage}{.62\textwidth}
			\centering
			\includegraphics[width=1\linewidth]{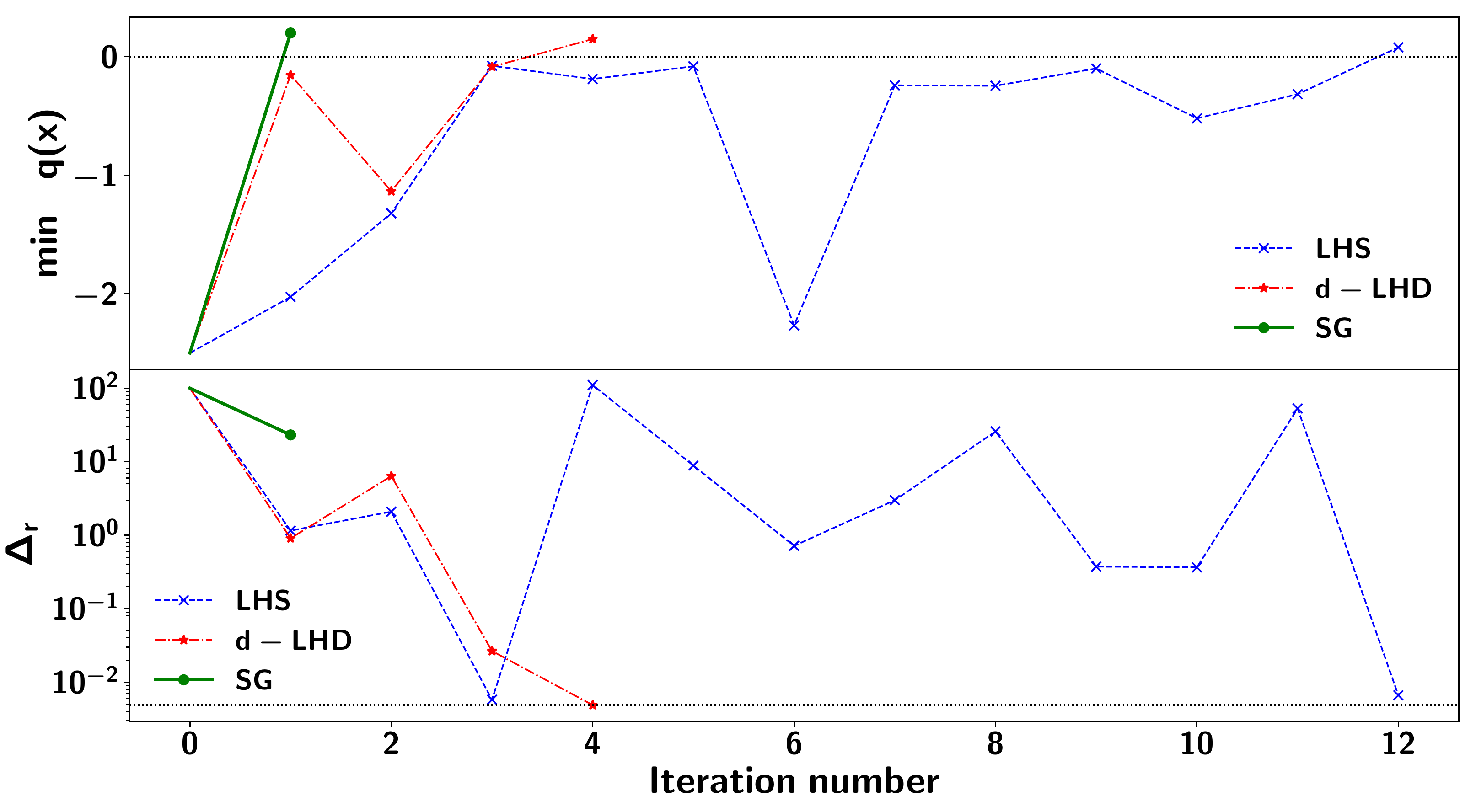}
			\captionof{figure}{Minimizers of \cref{eq:globalMin} (top) and testing error (bottom) per iteration of a run of \cref{A:Polyak} on \cref{fn:f17} interpolation data sampled by using all three strategies. In the top plot, \cref{A:Polyak} stops iterating when the minimizer is at least 0 (shown with a horizontal dotted line). In the bottom plot, the approximation obtained by using d-LHD sampled data has the lowest testing error (shown with a horizontal dotted line). The numbering of the iterations of \cref{A:Polyak} starts at 1. In both plots, the lines for all strategies start from the same extreme point (at iteration 0).}
			\label{fig:minimizererror}
		\end{minipage}
		\hfill
		\begin{minipage}{.355\textwidth}
			\centering
			\includegraphics[width=\textwidth]{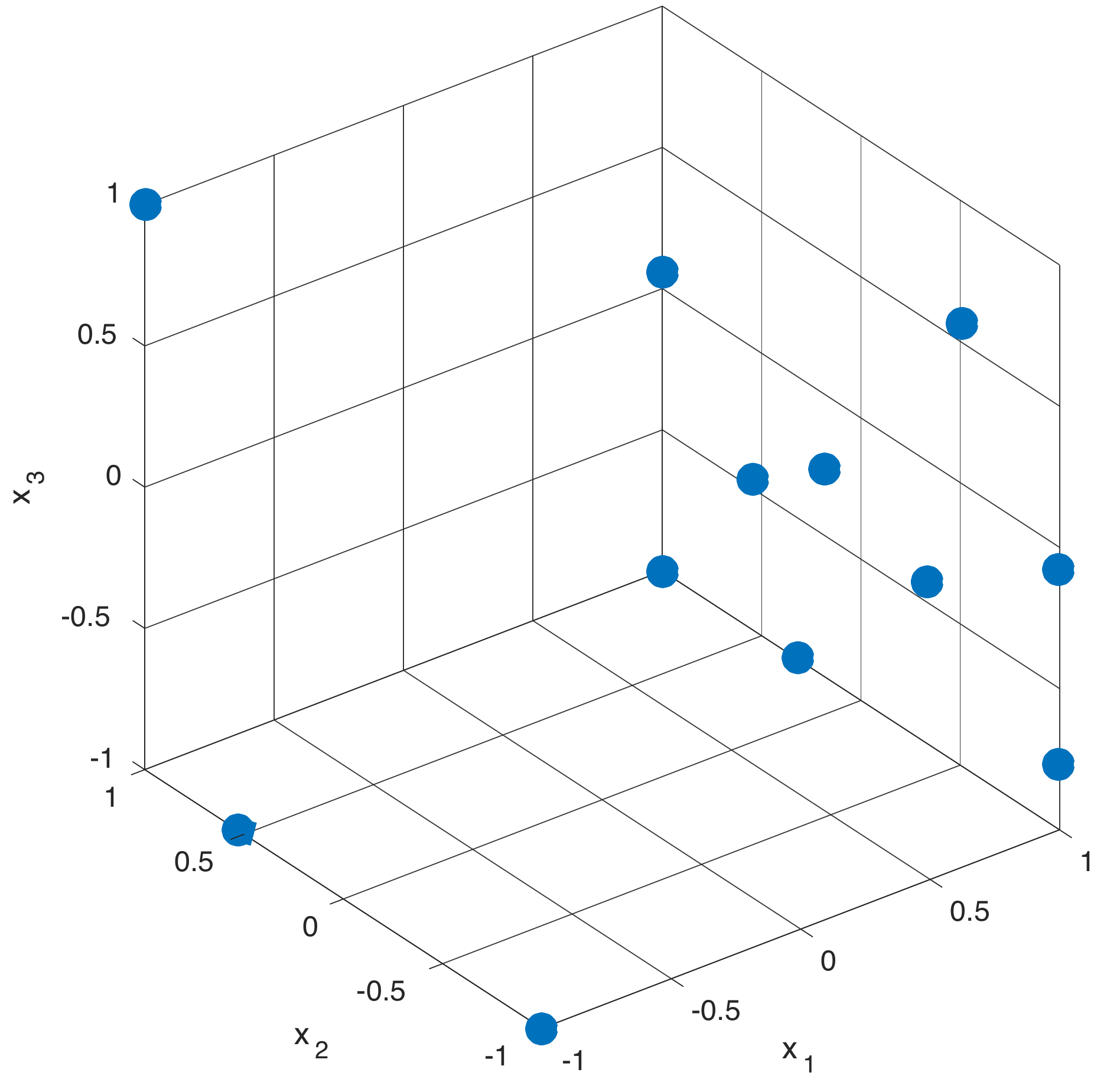}
			\caption{Minimizer of \cref{eq:globalMin} from all iterations of \cref{A:Polyak} on \cref{fn:f17} interpolation data. The data are sampled by using LHS. All minimizers lie on the face of the cubic domain.}
			\label{fig:minimizersonface}
		\end{minipage}%
	\end{figure}

	 When using the SG interpolation points, we observe that \cref{A:Polyak} takes only one iteration to converge to a pole-free rational approximation. This is shown in the top plot of \cref{fig:minimizererror}. 
	However, penalizing the coefficients of the monomial basis expansions of $p(x)$ and $q(x)$ results in a high testing error thereby deteriorating the approximation quality, as shown in the bottom plot of \cref{fig:minimizererror}, which is undesirable. 
	
\begin{figure}[htb!]
		\centering
		\includegraphics[width=.7\textwidth]{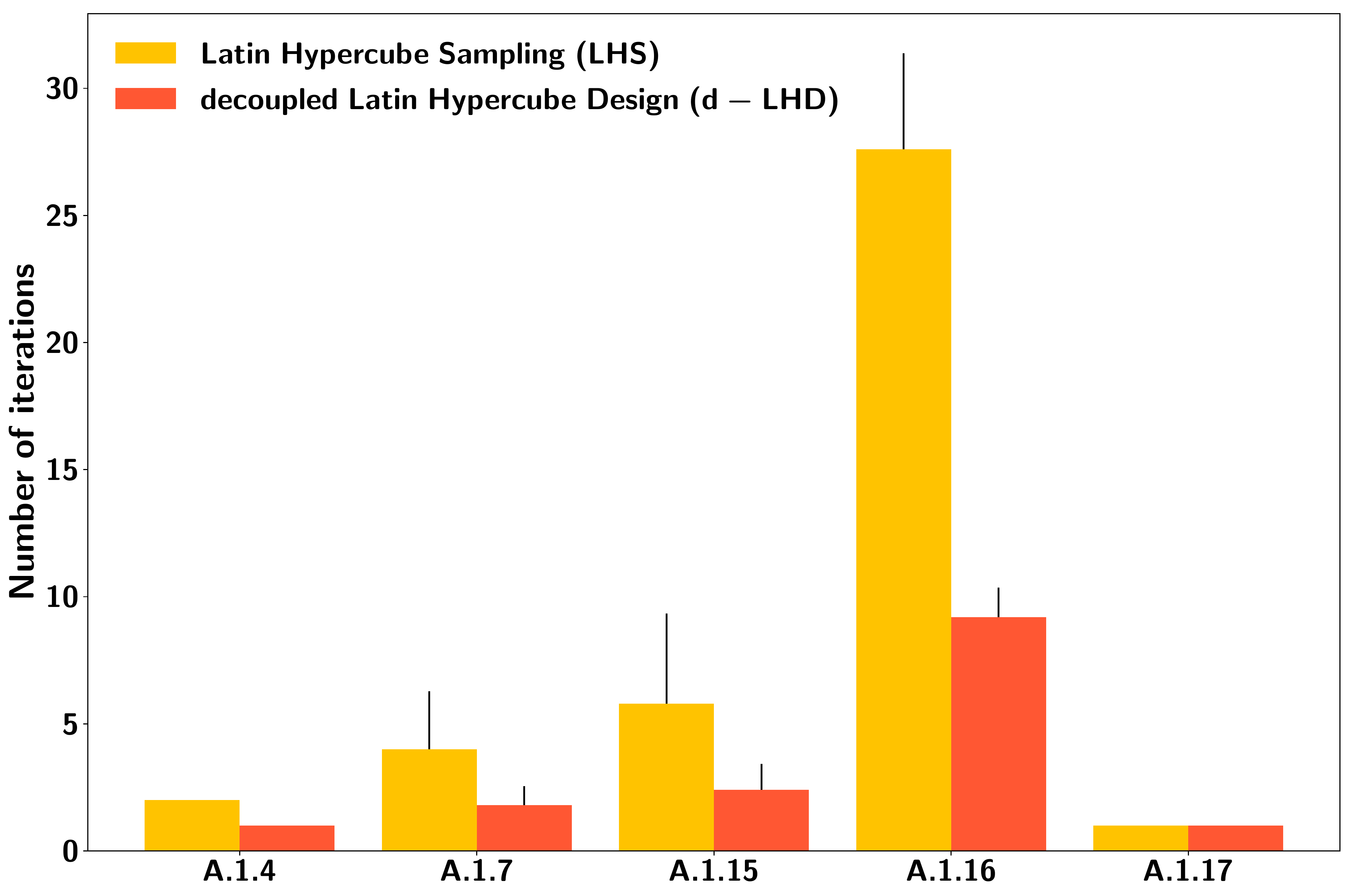}
		\caption{
			Number of iterations performed by \cref{A:Polyak} over noise-free data sampled by using LHS and d-LHD strategies. The standard deviation  is shown as a black vertical line. When the data are sampled by using SG, the number of iterations is 1 for all functions. As shown in the bottom plot of \cref{fig:minimizererror}, however,the testing error of the approximation obtained by using SG is much higher than those obtained using LHS and d-LHD. 
		}

		\label{fig:iteration}
	\end{figure}

	When using  LHS, each independent dimension is sampled by using an even sampling method, and then these samples are randomly combined to obtain the sample data. \cref{fig:lhssample} shows one such 2D LHS sample. The advantage of LHS is that the interpolation points are not collinear. This results in a Hessian matrix of full rank for the fitting problem in \cref{eq:rationalApproxFD} and hence does not require the regularization term to be added. However, the number of iterations of \cref{A:Polyak} over noise-free LHS data is on average five times the number of iterations over noise-free SG data (see \cref{fig:iteration}).
	
	We observe that almost all the minimizers of $q(x)$ found in each iteration of \cref{A:Polyak} lie on the faces of the domain as shown in \cref{fig:minimizersonface}. 
	Hence, when SG places a number of points on the face of the domain, the number of spurious poles is minimized, which on average requires fewer iterations of \cref{A:Polyak}. 	
	So ideally, we want to use a sampling strategy that covers the faces and the inside of the domain evenly such that the points are not collinear, thereby combining the best features of SG and LHS. 
	
	The authors of \cite{LEATHERMAN2017346} proposed  the maximin augmented nested Latin hypercube design sampling strategy  to maximize prediction accuracy.  In this strategy, the samples are constructed by augmenting nested LHDs with additional parameters using a modified smart swap algorithm such that the final design satisfies the maxmin property. 
	However, the required properties of the samples to satisfy our goal are simpler, and we therefore use our  decoupled Latin hypercube design. 
	We construct nested LHDs over all the $2n$ facets of the domain with dimension $n-1$; in other words, one of the dimensions in each face's sample is fixed. Because these samples are LHDs, the points are not collinear.
	In order to cover the inside of the domain, an augmented LHD is obtained inside the  $n$-dimensional hyper-rectangle. 
	These two steps are independent. Even though the samples on each face and inside the domain satisfy the maxmin property, we do not require that the final design satisfy the maxmin property.
	We call this sampling strategy decoupled Latin Hypercube Design (d-LHD).
	
	In d-LHD, the number of points sampled is still twice the degrees of freedom of the polynomials in each approximation, namely, $K = 2(\alpha_n(M) + \alpha_n(N))$.
	On the $2n$ faces, there is a $n-1$ dimensional rational function;  hence, the number of points sampled along each face of the domain is given as
	\begin{equation}
	K^{(fc)} = \frac{2(\alpha_{n-1}(M) + \alpha_{n-1}(N))}{2n}=\frac{(\alpha_{n-1}(M) + \alpha_{n-1}(N))}{n}, 
	\end{equation}
	and the number of points sampled inside the domain is given as\
	\begin{equation}
	K^{(in)} = K - 2n\cdot K^{(fc)}.
	\end{equation}
	Thus, the d-LHD  samples points on each face as well as the inside of the domain as  illustrated in \cref{fig:d-lhdsample}. 
	
	\begin{table}[htb!]
		\centering
		\caption{
			Number of iterations over all test functions in \cref{T:TPs} performed by \cref{A:Polyak} over noise-free data sampled using LHS and d-LHD strategies. 
			When the data is sampled by using SG, the number of iterations is 1 for all functions. 
			However, the testing error of the approximation obtained by using SG for all functions is much higher than those obtained by using LHS and d-LHD (see electronic supplement \cref{app:supppoleserrors}). 
			For almost all functions, \cref{A:Polyak} takes fewer iterations over data sampled using d-LHD, which is evident from the median and the geometric mean.  However, the average number of iterations over \cref{fn:f20} data sampled by using d-LHD is 80.4 whereas that over data sampled by using LHS is 26.4. This outlier makes the arithmetic mean look better for the LHS strategy.
			The corresponding results for each function can be found in the electronic supplement \cref{app:supptimesiters}.
		} \label{tab:iterationsummary}
		\begin{tabular}{|c|c|c|} 
			\cline{1-3}
			\multicolumn{1}{|c|}{Statistic}&\multicolumn{1}{|c|}{LHS}&\multicolumn{1}{|c|}{d-LHD}
			\\\hline\hline
		Arithmetic Mean&4.17&5.50
		\\\hline
		Geometric Mean&1.86&1.51
		\\\hline
		Median&1.20&1.00
		\\\hline
		Range&26.60&79.40
		\\\hline
		\end{tabular}
	\end{table}
	
	\cref{fig:iteration} compares the number of iterations performed by \cref{A:Polyak} when the function domains are sampled with LHS and d-LHD, respectively. \cref{tab:iterationsummary} shows statistics for the number of iterations performed over all test functions in \cref{T:TPs}. Fitting the approximation to data sampled using SG  takes only one iteration, but it causes a higher  testing error compared with the other two strategies (see \cref{fig:minimizererror}). Additionally, when using d-LHD to sample points, the number of iterations of \cref{A:Polyak} is almost always lower compared with LHS (see  \cref{fig:iteration}) without  compromising  the approximation quality (see the bottom plot of \cref{fig:minimizererror}).
	Hence, in the remainder of this section, we  present results for the approximations performed with data sampled by d-LHD. The results corresponding to the other sampling strategies can be found in the electronic supplement \cref{app:supppoleserrors,app:supptimesiters}.

\subsection {Comparison of Approximation Quality}

	In this section, we evaluate the ability of constraints in \cref{A:Polyak} to remove spurious poles and compare the quality of our rational approximations. More specifically, we examine the number of spurious poles detected and its effect on the testing error in the three rational approximation approaches.
	Then, we compare the quality of the rational approximations with the polynomial approximation by comparing their testing errors.
	
	\subsubsection{Ability to Remove Spurious Poles}
	
	In this section, we compare the number of spurious poles found in the three rational approximation approaches, since we are interested in separating the error due to these poles from the actual approximation error. 
	Detecting these poles is difficult, however, because multivariate rational approximations are more complicated than univariate ones. They may have unaccountably many singular points, and these singularities are typically never isolated. Hence, to perform this comparison, we find the testing points near poles or polelike points that have large function deviations. 
	As shown in \cref{fig:minimizersonface}, the minimizers of $q_l(x)$ in each iteration of \cref{A:Polyak} tend to be on the boundary of the domain.
	Hence, we choose the testing points randomly on the faces of the domain in addition to randomly chosen points inside of the domain. 
	For these testing points, we define $W_{r,t}$ as the index set of points whose absolute approximated value is much larger than the corresponding absolute value of the function indicating a possible spurious pole. More formally, we define  $W_{r,t}$ as:
\begin{equation}\label{eq:noofpoles}
W_{r,t} =  W^{(fc)}_{r,t} \cup W^{(in)}_{r,t} =
\left\{j | \frac{|r(x^{(j)})|}{\max\left(1,|f^{(fc)}_{max}|\right)} > t, \right\} \cup \left\{k | \frac{|r(x^{(k)})|}{\max\left(1,|f^{(in)}_{max}|\right)} > t, \right\},
\end{equation}
where $x^{(j)}$ and $x^{(k)}$ are testing points on the face and inside of the domain respectively, $j \in I^{(fc)},k \in I^{(in)}, I^{(fc)} \cup I^{(in)} = \{K,\dots,L-1\},  I^{(fc)} \cap I^{(in)} =   \emptyset $, $f_{max}^{(fc)} = \max |f_j|$, $f_{max}^{(in)} = \max |f_k|$, and $t$ is a large threshold.

	\begin{figure}[htb!]
		\centering
		\begin{subfigure}[b]{0.33\textwidth}
			\centering
			\includegraphics[width=1\textwidth]{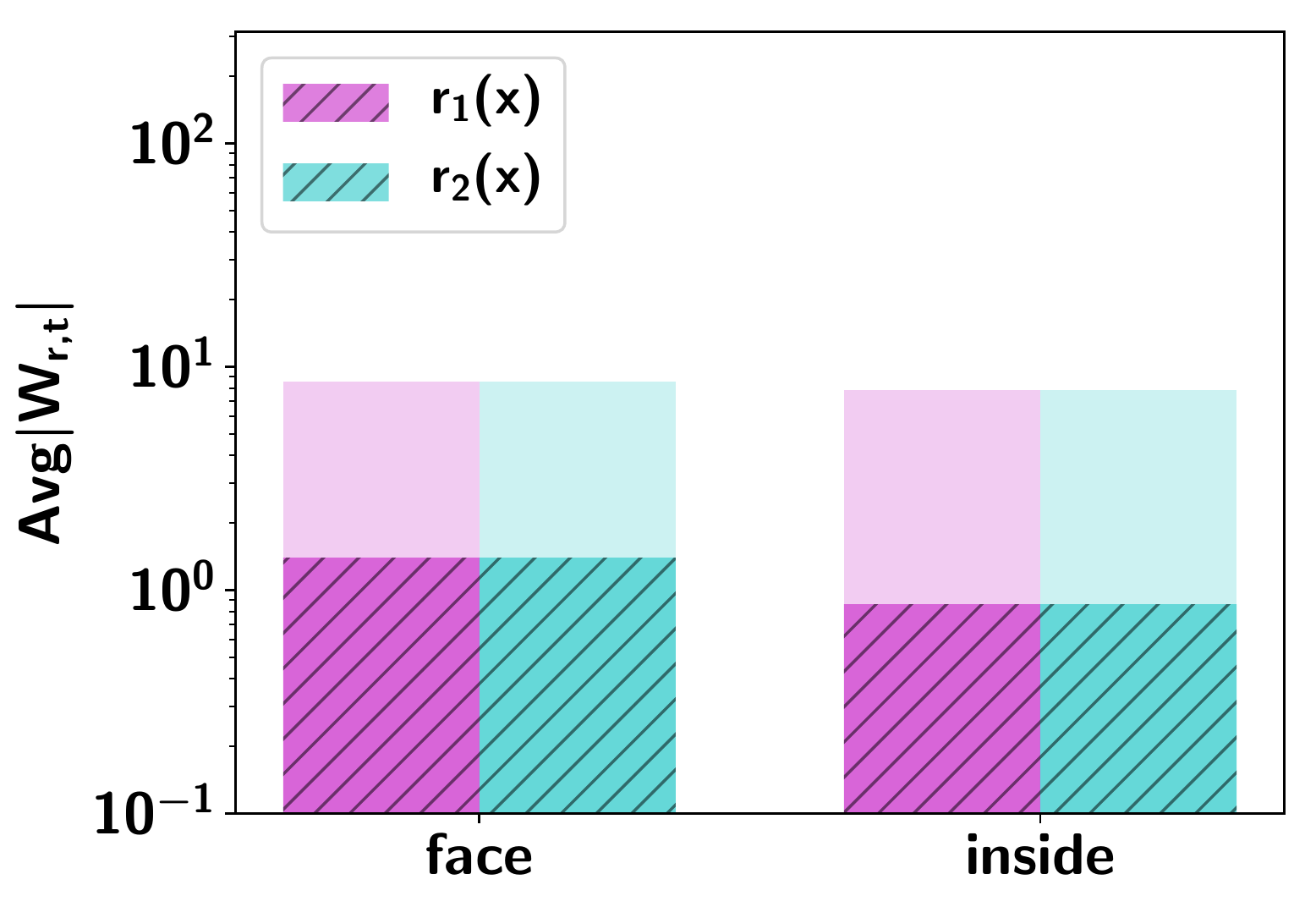}
			\caption%
			{{\small  $\epsilon=0$}}    
			\label{fig:polenonoise}
		\end{subfigure}
		\hfill
		\begin{subfigure}[b]{0.33\textwidth}  
			\centering 
			\includegraphics[width=1\textwidth]{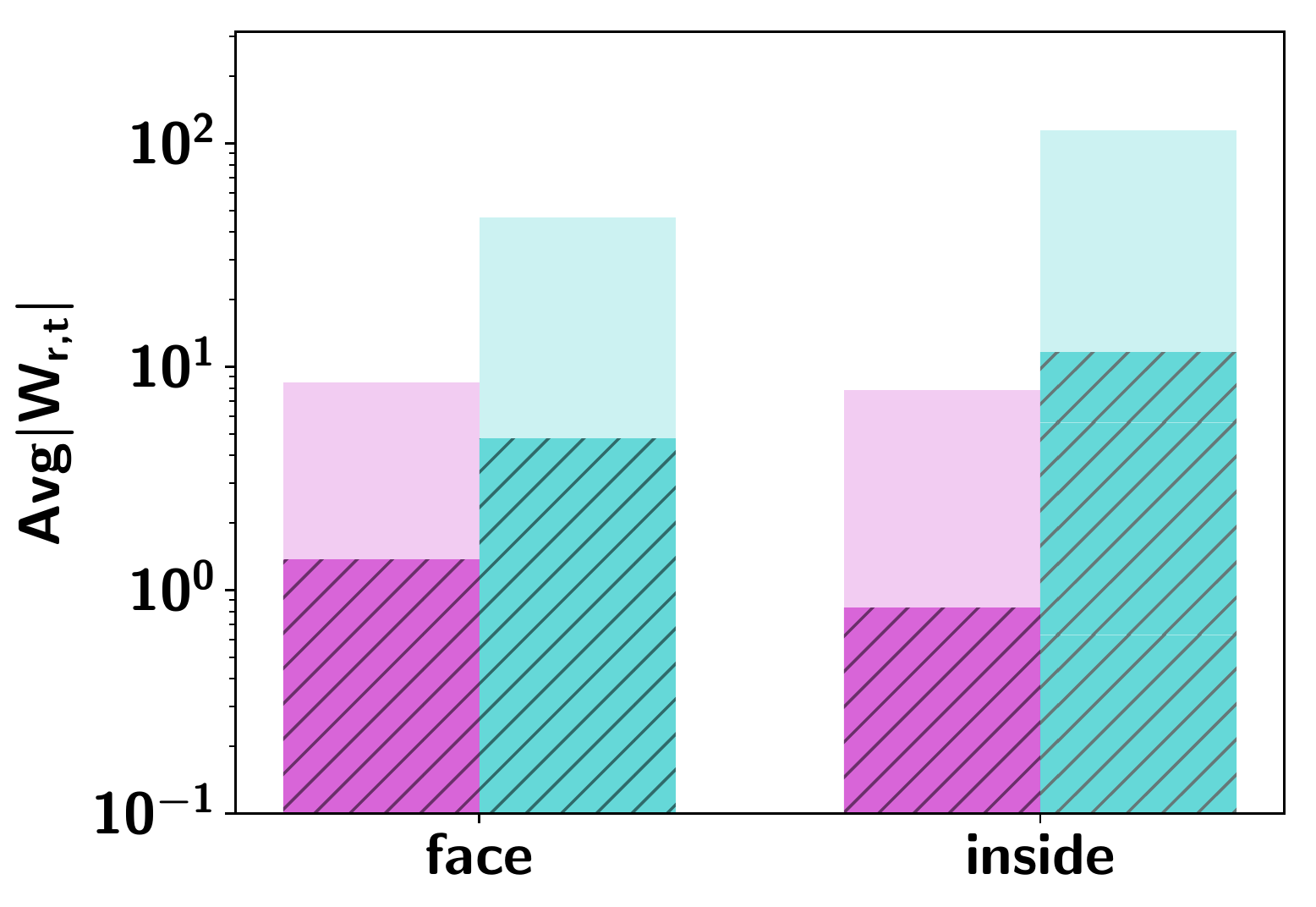}
			\caption[]%
			{{\small $\epsilon = 10^{-6}$}}    
			\label{fig:pole10-6noise}
		\end{subfigure}
		\hfill
		\begin{subfigure}[b]{0.33\textwidth}  
			\centering 
			\includegraphics[width=1\textwidth]{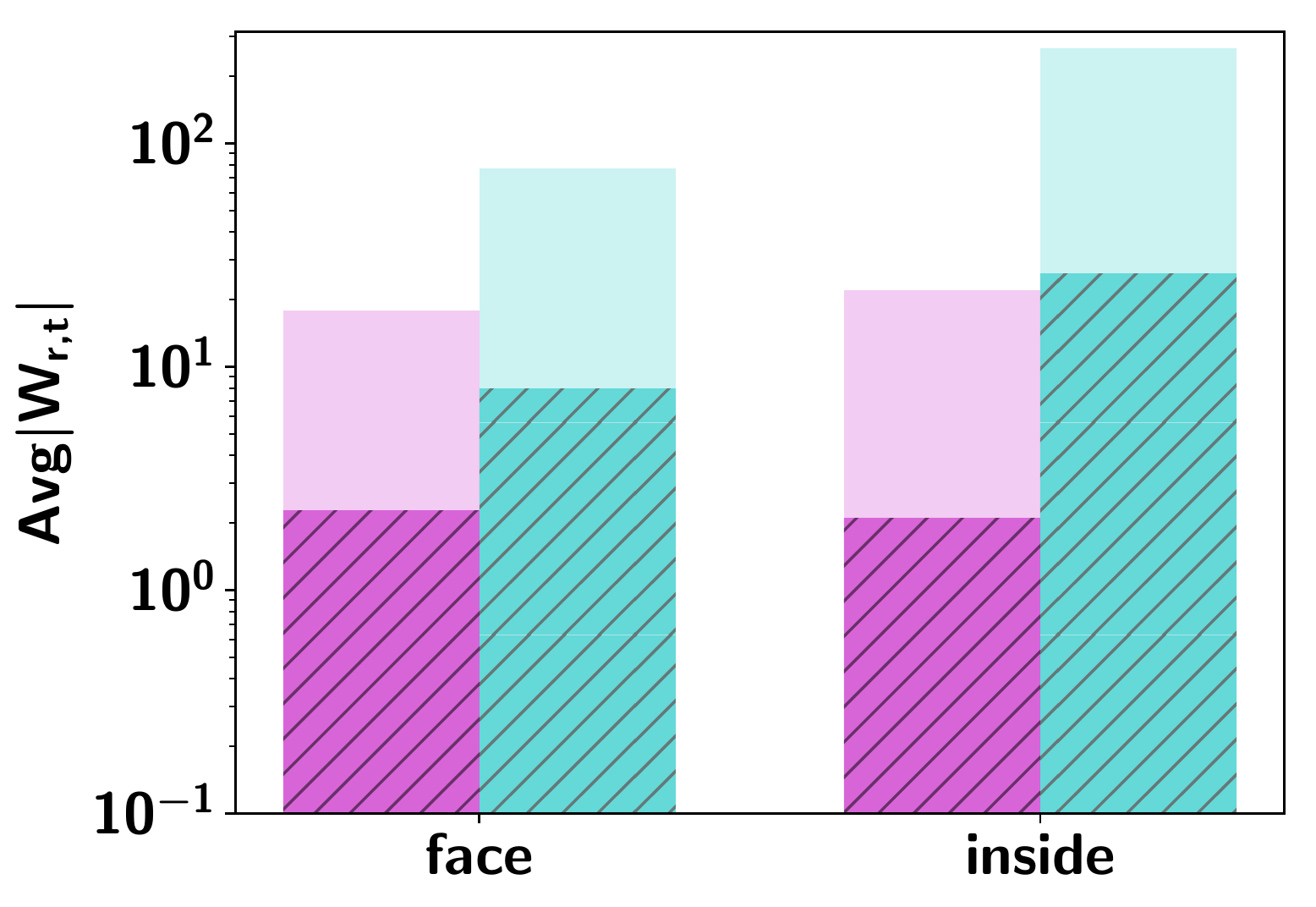}
			\caption[]%
			{{\small $\epsilon=10^{-2}$}}    
			\label{fig:pole10-2noise}
		\end{subfigure}
		\caption[]
		{
			Comparison of the average number of polelike points found over all functions in \cref{T:TPs} for different relative noise levels $\epsilon$.
			The data is sampled with d-LHD.
			Each bar represents the average number of polelike points found when $t\ge 10^2$.
			The average number of polelike points found when $t\ge 10^3$ is shown as a hatched bar.
			The average number of polelike points found when $10^2 \le t <10^3$ is shown as a faded bar.
			The number of polelike points found in $r_3(x)$ is 0 for all noise levels.
		}
%
%
		\label{fig:poleplot}
	\end{figure}

	\cref{fig:poleplot} shows the average number of polelike points found over all functions in \cref{T:TPs} when the interpolation data for these functions was sampled by using d-LHD. The number of polelike points per function in \cref{T:TPs} for d-LHD, SG, and LHS-based approximations is given in the electronic supplement \cref{app:supppoleserrors}. 
	We observe polelike points in $r_1(x)$ and $r_2(x)$ for noise-free and noisy interpolation data. 
	In contrast, the approximation, $r_3(x)$ does not have these polelike points. As discussed in \cref{sec:rasip}, this is due to the iterative removal of poles by \cref{A:Polyak} by design, thereby giving a pole-free $r_3(x)$.
	When no noise is added to the interpolation data, that is, when $\epsilon=0$, the number of pole-like points found on the faces of the domain is larger than inside of the domain in $r_1(x)$ and $r_2(x)$. 
	This difference is more prominent when the interpolation data are sampled by using LHS. 
	The number of polelike points found in $r_1(x)$ and $r_2(x)$ on the face is 24\% higher than those found on the inside whereas this difference is only 8.5\% when the interpolation data are sampled by using d-LHD.
	The reason is that LHS samples fewer interpolation points on the faces of the domain, 
	causing the LHS-based approximations to be less accurate on the boundary of the domain, especially when the  function domain is in close proximity to the true poles.
	Also, this result is consistent with  our earlier observation that the minimizers found in each iteration of \cref{A:Polyak} tend to lie on the face of the domain. 
	
	In the presence of noise, that is, when $\epsilon\neq 0$ the number of polelike points found in $r_1(x)$ and $r_2(x)$ increase inside as well as on the faces of the domain.
	As above, $r_3(x)$ does not suffer from spurious poles.
	The number of polelike points found for $r_2(x)$ is much higher than for $r_1(x)$. 
	As discussed before, the reason is that multivariate rational approximations are more complicated than univariate ones. Their singularities are never isolated; and even if we eliminate unnecessary degrees of freedom, they will still, in general, have unaccountably many singular points. 
	This problem is only compounded when the input data are noisy. Thus, we cannot hope that the degree-reduction approach will work as well in the multivariate case as it does in the univariate case.
	On the other hand, the optimization approach by design eliminates poles in $D$.

	\subsubsection{Comparison of the Testing Error}
	
	\begin{figure}[htb!]
		\centering
		\begin{subfigure}[b]{0.33\textwidth}
			\centering
			\includegraphics[width=1\textwidth]{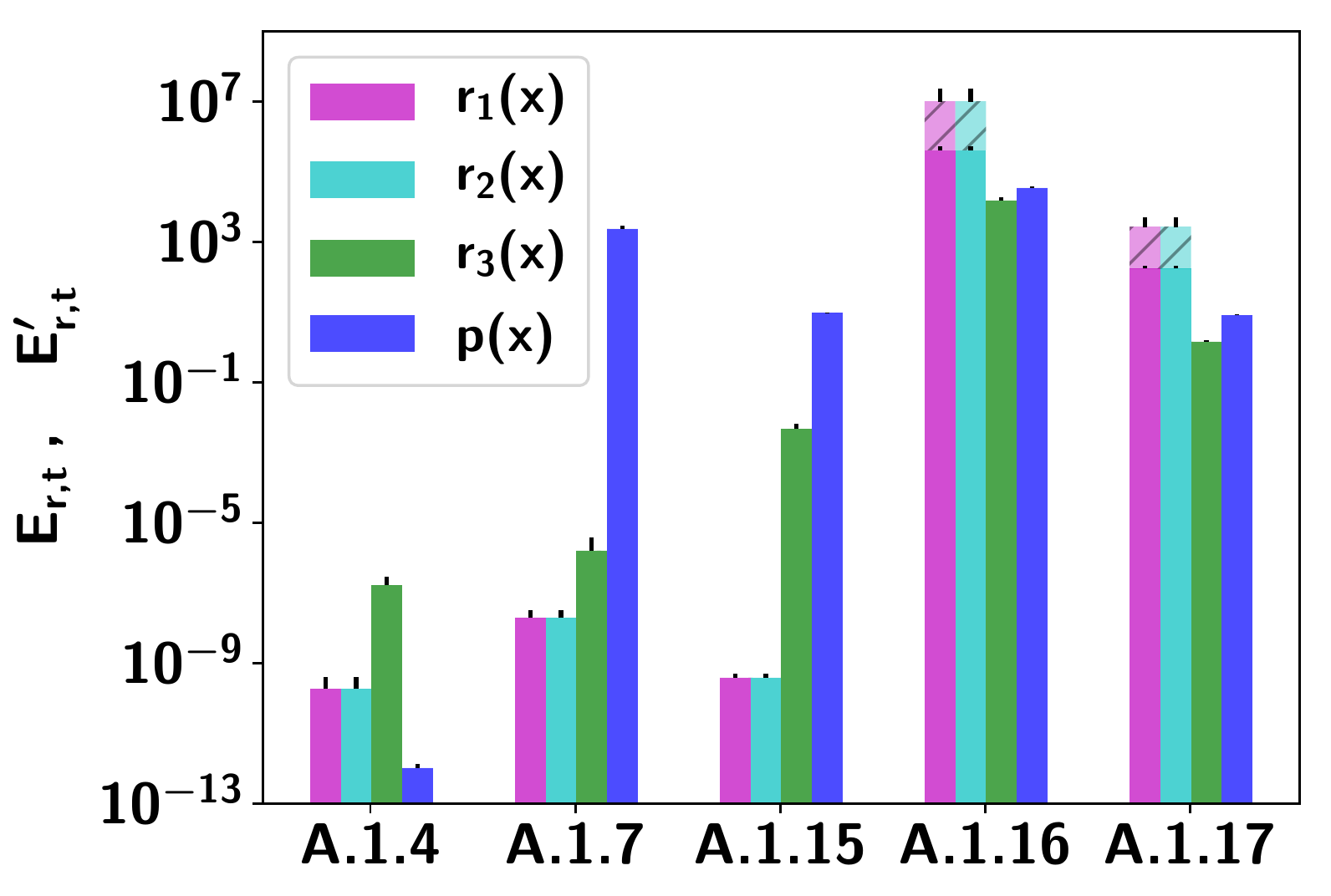}
			\caption%
			{{\small  $\epsilon=0$}}    
			\label{fig:errornonoise}
		\end{subfigure}
		\hfill
		\begin{subfigure}[b]{0.33\textwidth}  
			\centering 
			\includegraphics[width=1\textwidth]{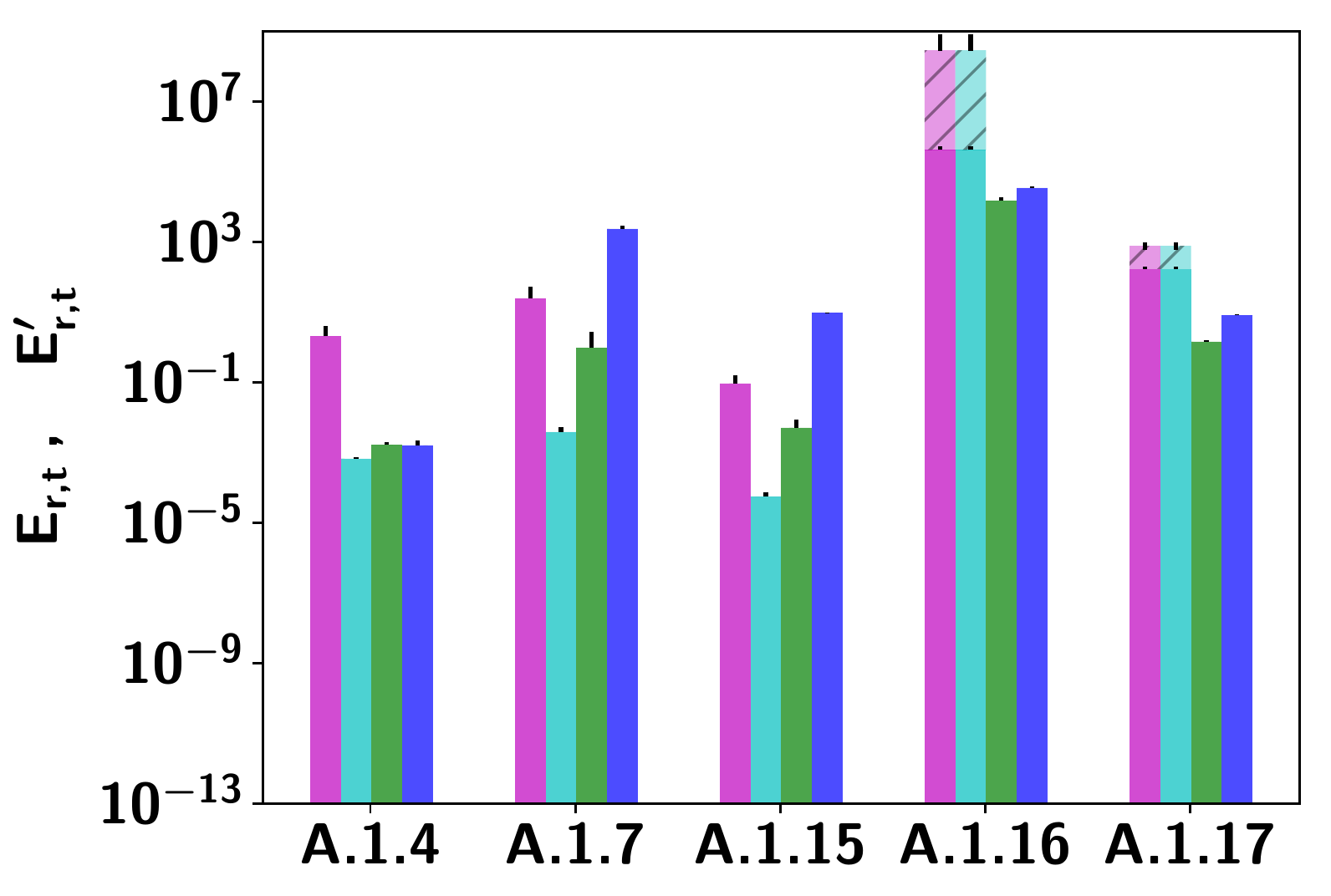}
			\caption[]%
			{{\small $\epsilon = 10^{-6}$}}    
			\label{fig:error10-6noise}
		\end{subfigure}
		\hfill
		\begin{subfigure}[b]{0.33\textwidth}  
			\centering 
			\includegraphics[width=1\textwidth]{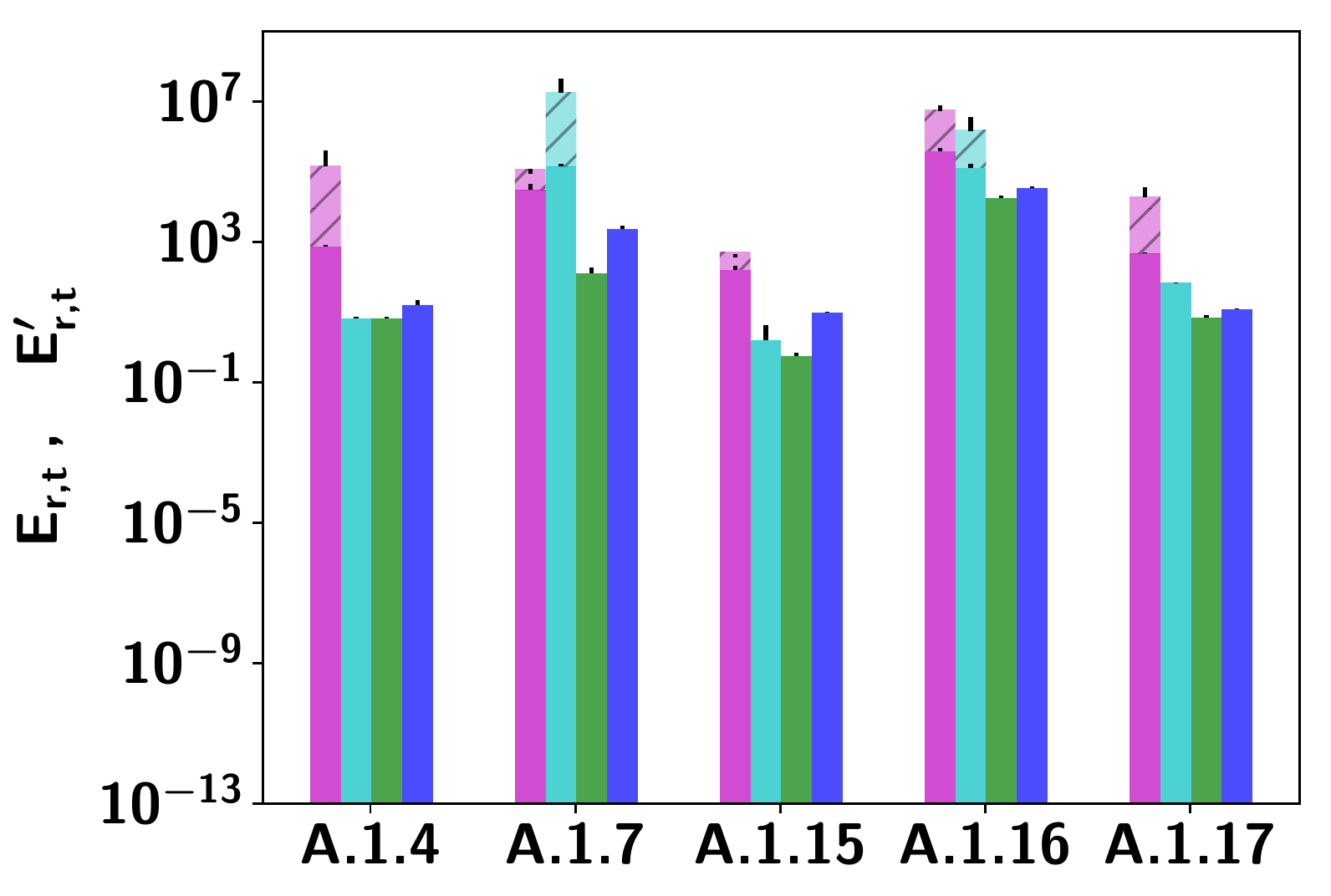}
			\caption[]%
			{{\small $\epsilon=10^{-2}$}}    
			\label{fig:error10-2noise}
		\end{subfigure}
		\caption[]
		{Comparison of the quality of rational and polynomial approximations. 
			The data are sampled with d-LHD, the threshold is $t=10^2$  and $\epsilon$ is the level of relative noise added to the data. 
			For approximations where pole-like points were found, the average error due to pole-like points is shown as a faded hatched bar.
			The average errors not due to polelike points are shown as  solid bars and are superimposed on the faded hatched bar where applicable.
			All the faded hatched bars, when not 0, are taller than the solid bars.
			The standard deviation  is shown as ablack vertical line. 
			\Cref{fn:f4} is an exponential function, \cref{fn:f8} is rational function, \cref{fn:f17} is a Breit-Wigner function, and \cref{fn:f18,fn:f19} are functions whose denominator is a polynomial. 
		}
		\label{fig:errorplot}
	\end{figure}

	\begin{table}[htb!]
				\footnotesize
		\centering
		\caption{
			Testing error ($\Delta_r$) for all test functions in \cref{T:TPs} of rational and polynomial approximations. 
			The data are sampled with d-LHD, and $\epsilon$ is the level of relative noise added to the data. 
			Since the scale of the error for each function is different, the error is first normalized to a 0-1 scale before calculating each statistic over all functions.  
			The corresponding results for each function can be found in the electronic supplement \cref{app:supppoleserrors}.
		} \label{tab:errorsummary}
				\setlength{\tabcolsep}{2.8pt}
		\begin{tabular}{|*{13}{c|}} 
			\cline{2-13}
			\multicolumn{1}{c}{}&\multicolumn{4}{|c|}{$\epsilon=0$}&\multicolumn{4}{|c|}{$\epsilon=10^{-6}$}&\multicolumn{4}{|c|}{$\epsilon=10^{-2}$}
			\\\cline{1-13}
			\multicolumn{1}{|c|}{Statistic}&
			\multicolumn{1}{|c|}{$r_1(x)$}&\multicolumn{1}{|c|}{$r_2(x)$}&\multicolumn{1}{|c|}{$r_3(x)$}&\multicolumn{1}{|c|}{$p(x)$}
			&\multicolumn{1}{|c|}{$r_1(x)$}&\multicolumn{1}{|c|}{$r_2(x)$}&\multicolumn{1}{|c|}{$r_3(x)$}&\multicolumn{1}{|c|}{$p(x)$}
			&\multicolumn{1}{|c|}{$r_1(x)$}&\multicolumn{1}{|c|}{$r_2(x)$}&\multicolumn{1}{|c|}{$r_3(x)$}&\multicolumn{1}{|c|}{$p(x)$}
			\\\hline\hline
			
			\makecell{Arithmetic \\Mean} &6.34E-02&6.34E-02&5.36E-02&8.38E-02&6.43E-02&5.33E-02&5.37E-02&8.38E-02&6.23E-02&7.49E-02&5.80E-02&8.37E-02
			\\\hline
			
			Median&1.88E-15&1.88E-15&4.48E-08&3.19E-04&1.85E-08&5.86E-07&5.16E-06&3.19E-04&3.65E-03&2.63E-03&8.71E-04&5.72E-04
			\\\hline
			
		\end{tabular}
	\end{table}

    To better compare the testing error, we divide it into two parts: the component due to poles and the remainder.
	Given the definition of $W_{r,t}$ in \cref{eq:noofpoles}, the error due to polelike points is defined as

	\begin{equation}\label{eq:errorfrompoles}
	E_{r,t} = \left[\sum\limits_{j \in W_{r,t}} (r(x^{(j)}) - f_j )^2\right]^{1/2}, 
	\end{equation}
	and the error not due to polelike points is
	\begin{equation}\label{eq:errornotfrompoles}
	E'_{r,t} = \left[\Delta_r^2 - E_{r,t}^2\right]^{1/2}.
	\end{equation}
The testing error for the three rational approximations $r_1(x)$, $r_2(x)$, and $r_3(x)$, as well as the polynomial approximation $p(x)$ is given in \cref{fig:errorplot} and \cref{tab:errorsummary}. The data for these plots are given in the electronic supplement \cref{app:supppoleserrors}. 
In order to ensure a fair comparison, the degrees of freedom are the same among the rational approximations. The degrees of freedom of the polynomial approximation are at least as large as those of the rational approximation.
For the rational approximations, the error due to pole-like points is given for the threshold value of $t = 10^2$.
From \cref{tab:errorsummary}, we observe that the approximation $r_3(x)$ performs best overall for all functions and all noise levels. More specifically, the approximation $r_3(x)$ has the lowest error when there is no noise or there are high levels of noise in the interpolation data. However, the quality of the approximation of $r_2(x)$ matches that of $r_3(x)$ when the level of noise is low ($\epsilon = 10^{-6}$). The reason is that the degree reduction in \cref{ALG:MVVandQR} is able to reduce poles and and give a better-quality approximation for low noise levels.
From \cref{fig:errorplot}, we observe that whenever polelike points are found, their contribution to the testing error is high, as defined in \cref{eq:errorfrompoles}. 
 The polynomial approximation yields a lower testing error than rational approximations do
 for the noise-free case of \cref{fn:f4} because  \cref{fn:f4} is approximated by a polynomial. 
 Conversely, for rational functions such as \cref{fn:f7}, the rational approximations over noise-free data yield better testing errors than the polynomial approximation does. 
 Moreover, for rational functions, the errors without polelike points in approximations $r_1(x)$ and $r_2(x)$ is on the order of $10^{-8}$ and  is lower than $10^{-6}$ for $r_3(x)$. 
 We believe the reason is that the approximations $r_1(x)$ and $r_2(x)$ are obtained from \cref{ALG:MVVandQR}, whose orthonormal basis implementation is numerically more accurate than the constrained optimization approach of \cref{A:Polyak} in the monomial basis.
 We also observe that the testing errors of the approximations of noise-free data of \cref{fn:f17} show trends similar to those described above for the rational functions. 
 The reason is that the denominator of this function approximates to a polynomial of degree 4 and has the unit of physical energy $E^4$. Since \cref{fn:f17} is unitless, the numerator also has a unit of $E^4$. Thus, the entire function can be approximated by a rational function with numerator of degree 4 and denominator of degree 4~\citep{PhysRev.49.519,Bohm:2004zi}.
 For other functions, the approximation $r_3(x)$ over noise-free data performs better than the other approximations due to the lack of spurious poles as well as due to better goodness of fit.

Generally all approximations for noise-free data are better than for noisy data. 
For $r_1(x)$ and $r_2(x)$, one reason is the higher number of  polelike points in the noisy data case. 
Another reason  is the poor quality of the fit of the approximations to the data. 
This is because the degrees of freedom are the same in both the noise-free and the noisy data cases. 
Hence, for higher levels of noise, the approximation underfits the data because there may not be any spare degrees of freedom to fit the data and the noise.
We would prefer to prevent overfitting the data, but for high noise  levels, more degrees of freedom may be required in order to better fit the data.

\begin{figure}[htb!]
	\centering
	\includegraphics[width=0.8\textwidth]{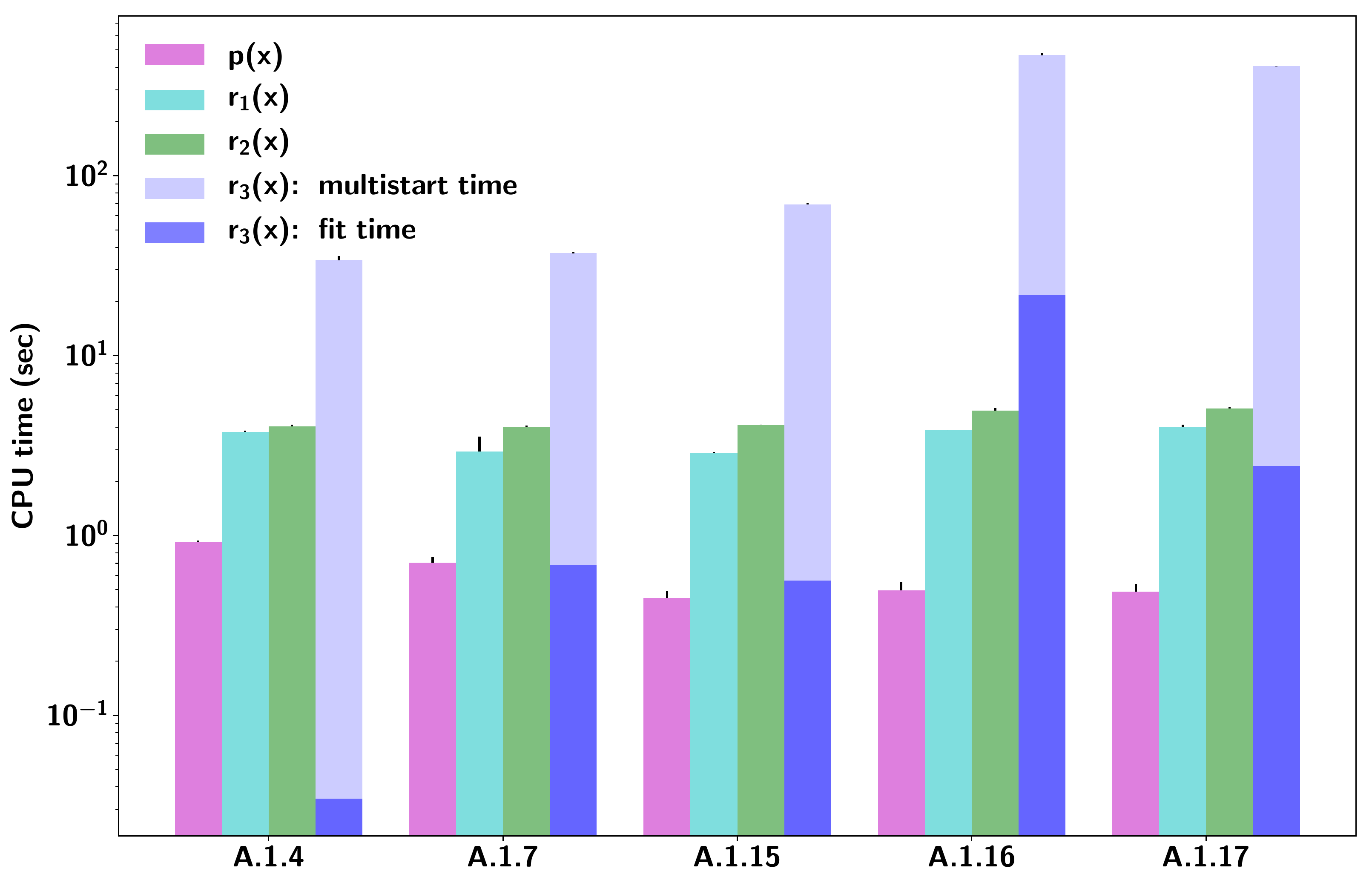}
	\caption{Total CPU time taken by the four approximation approaches when the interpolation data is sampled by using d-LHD and is noise free. 
		Each bar is the average CPU time, and the error bars at the top of each bar are the standard deviations. 
		The total time taken by \cref{A:Polyak} is shown as its fit time and multistart time to perform the global optimization of $q_l(x)$. 
		\Cref{fn:f4} is an exponential function, \cref{fn:f8} is a rational function, \cref{fn:f17} is a Breit-Wigner function, and \cref{fn:f18,fn:f19} are functions whose denominator is a polynomial. 
		The multistart time clearly dominates the total time taken by \cref{A:Polyak}. 
	}
	\label{fig:cputimeplot}
\end{figure}

	\begin{table}[htb!]
	\centering
	\caption{
		Total CPU time over all test functions in \cref{T:TPs} for all four approximation approaches when the interpolation data is sampled using d-LHD and is noise free. 
		The total time taken by \cref{A:Polyak} is shown as its fit time and multistart time to perform the global optimization of $q_l(x)$. 
		\cref{A:Polyak} is  more expensive than the other approaches  are, and this time is clearly dominated by the multistart time. 
		 The corresponding results for each function can be found in the electronic supplement \cref{app:supptimesiters}.
	} \label{tab:cputimesummary}
	\begin{tabular}{|*{6}{c|}} 
		\hline
		\multicolumn{1}{|c}{Statistic}&\multicolumn{1}{|c|}{$p(x)$}&\multicolumn{1}{|c|}{$r_1(x)$}&\multicolumn{1}{|c|}{$r_2(x)$}&\multicolumn{1}{|c|}{$r_3(x)$: Fit Time}&\multicolumn{1}{|c|}{$r_3(x)$: Multistart Time}
		\\\hline\hline
	Arithmetic Mean&0.66&3.17&4.00&8.60&88.95
	\\\hline
	Geometric Mean&0.63&3.11&3.94&0.19&31.44
	\\\hline
	Median&0.62&3.30&4.02&0.04&15.85
	\\\hline

	\end{tabular}
\end{table}

\subsection{Computational Effort of Computing Approximations}
In this section, we compare the computational effort required to compute all four approximations. 
\cref{fig:cputimeplot} and \cref{tab:cputimesummary} show the total CPU time taken by the four approximation approaches when the interpolation data are sampled by using d-LHD. 
In \cref{fig:cputimeplot} each bar is the average CPU time, and the error bars at the top of each bar indicate the standard deviation. 
The time taken by \cref{A:Polyak} is split into the time taken to fit the data by solving \cref{eq:rationalApproxFD} and the time to perform the global minimization of $q(x)$ by using the multistart approach across all iterations. 
Because the CPU times are generally  consistent across the different noise levels, we show the results only for the  noise-free case. 
The CPU times for all functions, sampling strategies and noise levels are given in the electronic supplement \cref{app:supptimesiters}.

We observe that the multistart time of \cref{A:Polyak} clearly dominates the total CPU time. The reason is that the fitting time grows with the number of iterations of \cref{A:Polyak}, as shown in \cref{fig:iteration}. 
However, the multistart time increases exponentially with the degrees of freedom, as is expected because the global optimization cost grows exponentially. 
On the other hand, the compute cost of the other approaches grows more slowly with the degrees of freedom because these algorithms are polynomial in time.
Despite the computation overhead, the cost of obtaining $r_3(x)$ may be negligible when they are used as surrogates for the expensive simulations of the physics processes, as we show in the next section.

\subsection{Summary of Computational Results}


From our experiments, we conclude that among the sampling strategies considered, d-LHD performs the best when the goal is to fit a rational approximation to data generated from a black box. We found this result to be true even when d-LHD was compared with the uniform random sampling of points over the domain.
The d-LHD method samples both on the faces and on the inside of the domain evenly and requires  only a few  interpolation points more than the degrees of freedom of the approximation. 
This is especially useful for applications whose function evaluations  are computationally  extremely expensive (minutes to hours per evaluation).
We also find that the approximations based on d-LHD-generated samples  require overall fewer  iterations of \cref{A:Polyak} and produce better-quality approximations than the LHS-based approximations do. 

The approximation approach using \cref{ALG:MVVandQR} with and without degree reduction is computationally more efficient than the approach using \cref{A:Polyak}.
However, our goal was do develop an approximation method for computationally expensive simulations that performs overall well when the underlying simulation function is unknown (black box). Thus, the computational overhead of the algorithms is negligible; and when applied to a true black-box simulation, \cref{A:Polyak} is more likely to give low errors, in particular when the interpolation data are noisy. More specifically, no spurious poles are found in $r_3(x)$ for nonrational functions as well as noisy problems, and the goodness of fit of $r_3(x)$ is much better than $r_1(x)$ or $r_2(x)$.
We note here that these claims are based on the assumption that  the data are sampled over a domain that does not include any true poles.

The approximation using \cref{A:Polyak} is computationally more expensive than the other approaches because it solves a harder problem of removing the poles iteratively.
Most of this expense is due to  the multistart optimization whose time grows exponentially with the degrees of freedom since it is tasked to the perform global minimization of $q(x)$.
As we will see in the next section, however, the additional time to get a pole-free and a superior quality of approximation is a small price to pay considering that this approximation replaces expensive simulations of the physics processes.

%% file: apps.tex
\section{Rational Approximation for High Energy Physics}\label{SEC:HEP}
\begin{figure}[htb!]
    \centering
    \includegraphics[width=.68\textwidth]{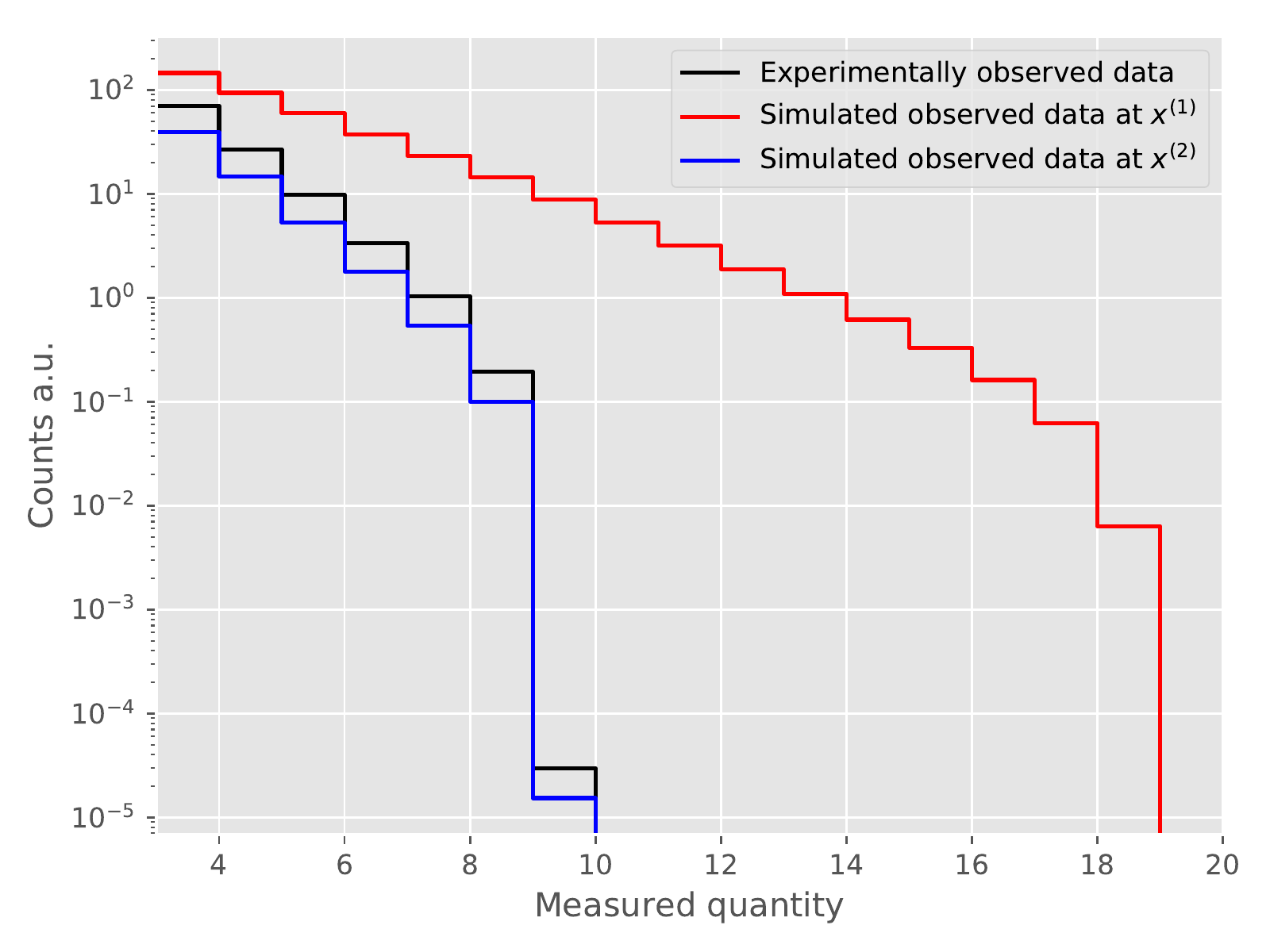}
    \caption{Illustration of a typical problem setup in particle physics. Shown are
        histogram skylines of observable experimental data (black) with observable
        predictions coming from simulations at different points $x^{(k)}$ in the same parameter space (blue and red).
        Numerical comparisons of the experimental with the simulated quantities are typically used to infer quantitative
        statements on the (unobservable) parameters of the simulation. Here, for example,, one would be interested in finding parameter points $x$ such that the corresponding simulation
       prediction resembles the experimental observation.}
    \label{fig:hepexample}
\end{figure}

A common problem in HEP is to infer information on \emph{unobservable}
parameters, $x$, from experimentally measured data, $d$. Typically, this is 
achieved by using a dedicated physics simulation program and statistical
measures. Since it is particularly well suited for rational approximations, we
will discuss a measure called ``binned likelihood'' in which the
measured and the simulated data take on the form of histograms with the same
binning \citep{Barlow:1993dm}. The binning of the histogram is driven by experimental
constraints such as how precise the quantity in question can actually be
measured. An illustrative example of the problem setup is shown in
\cref{fig:hepexample}.

In our example, the simulation predicts how postulated dark matter particles
interact with a Xenon-based detector in a process called ``direct detection''
(see \cite[Sec. 26]{PhysRevD.98.030001}). 
Our physics simulation has three parameters, $x=(m_\chi, c_+, c_\pi)$, which represent
the dark matter particle mass and two couplings to ordinary matter. The
parameter domain is $[10,100]\times[0.0001,0.001]\times[0.001,0.1]$.
The particle mass has the dimension of $\mathrm{GeV/c^2}$, while the couplings are dimensionless.

\sloppy We note that at the time of writing, no experimental result on the direct detection
of dark matter has been published.   Here, we assume a signal consistent with a dark matter mass $m_\chi=$ 10 $\mathrm{GeV/c^2}$ and an interaction strength large enough to produce approximately 100 events in future xenon detectors.   The simulated experimental data for each bin is:  $\{d_1, d_2, \ldots, d_6\} =\{ 70.4 , 26.7 , 9.8 , 3.4 , 1.0 , 0.2 \}$.
The values $d_b,  b=1,\ldots,6$,  are simulated by using a specific framework of the generalized
spin-independent response to dark matter in direct detection
experiments~\citep{Hoferichter:2016nvd,Cerdeno:2018bty}. 

At a given point $x^{(k)}$ in the parameter domain, we define the likelihood function
$\mathcal{L}(x^{(k)}|d_1, d_2, \ldots, d_6)$ as the product of  independent Poisson processes over all bins
(generalized for noninteger variables):
\begin{equation}\label{eq:likelihood}
    \mathcal{L}(x^{(k)}|d_1, d_2, \ldots, d_6) = \prod_{b=1}^6 \frac{N_b(x^{(k)}) d_b e^{N_b(x^{(k)})}}{\Gamma(d_b+1)},
\end{equation}
where the $N_b(x^{(k)})$ denote the simulated quantities for a point $x^{(k)}$
that correspond to the $d_b$; in other words, one simulation returns the values for all bins.   

By numerically maximizing \cref{eq:likelihood}, we infer information
about the parameters $x$ or regions of the parameter space  that yields
simulated data consistent with their experimentally observed counterpart.  We
use MultiNest~\citep{Feroz:2008xx,Feroz:2013hea,pymultinest} for this purpose,
which requires the evaluation of \cref{eq:likelihood} at tens of
thousands\footnote{The dimension of the problem and the convergence criteria of
    the MultiNest algorithm strongly influence the number of required function
calls.} of $x^{(k)}$ to succeed. The computational cost of this operation is
driven by the cost to obtain $N_b(x^{(k)})$ and can be  substantial.

In the following we will discuss  how  rational approximations can be used
to significantly reduce the required CPU cost of maximizing the likelihood. We will
show results for rational approximations of degree $M=4,N=4$ as well as polynomial approximations
of degree 7.\footnote{The degree is chosen such that the number of coefficients is comparable to the number of coefficients used in the
rational approximations.}

First, we calculate separate rational approximations $r_b(x)$ that approximates $N_b(x)$ for each bin $b$.   This calculation requires evaluating the exact simulation  at sufficiently many
training points. We use $N_\text{point}=500$ points sampled using the Latin hypercube method from the
parameter space, $x^{(k)}, k=1,\cdots,N_\text{points}$, at which we evaluate
 $N_b(x^{(k)}), k=1,\cdots,N_\text{points}$.  We compute the $r_b$ from the input-output data pairs
\begin{equation}\label{eq:ipoldata}
   \left\{\left( x^{(k)}, N_b(x^{(k)})\right)\right\}_{k=1}^{N_\text{points}} \text{ for }{b=1,\ldots,N_\text{bins}}.
\end{equation}
By replacing the expensive simulations to obtain  $N_b(x^{(k)})$ with cheap-to-evaluate rational (or polynomial) approximations
$r_b$ in \cref{eq:likelihood} we can define an approximate likelihood:
\begin{equation}\label{eq:likelihood-app}
    \mathcal{L}(x^{(k)}|d_1, d_2, \ldots, d_{N_{\text{bins}}}) \approx \mathcal{\tilde{L}}(x^{(k)}|d_1, d_2, \ldots, d_{N_{\text{bins}}}) = \prod_{b=1}^{N_{\text{bins}}} \frac{r_b(x^{(k)}) d_b e^{r_b(x^{(k)})}}{\Gamma(d_b+1)}.
\end{equation}
\sloppy The maximization of \cref{eq:likelihood} and \cref{eq:likelihood-app} requires
about 30,000 evaluations of $\mathcal{L}(x^{(k)}|d_1, d_2, \ldots, d_{N_{\text{bins}}})$ and
$\mathcal{\tilde{L}}(x^{(k)}|d_1, d_2, \ldots, d_{N_{\text{bins}}})$, respectively. The run time of the
latter is, however, about a factor 50 faster (\cref{tab:timeloli}).

\begin{figure}[htb!]
    \centering
    \begin{subfigure}[b]{0.475\textwidth}
        \centering
        \includegraphics[width=.9\textwidth]{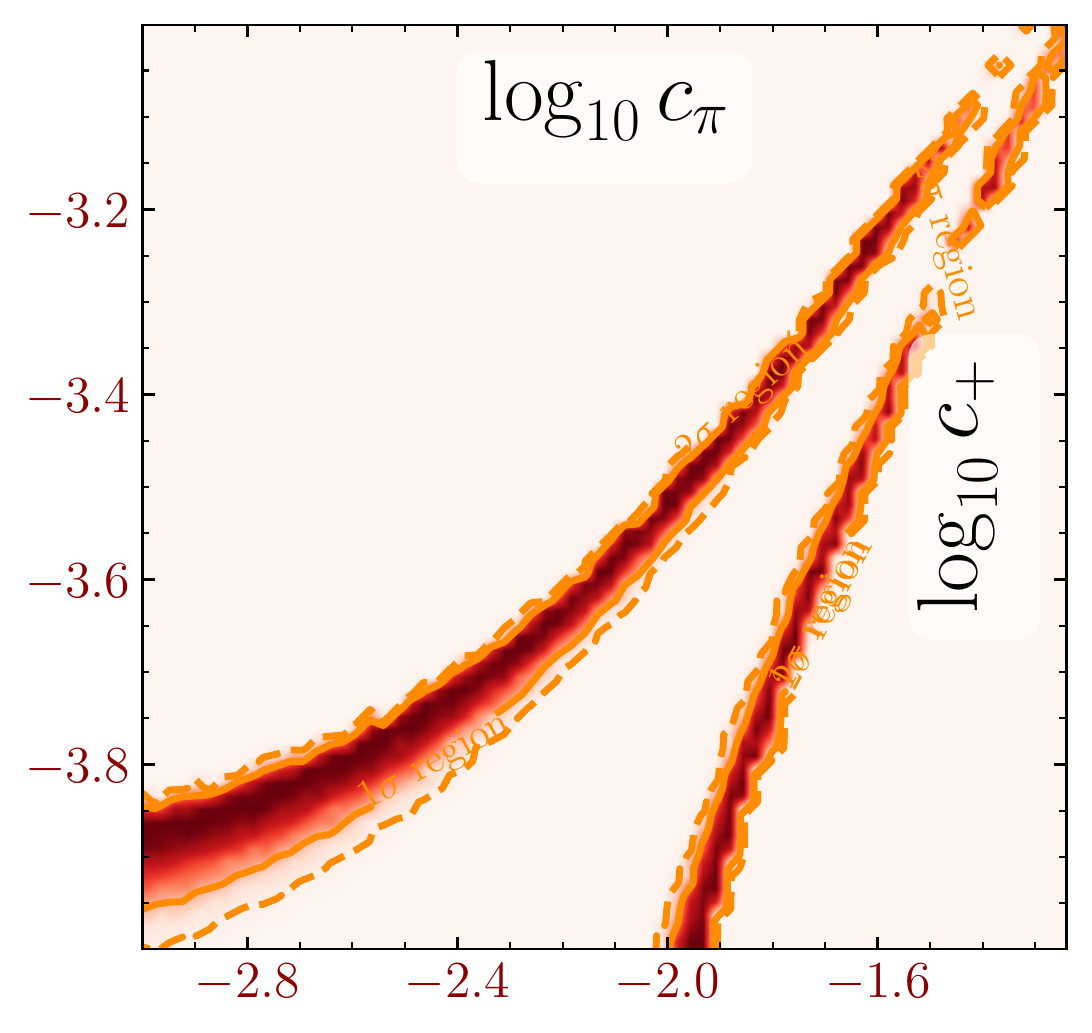}
        \caption%
        {{\small Result with full simulation}}    
        \label{fig:nest-si-full}
    \end{subfigure}
    \hfill
    \begin{subfigure}[b]{0.475\textwidth}  
        \centering 
        \includegraphics[width=.9\textwidth]{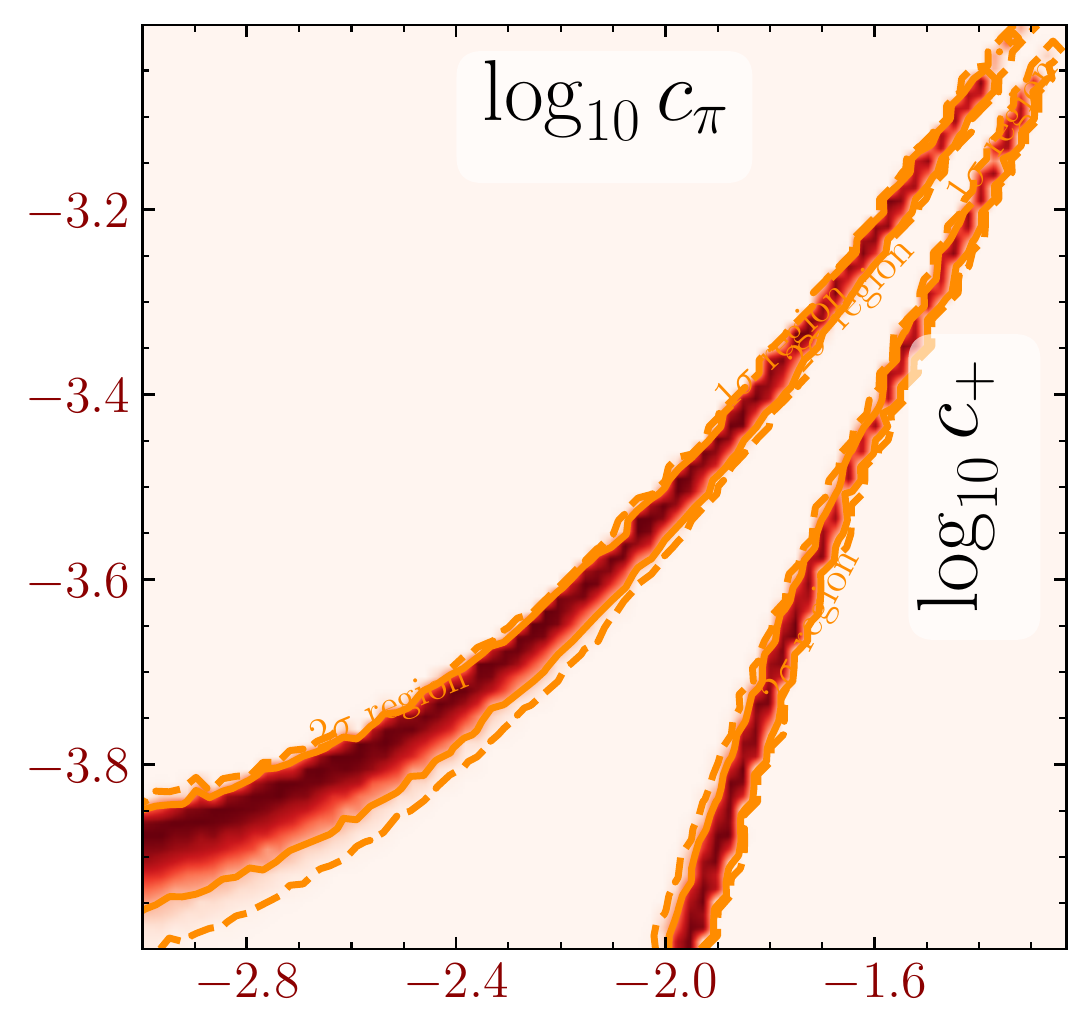}
        \caption[]%
        {{\small Result with pole-free rational approximation}}    
        \label{fig:nest-si-sip}
    \end{subfigure}
    \vskip\baselineskip
    \begin{subfigure}[b]{0.475\textwidth}   
        \centering 
        \includegraphics[width=.9\textwidth]{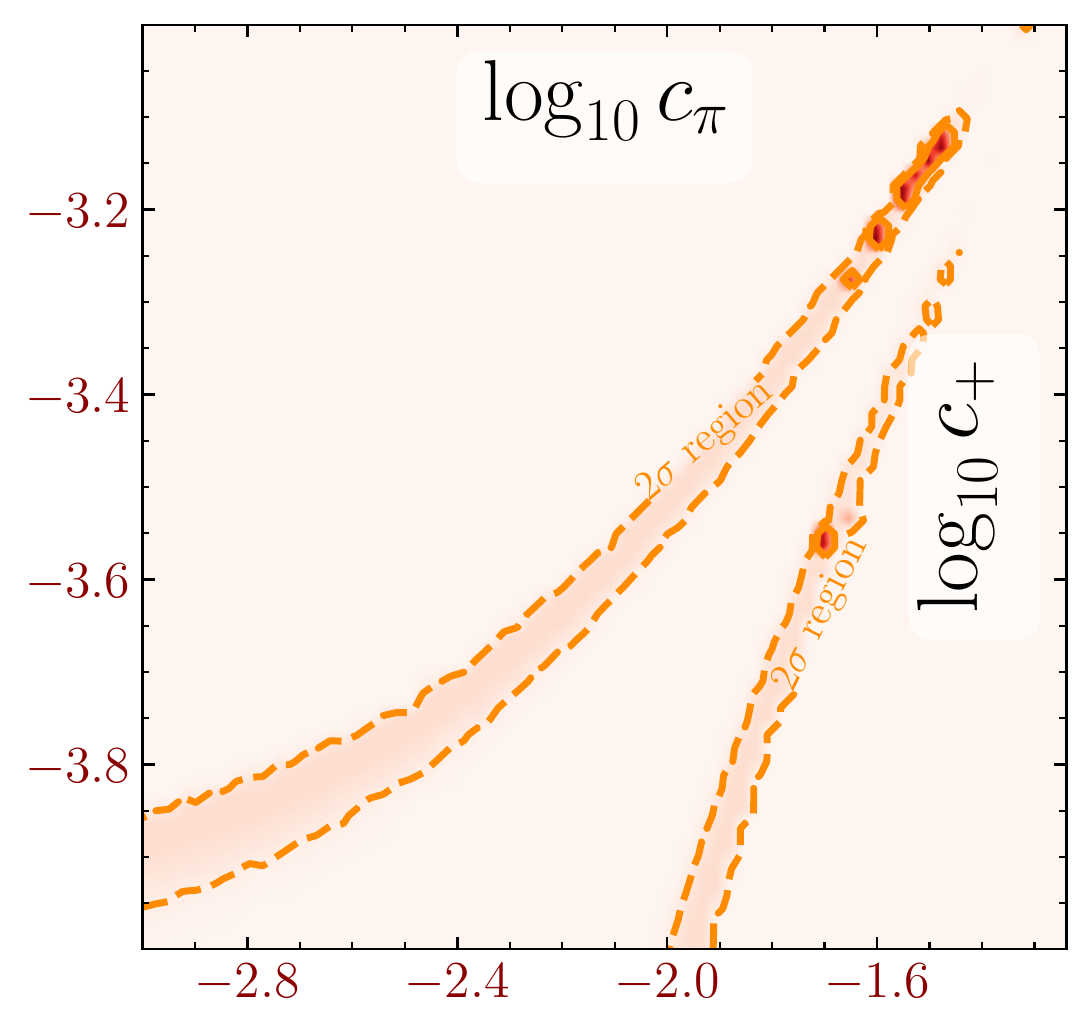}
        \caption[]%
        {{\small Result with non-pole-free rational approximation}}
        \label{fig:nest-si-la}
    \end{subfigure}
    \quad
    \begin{subfigure}[b]{0.475\textwidth}   
        \centering 
        \includegraphics[width=.9\textwidth]{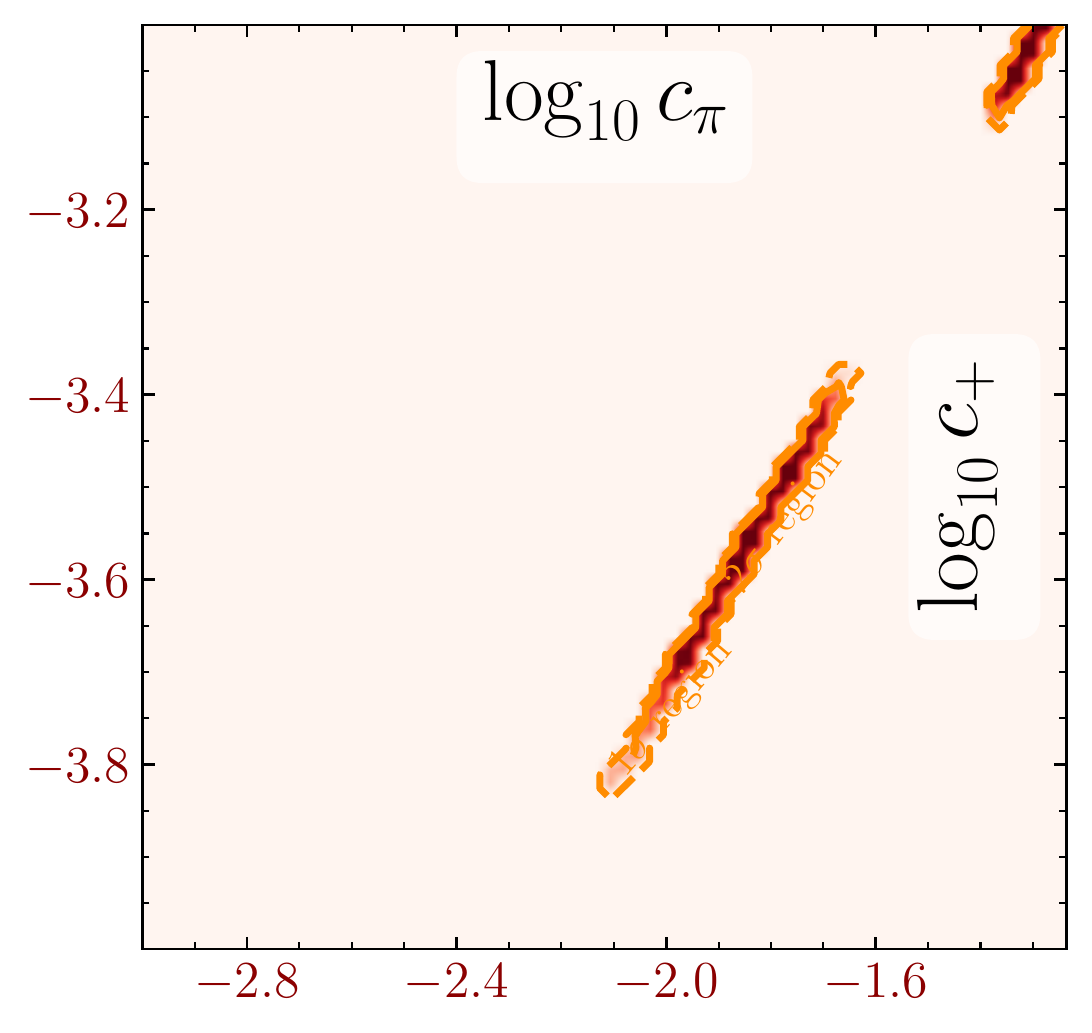}
        \caption[]%
        {{\small Result with polynomial approximation}}    
        \label{fig:nest-si-poly}
    \end{subfigure}
    \caption[]
    {Two-dimensional profile-likelihood projections of a 3-dimensional
        parameter space with superplot~\citep{Fowlie:2016hew}. Regions of higher likelihood are shown  darker. The data are
        normalized to the maximum likelihood observed before plotting. We
        compare the result obtained with the full physics simulation
        (\cref{fig:nest-si-full}) to the result obtained when using pole-free
        rational approximations ($M=4,N=4$) calculated with the semi-infinite
        approach (\cref{fig:nest-si-sip}). \cref{fig:nest-si-la} shows the
        effect of poles in the relevant parameter domain. The poles are visible as
        dark dots. For completeness, \cref{fig:nest-si-poly} shows the result when
        using polynomial approximations with a similar number of coefficients as in
        \cref{fig:nest-si-sip}.}

\label{fig:nest-si}
\end{figure}

\begin{table}[!h]
\caption{Comparison of computational cost when maximizing the likelihood \cref{eq:likelihood} using the true simulation and maximizing the approximate likelihood \cref{eq:likelihood-app} using a rational approximation for the data in each bin.}
    \label{tab:timeloli}
    \centering
    \begin{tabular}{l|c c}
        & Likelihood evaluations & total run-time~[s] \\ \hline
      using full simulation \cref{fig:nest-si-full} & 29459 & 14594 \\
      using $r_b$  with \cref{A:Polyak}, $M=4,N=4$ & 29612 & 288
    \end{tabular}
\end{table}
To demonstrate that the results obtained with the rational approximations $r_b$ are in agreement
with the full simulation, we present our results in terms of two-dimensional
profile-likelihood projections (\cref{fig:nest-si}). We limit the discussion to
the projection onto the $c_\pi-c_+$ plane because it exhibits the most interesting
pattern. In those plots, dark regions indicate higher likelihood values and
therefore high level of compatibility with experimentally observed data.

The top-left plot (\cref{fig:nest-si-full}) shows the result obtained with the
full simulation \cref{eq:likelihood}. We observe two ridges of equal
likelihood, meaning that there are very different parameter combinations that
are equally in agreement with the experimentally observed data. This is our
ground truth for comparison with the approximation based results.

The result obtained with pole-free (\cref{A:Polyak}) rational approximations in
\cref{fig:nest-si-sip} is in excellent agreement with the ground truth both
qualitatively and quantitatively. The rational approximations obtained with
\cref{EQN:LSQ-LA2}, shown in \cref{fig:nest-si-la}, demonstrates the impact of
spurious poles: although we find qualitative similarities with the ground truth,
the poles that are present in some of the $r_b$ lead to a complete distortion of the evaluated
likelihoods and therefore to a quantitatively wrong interpretation.
For completeness, we show in \cref{fig:nest-si-poly} the result obtained with polynomials of order 7.
It clearly shows the advantage of rational approximations since the polynomial approximations are apparently not able to capture the
true likelihood at all. Thus, the information inferred by using the polynomial
approximation would be misleading.


%% file: concl.tex
\section{Conclusions}\label{SEC:Conclusion}
We have presented two approaches for computing rational approximations for computationally expensive black-box functions. Our first approach uses linear algebra to construct the rational approximation, but it does not guarantee that the approximation is pole free. Our second approach exploits a  semi-infinite optimization problem formulation that leads to accurate rational approximations without poles.

Our numerical study shows that the selection of interpolation points for fitting the approximations has a major impact on the approximation error and the number of iterations taken by the pole-free rational approximation. We find that a Latin hypercube design that is augmented with sample points on the boundary of the  parameter domain leads to improved approximations more efficiently.  We hypothesize  that this is due to close proximity of the function domain to the true poles and the approximations fare poorly without the sample points on the boundary. We showed that for a variety of analytic fast-to-compute test problems with and without noise the rational approximations generally perform better than the polynomial approximations do. The result was further confirmed by approximating data generated from an expensive  HEP simulation. The polynomial did not capture the true underlying functional relationship at all. Thus, for black-box simulations whose true underlying functional forms are unknown, using a polynomial may lead to incorrect conclusions.

An outstanding challenge for using rational approximations is the determination of the ``correct'' polynomial degrees in the numerator and denominator. We have experimented with a heuristic method to determine these degrees, but noisy data pose an additional challenge, and more research is needed. 

The structural constraint considered in the pole-free rational approximation is mitigating poles through enforcing non-negativity of the denominator $q(x)$. Other structural constraints arise in the solution of chance-constraint optimization, where we wish to approximate an empirical cumulative density function. By construction, the function should be monotonic, which again imposes a constraint on the rational approximation. Such structural constraints should also be modeled in the future.

The rational approximations require a minimum number of interpolation points to fit the model. One drawback is in obtaining these interpolation points, since the number of points required increases significantly with the number of parameters and the degrees of the polynomials of the approximation. Hence, obtaining these interpolation points may become computationally too expensive. 
Additionally, the multistart global optimization of the denominator $q(x)$ in the pole-free rational approximation will become computationally significantly more expensive as the number of parameters increase.  We have tested problems with up to 7 parameters, but especially in high energy physics dozens of parameters are commonly encountered. 
Thus, the question of scalability of the proposed rational approximation approaches must be addressed in the future.

%% file: testProblems.tex
\section{Description of Test Problems}\label{SEC:TPs}

The functional forms of the test problems that we use in our numerical experiments are shown in Table~\ref{T:TPs}. 



\begin{table}[htb!]
	\centering
	\small
	\caption{Description of fast-to-compute test problems. Here, $n$ is the number of variables, $M$ is the degree of the numerator, $N$ is the degree of the denominator, $f$ is the functional form, and \textit{Domain} for each dimension is the interval in which no  poles exist.
If either the numerator or denominator is not a polynomial, then the entry for $M$ or $N$ is a dash, respectively. } \label{T:TPs}
	\begin{tabular}{|N|c|c|c|c|c|c|}
		\hline
		\multicolumn{1}{|c|}{No.} & Description & $n$ &$M$&$N$& $f$ & \textit{Domain}  \\ \hline\hline 
		\label{fn:f1} & \makecell{Function whose denominator \\ is a polynomial}  & 2 & -- & 4  & 
			$\dps \frac{e^{x_1x_2}}{(x_1^2-1.44)(x_2^2-1.44)}$ & $x \in [-1,1]^2$\\ [0.2em] \hline
			
		\label{fn:f2} & Log function  & 2 & -- & -- & 
			$\dps \log(2.25-x_1^2-x_2^2)$ & $x \in [-1,1]^2$ \\ [0.2em] \hline
			
		\label{fn:f3} & Hyperbolic tangent function  & 2 & -- & --  & 
			$\dps \tanh(5(x_1-x_2))$ & $x \in [-1,1]^2$ \\ [0.2em] \hline
			
		\label{fn:f4} & Exponential function  & 2 & -- & -- & 
			$\dps e^{\frac{-(x_1^2+x_2^2)}{1000}}$ & $x \in [-1,1]^2$ \\ [0.2em] \hline
			
		\label{fn:f5} & Absolute value function  & 2 & -- & -- & 
			$\dps |(x_1-x_2)|^3$ & $x \in [-1,1]^2$ \\ [0.2em] \hline
			
		\label{fn:f7} & \makecell{Rational function}  & 2 & 3 & 3  & 
			$\dps \frac{x_1+x_2^3}{x_1x_2^2+1}$ & $x\in [0,1]^2$ \\ [0.2em] \hline
			
		\label{fn:f8} & \makecell{Rational function} & 2 &2& 2 & 
			$\dps \frac{x_1^2+x_2^2+x_1-x_2-1}{(x_1-1.1)(x_2-1.1)}$ & $x \in [-1,1]^2$ \\ [0.2em] \hline
			
		\label{fn:f9} & \makecell{Rational function}  & 2 &4 &4  & 
			$\dps \frac{x_1^4+x_2^4+x_1^2x_2^2+x_1x_2}{(x_1^2-1.1)(x_2^2-1.1)}$ & $x \in [-1,1]^2$ \\ [0.2em] \hline
			
		\label{fn:f10} & \makecell{Rational function} & 4 & 2& 2 & 
			$\dps \frac{x_1^2+x_2^2+x_1-x_2+1}{(x_3-1.5)(x_4-1.5)}$ & $x \in [-1,1]^4$ \\ [0.2em] \hline
			
		\label{fn:f12} & \makecell{Rational function}  & 2 & 2 &3  & 
			$\dps \frac{x_1^2 + x_2^2 + x_1 - x_2 - 1}{x_1^3+x_2^3+4}$ & $x \in [-1,1]^2$ \\ [0.2em] \hline
			
		\label{fn:f13} & \makecell{Rational function}  & 2 & 3 &2  & 
			$\dps \frac{x_1^3+x_2^3}{x_1^2 + x_2^2 + 3}$ & $x \in [-1,1]^2$\\ [0.2em] \hline
			
		\label{fn:f14} & \makecell{Rational function}  & 2 & 4 & 4  & 
			$\dps \frac{x_1^4 + x_2^4 + x_1^2x_2^2 + x_1x_2}{x_1^2x_2^2-2x_1^2-2x_2^2+4}$ &$x \in [-1,1]^2$ \\ [0.2em] \hline
			
		\label{fn:f15} & \makecell{Rational function}  & 2 & 3 &4  & 
			$\dps \frac{x_1^3+x_2^3}{x_1^2x_2^2-2x_1^2-2x_2^2+4}$ &$x \in [-1,1]^2$ \\ [0.2em] \hline
			
		\label{fn:f16} & \makecell{Rational function} & 2 & 4 & 3  & 
			$\dps \frac{x_1^4 + x_2^4 + x_1^2x_2^2 + x_1x_2}{x_1^3+x_2^3+4}$ &$x \in [-1,1]^2$ \\ [0.3em] \hline
			
		\label{fn:f17} & Breit-Wigner function  & 3 & -- & --  & 
			\makecell{$\frac{2\sqrt{2}M\Gamma\gamma}{(\pi\sqrt{M^2+\gamma})[(E^2-M^2)^2+M^2\Gamma^2]}$\\ where $\gamma = \sqrt{M^2(M^2+\Gamma^2)}$} & \makecell{$E \in [80,100],$\\ $\Gamma \in [5,10],$\\ $M \in [90,93] $}\\ [0.2em] \hline
			
		\label{fn:f18} & \makecell{Function whose denominator \\ is a polynomial} & 4 & -- & 4 & 
			$\dps \frac{\tan^{-1}{(x_1)}+\dots+\tan^{-1}{(x_4)}}{x_1^2x_2^2-x_1^2-x_2^2+1}$ & $x \in [-0.95,0.95]^4$\\ [0.4em] \hline
			
		\label{fn:f19} & \makecell{Function whose denominator \\ is a polynomial}  & 4 & -- & 2  & 
			$\dps \frac{e^{x_1x_2x_3x_4}}{x_1^2+x_2^2-x_3x_4+3}$ & $x \in [-1,1]^4$\\ [0.4em] \hline
			
		\label{fn:f20} & Sinc function  & 4 & -- & -- & 
			$\dps 10\prod_{i=1}^{4}\frac{\sin{x_i}}{x_i}$ & $x \in [10^{-6},4\pi]^4$\\ [0.2em] \hline
			
		\label{fn:f21} & Sinc function  & 2 &-- & -- & 
			$\dps 10\frac{\sin{x_1}}{x_1}\frac{\sin{x_2}}{x_2}$ & $x \in [10^{-6},4\pi]^2$\\ [0.2em] \hline
			
		\label{fn:f22} & \makecell{Polynomial function}  & 2 & 2 & --  & 
			$\dps x_1^2 + x_2^2+ x_1x_2 - x_2+1$ & $x\in[-1,1]^2$\\ [0.2em] \hline
		
	\end{tabular}
\end{table}

%

%% file: multistartstrategy.tex
\section{Checking for Poles}\label{app:globaloptprob}

In this section, we discuss practical ways to solve the global optimization subproblem in line 4 of \cref{A:Polyak}. 
We compare different strategies to perform the global minimization of $q_l(x)$ to detect poles in $D$ for a number of fast-to-compute test functions described in \cref{T:TPs}. 
The benchmark for the comparison is the Baron global optimization solver for nonlinear and mixed-integer nonlinear problems \citep{ts:05, sahinidis:baron:17.8.9}. 
The other strategies include ``singlestart,'' in which we choose  one  point randomly  from $D$ as starting point for the optimization; ``multistart,'' which starts multiple optimizations  from different  points in $D$; and ``sampling,'' where $q_l(x)$ is evaluated at multiple random points in $D$ to check if any evaluation of $q_l(x) < 0$.
We allow multistart and sampling  to run for the same amount of time as Baron to ensure a fair comparison of these approaches. 
However,  multistart and sampling stop as soon as the first $x$ with $q_l(x) < 0$ is detected and do not continue toward finding the global minimum.

The results from this comparison are summarized in \cref{tab:baronComparison}.
We observe that multistart detects poles almost as well as Baron in a much shorter time. 
The reason is that multistart stops as soon as some $x$ with $q_l(x) < 0$ is detected, whereas Baron tries to solve the problem to optimality in each iteration of \cref{A:Polyak}.
Also, multistart detects poles almost as well as Baron does when the time taken by both approaches is the same. 
Therefore, to set a suitable time limit for multistart a priori, we estimate the amount of time Baron would take to solve the problem given the number of nonlinearities. 
Then we compute the number of multistart iterations that can be completed within this time.
The goal of this heuristic is to minimize the occurrence of poles in $q(x)$ without spending the effort required to run Baron. 
The number of multistart iterations needed is approximately an exponential function, $\phi$, of the number of nonlinearities, nnl,
when multistart ran for the same time as Baron. 
\begin{equation}\label{eq:multistart}
\phi(\mathrm{nnl}) = 2042.023e^\mathrm{\textit{0.029}nnl}
\end{equation} 

\begin{table}[htb!]
	\caption{Comparing global optimization strategies for detecting poles. Here,  $n$ denotes the number of variables, \textit{nnl} is the number of nonlinearities in $q_l(x)$ that is obtained by subtracting the constant and linear terms from the total degrees of freedom of $q_l(x)$. 
\textit{Time} is the CPU time in seconds, and \textit{\%FN} is the percentage of false negatives for detecting  poles, thst is, when Baron identifies the existence of a pole, while the corresponding other method did not. 
The results are more informative for problems where $n>2$. Hence, results are only obtained for some functions with $n=2$.} \label{tab:baronComparison}
		\centering
		\begin{tabular}{ccc|c|c|c|c|c|c|c|}
			\cline{4-10}
			&&& Baron&
			\multicolumn{2}{|c|}{Singlestart}
			&\multicolumn{2}{|c|}{Multistart}
			&\multicolumn{2}{|c|}{Sampling}\\
			\cline{1-10}
			\multicolumn{1}{|c|}{Function No.}
			&\multicolumn{1}{|c|}{$n$}
			&\multicolumn{1}{|c|}{\shortstack{\textit{nnl}}}
			&\shortstack{\textit{Time}}
			&\shortstack{\textit{\% FN}}
			&\shortstack{\textit{Time}}
			&\shortstack{\textit{\% FN}}
			&\shortstack{\textit{Time}}
			&\shortstack{\textit{\% FN}}
			&\shortstack{\textit{Time}}
			\\
			\cline{1-10}
			\multicolumn{1}{|c|}{\ref{fn:f14}}&\multicolumn{1}{|c|}{2}&\multicolumn{1}{|c|}{7}&0.0809&0.68&0.0021&0.00&0.0319&1.35&0.0306\\
			\multicolumn{1}{|c|}{\ref{fn:f15}}&\multicolumn{1}{|c|}{2}&\multicolumn{1}{|c|}{7}&0.0575&2.67&0.0017&0.00&0.0541&2.00&0.0539\\
			\multicolumn{1}{|c|}{\ref{fn:f16}}&\multicolumn{1}{|c|}{2}&\multicolumn{1}{|c|}{7}&0.0564&1.29&0.0018&0.00&0.0506&0.65&0.0503\\
			\multicolumn{1}{|c|}{\ref{fn:f17}}&\multicolumn{1}{|c|}{3}&\multicolumn{1}{|c|}{16}&0.1066&9.66&0.0057&0.00&0.0742&1.70&0.0743\\
			\multicolumn{1}{|c|}{\ref{fn:f18}}&\multicolumn{1}{|c|}{4}&\multicolumn{1}{|c|}{30}&0.2653&23.63&0.0082&1.10&0.0757&19.23&0.1270\\
			\multicolumn{1}{|c|}{\ref{fn:f19}}&\multicolumn{1}{|c|}{4}&\multicolumn{1}{|c|}{30}&0.1202&0.00&0.0051&0.00&0.0539&5.00&0.0756\\
			\multicolumn{1}{|c|}{\ref{fn:f20} in 7D}&\multicolumn{1}{|c|}{7}&\multicolumn{1}{|c|}{112}&259.0448&0.29&0.0078&0.00&0.3549&2.20&0.4579\\
			\hline
		\end{tabular}
\end{table}




%% file: suppPolesAndError.tex
\section{Number of Pole-like Points and Errors}\label{app:supppoleserrors}

In \cref{fig:poleplot}, we showed the average number of pole-like points found over all functions in \cref{T:TPs} in each rational approximation for different noise levels when the interpolation data was sampled using d-LHS. Additionally, in \cref{fig:poleplot}, we compared the quality of the four approximation approaches for five typical test functions from \cref{T:TPs} whose interpolation data is sampled using d-LHS. We summarized these results over all functions in \cref{tab:errorsummary}. 
In this section we present detailed results of the pole-like points and error (due to poles, not due to poles and total testing error) found in the approximation approaches of all the functions from \cref{T:TPs}. These detailed results are given for when the interpolation data is sampled using the three strategies of SG, LHS, d-LHD.  The results are presented for noise-free data ($\epsilon = 0$) as well as for data with the relative noise level of $\epsilon = 10^{-6}$ and $\epsilon = 10^{-2}$.

The number of pole-like points and error results for approximating noise-free interpolation data of all functions from \cref{T:TPs} are shown in \cref{tab:supoleerror0}. In this table, $p(x)$ is the polynomial approximation, (b) $r_1(x)$ is the rational approximation using \cref{ALG:MVVandQR} without degree reduction, (c) $r_2(x)$ is the rational approximation using \cref{ALG:MVVandQR} with the degree reduction described in \cref{ALG:ReduceDegree}, and (d) $r_3(x)$ is the rational approximation using \cref{A:Polyak}. 
The number of pole-like points found on the face and inside the domain (see \cref{eq:noofpoles}) is given. 
The error due to poles (see \cref{eq:errorfrompoles}), not due to poles (see \cref{eq:errornotfrompoles}), and total testing error (see \cref{eq:testError}) is also given. 
Because pole-like points and errors related to these points are only applicable for rational approximations, we only give these results for $r_1(x), r_2(x)$, and $r_3(x)$, and put a ``'-" in its place for $p(x)$ instead.
\textit{Function No.} is the number of the function from \cref{T:TPs}. 
\textit{Sample Type} is the name of the sampling strategy used to pick interpolation points from the domain. 
\textit{AM} is the arithmetic mean and \textit{MD} is the median of the results from running the experiment over LHS and d-LHD samples with different random seeds. 
Since SG is deterministic in picking points from the domain, a ``-" is placed for \textit{MD} whenever the sampling strategy is SG.
Similarly, the timing and iteration results for approximating noisy data of $\epsilon=10^{-6}$ and $\epsilon=10^{-2}$ is given in \cref{tab:supoleerror10-6} and \cref{tab:supoleerror10-2}, respectively.

{\setlength{\tabcolsep}{4pt}
	\begin{center}
		\begin{footnotesize}
\begin{longtable}{|*{11}{c|}}

				\caption{
					Number of pole-like points found on the face and inside the domain as well as errors found in all four approximations of functions from \cref{T:TPs}. 
					The interpolation data for these functions are noise-free, i.e., $\epsilon=0$. 
					Number of pole-like points and the error due to poles and not due to poles are given for $r_1(x), r_2(x)$, and $r_3(x)$ only. A  ``-" is used in its place for $p(x)$.
					The number of pole-like points, error due to poles and not due to poles is given for the threshold value of $t=10^2$.
					\textit{Function No.} is the number of the function from \cref{T:TPs}. 
					\textit{Sample Type} is the name of the sampling strategy used to pick interpolation points from the domain. 
					\textit{AM} is the arithmetic mean and \textit{MD} is the median.
					Since SG is deterministic in picking points from the domain, a ``-" is placed for \textit{MD} whenever the sampling strategy is SG.
				} 
				\label{tab:supoleerror0}\\
				\cline{4-11}
				\multicolumn{3}{c|}{}&\multicolumn{8}{|c|}{$\epsilon=0$}
				\\
				\cline{1-11}
				\multirow{2}{*}{\makecell{\textit{Function} \\ \textit{No.}}} 
				&\multirow{2}{*}{\makecell{\textit{Sample}\\ \textit{Type}}} 
				&\multirow{2}{*}{} 
				&\multicolumn{2}{|c|}{$r_1(x)$}
				&\multicolumn{2}{|c|}{$r_2(x)$}
				&\multicolumn{2}{|c|}{\makecell{$r_3(x)$}}
				&\multicolumn{2}{|c|}{$p(x)$}
				\\
				\cline{4-11}
				&&
				&\textit{AM} & \textit{MD}
				&\textit{AM} & \textit{MD}
				&\textit{AM} & \textit{MD}
				&\textit{AM} & \textit{MD}
				\\\hline
				\endfirsthead
				
				\cline{4-11}
				\multicolumn{3}{c|}{}&\multicolumn{8}{|c|}{$\epsilon=0$}
				\\
				\cline{1-11}
				\multirow{2}{*}{\makecell{\textit{Function} \\ \textit{No.}}} 
				&\multirow{2}{*}{\makecell{\textit{Sample}\\ \textit{Type}}} 
				&\multirow{2}{*}{} 
				&\multicolumn{2}{|c|}{$r_1(x)$}
				&\multicolumn{2}{|c|}{$r_2(x)$}
				&\multicolumn{2}{|c|}{\makecell{$r_3(x)$}}
				&\multicolumn{2}{|c|}{$p(x)$}
				\\
				\cline{4-11}
				&&
				&\textit{AM} & \textit{MD}
				&\textit{AM} & \textit{MD}
				&\textit{AM} & \textit{MD}
				&\textit{AM} & \textit{MD}
				\\\hline
				
				\endhead
				
				\multirow{15}{*}{\ref{fn:f1}}&\multirow{5}{*}{SG}&$|W^{(fc)}_{r,10^2}|$&0&-&0&-&0&-&-&-
				\\*  \cline{3-11}
				&&$|W^{(in)}_{r,10^2}|$&0&-&0&-&0&-&-&-
				\\*  \cline{3-11}
				&&$E_{r,10^2}$&0&-&0&-&0&-&-&-
				\\*  \cline{3-11}
				&&$E'_{r,10^2}$&4.64&-&4.64&-&5.30E+02&-&-&-
				\\*  \cline{3-11}
				&&$\Delta_r$&4.64&-&4.64&-&5.30E+02&-&71.33&-
				\\*  \cline{2-11}
				&\multirow{5}{*}{LHS}&$|W^{(fc)}_{r,10^2}|$&0.60&0&0.60&0&0&0&-&-
				\\*  \cline{3-11}
				&&$|W^{(in)}_{r,10^2}|$&0.20&0&0.20&0&0&0&-&-
				\\*  \cline{3-11}
				&&$E_{r,10^2}$&1.84E+03&0&1.84E+03&0&0&0&-&-
				\\
\hline

			\end{longtable}
		\end{footnotesize}
	\end{center}
}

{\setlength{\tabcolsep}{4pt}
	\begin{center}
		\begin{footnotesize}
			\begin{longtable}{|*{11}{c|}}

				\caption{
					Number of pole-like points found on the face and inside the domain as well as errors found in all four approximations of functions from \cref{T:TPs}. 
					The interpolation data for these functions have a relative noise level of $\epsilon=10^{-6}$.
					Number of pole-like points and the error due to poles and not due to poles are given for $r_1(x), r_2(x)$, and $r_3(x)$ only. A  ``-" is used in its place for $p(x)$.
					The number of pole-like points, error due to poles and not due to poles is given for the threshold value of $t=10^2$.
					\textit{Function No.} is the number of the function from \cref{T:TPs}. 
					\textit{Sample Type} is the name of the sampling strategy used to pick interpolation points from the domain. 
					\textit{AM} is the arithmetic mean and \textit{MD} is the median.
					Since SG is deterministic in picking points from the domain, a ``-" is placed for \textit{MD} whenever the sampling strategy is SG.
				} 
				\label{tab:supoleerror10-6}\\
				\cline{4-11}
				\multicolumn{3}{c|}{}&\multicolumn{8}{|c|}{$\epsilon=10^{-6}$}
				\\
				\cline{1-11}
				\multirow{2}{*}{\makecell{\textit{Function} \\ \textit{No.}}} 
				&\multirow{2}{*}{\makecell{\textit{Sample}\\ \textit{Type}}} 
				&\multirow{2}{*}{} 
				&\multicolumn{2}{|c|}{$r_1(x)$}
				&\multicolumn{2}{|c|}{$r_2(x)$}
				&\multicolumn{2}{|c|}{\makecell{$r_3(x)$}}
				&\multicolumn{2}{|c|}{$p(x)$}
				\\
				\cline{4-11}
				&&
				&\textit{AM} & \textit{MD}
				&\textit{AM} & \textit{MD}
				&\textit{AM} & \textit{MD}
				&\textit{AM} & \textit{MD}
				\\\hline
				\endfirsthead
				
				\cline{4-11}
				\multicolumn{3}{c|}{}&\multicolumn{8}{|c|}{$\epsilon=10^{-6}$}
				\\
				\cline{1-11}
				\multirow{2}{*}{\makecell{\textit{Function} \\ \textit{No.}}} 
				&\multirow{2}{*}{\makecell{\textit{Sample}\\ \textit{Type}}} 
				&\multirow{2}{*}{} 
				&\multicolumn{2}{|c|}{$r_1(x)$}
				&\multicolumn{2}{|c|}{$r_2(x)$}
				&\multicolumn{2}{|c|}{\makecell{$r_3(x)$}}
				&\multicolumn{2}{|c|}{$p(x)$}
				\\
				\cline{4-11}
				&&
				&\textit{AM} & \textit{MD}
				&\textit{AM} & \textit{MD}
				&\textit{AM} & \textit{MD}
				&\textit{AM} & \textit{MD}
				\\\hline
				
				\endhead
				
				\multirow{15}{*}{\ref{fn:f1}}&\multirow{5}{*}{SG}&$|W^{(fc)}_{r,10^2}|$&0&-&0&-&0&-&-&-
				\\*  \cline{3-11}
				&&$|W^{(in)}_{r,10^2}|$&0&-&0&-&0&-&-&-
				\\*  \cline{3-11}
				&&$E_{r,10^2}$&0&-&0&-&0&-&-&-
				\\*  \cline{3-11}
				&&$E'_{r,10^2}$&9.14&-&9.14&-&5.30E+02&-&-&-
				\\*  \cline{3-11}
				&&$\Delta_r$&9.14&-&9.14&-&5.30E+02&-&71.33&-
				\\*  \cline{2-11}
				&\multirow{5}{*}{LHS}&$|W^{(fc)}_{r,10^2}|$&1.40&1&1.40&1&0&0&-&-
				\\*  \cline{3-11}
				&&$|W^{(in)}_{r,10^2}|$&0.40&0&0.40&0&0&0&-&-
				\\*  \cline{3-11}
				&&$E_{r,10^2}$&1.20E+04&8.30E+03&1.20E+04&8.30E+03&0&0&-&-
				\\
\hline

			\end{longtable}
		\end{footnotesize}
	\end{center}
}

{\setlength{\tabcolsep}{4pt}
	\begin{center}
		\begin{footnotesize}
			\begin{longtable}{|*{11}{c|}}

				\caption{
					Number of pole-like points found on the face and inside the domain as well as errors found in all four approximations of functions from \cref{T:TPs}. 
					The interpolation data for these functions have a relative noise level of $\epsilon=10^{-2}$. 
					Number of pole-like points and the error due to poles and not due to poles are given for $r_1(x), r_2(x)$, and $r_3(x)$ only. A  ``-" is used in its place for $p(x)$.
					The number of pole-like points, error due to poles and not due to poles is given for the threshold value of $t=10^2$.
					\textit{Function No.} is the number of the function from \cref{T:TPs}. 
					\textit{Sample Type} is the name of the sampling strategy used to pick interpolation points from the domain. 
					\textit{AM} is the arithmetic mean and \textit{MD} is the median.
					Since SG is deterministic in picking points from the domain, a ``-" is placed for \textit{MD} whenever the sampling strategy is SG.
				} 
				\label{tab:supoleerror10-2}\\
				\cline{4-11}
				\multicolumn{3}{c|}{}&\multicolumn{8}{|c|}{$\epsilon=10^{-2}$}
				\\
				\cline{1-11}
				\multirow{2}{*}{\makecell{\textit{Function} \\ \textit{No.}}} 
				&\multirow{2}{*}{\makecell{\textit{Sample}\\ \textit{Type}}} 
				&\multirow{2}{*}{} 
				&\multicolumn{2}{|c|}{$r_1(x)$}
				&\multicolumn{2}{|c|}{$r_2(x)$}
				&\multicolumn{2}{|c|}{\makecell{$r_3(x)$}}
				&\multicolumn{2}{|c|}{$p(x)$}
				\\
				\cline{4-11}
				&&
				&\textit{AM} & \textit{MD}
				&\textit{AM} & \textit{MD}
				&\textit{AM} & \textit{MD}
				&\textit{AM} & \textit{MD}
				\\\hline
				\endfirsthead
				
				\cline{4-11}
				\multicolumn{3}{c|}{}&\multicolumn{8}{|c|}{$\epsilon=10^{-2}$}
				\\
				\cline{1-11}
				\multirow{2}{*}{\makecell{\textit{Function} \\ \textit{No.}}} 
				&\multirow{2}{*}{\makecell{\textit{Sample}\\ \textit{Type}}} 
				&\multirow{2}{*}{} 
				&\multicolumn{2}{|c|}{$r_1(x)$}
				&\multicolumn{2}{|c|}{$r_2(x)$}
				&\multicolumn{2}{|c|}{\makecell{$r_3(x)$}}
				&\multicolumn{2}{|c|}{$p(x)$}
				\\
				\cline{4-11}
				&&
				&\textit{AM} & \textit{MD}
				&\textit{AM} & \textit{MD}
				&\textit{AM} & \textit{MD}
				&\textit{AM} & \textit{MD}
				\\\hline
				
				\endhead
				
				\multirow{15}{*}{\ref{fn:f1}}&\multirow{5}{*}{SG}&$|W^{(fc)}_{r,10^2}|$&1&-&0&-&0&-&-&-
				\\*  \cline{3-11}
				&&$|W^{(in)}_{r,10^2}|$&0&-&0&-&0&-&-&-
				\\*  \cline{3-11}
				&&$E_{r,10^2}$&2.06E+03&-&0&-&0&-&-&-
				\\*  \cline{3-11}
				&&$E'_{r,10^2}$&1.03E+03&-&1.36E+03&-&5.31E+02&-&-&-
				\\*  \cline{3-11}
				&&$\Delta_r$&2.30E+03&-&1.36E+03&-&5.31E+02&-&73.22&-
				\\*  \cline{2-11}
				&\multirow{5}{*}{LHS}&$|W^{(fc)}_{r,10^2}|$&100&103&1.61E+02&157&0&0&-&-
				\\*  \cline{3-11}
				&&$|W^{(in)}_{r,10^2}|$&30&19&29&25&0&0&-&-
				\\*  \cline{3-11}
				&&$E_{r,10^2}$&1.56E+06&3.73E+05&9.12E+05&5.71E+05&0&0&-&-
				\\
\hline

			\end{longtable}
		\end{footnotesize}
	\end{center}
}

%% file: suppTimesAndIterations.tex
\section{CPU Times and Number of Iterations}\label{app:supptimesiters}

In \cref{fig:cputimeplot} we gave the total CPU time taken by the four approximation approaches for five typical test functions from \cref{T:TPs} when the interpolation data is sampled using d-LHD and is noise-free. 
Additionally, in \cref{fig:iteration} we showed the number of iterations  taken by \cref{A:Polyak} to approximate noise-free data of the same five functions sampled using LHS and d-LHD strategies. 
We showed summary results of the the CPU times and number of iterations over all functions from \cref{T:TPs} in \cref{tab:cputimesummary} and \cref{tab:iterationsummary}, respectively.
In this section, we present detailed results of the CPU time required to compute the four approximations of all functions from \cref{T:TPs} sampled using the three strategies of SG, LHS, D-LHD. Additionally, we also give detailed results of the number of iterations taken by \cref{A:Polyak} to approximate these functions. These results are presented for noise-free data ($\epsilon = 0$) as well as for data with the relative noise of $\epsilon = 10^{-6}$ and $\epsilon = 10^{-2}$.

The timing and iteration results for approximating noise-free interpolation data of all functions from \cref{T:TPs} is shown in \cref{tab:supptimeiter0}. In this table, $p(x)$ is the polynomial approximation, (b) $r_1(x)$ is the rational approximation using \cref{ALG:MVVandQR} without degree reduction, (c) $r_2(x)$ is the rational approximation using \cref{ALG:MVVandQR} with the degree reduction described in \cref{ALG:ReduceDegree}, and (d) $r_3(x)$ is the rational approximation using \cref{A:Polyak}. \textit{Function No.} is the number of the function from \cref{T:TPs}. 
\textit{Sample Type} is the name of the sampling strategy used to pick interpolation points from the domain. 
For $r_3(x)$, the total time is given as its fit time and multistart time. 
The number of iterations taken by \cref{A:Polyak} to obtain $r_3(x)$ is also given.  
$M$ is the arithmetic mean and $SD$ is the standard deviation of the results from running the experiment over LHS and d-LHD samples with different random seeds. Since SG is deterministic in picking points from the domain, a ``-" is placed for \textit{SD} whenever the sampling strategy is SG.
Similarly, the timing and iteration results for approximating noisy data of $\epsilon=10^{-6}$ and $\epsilon=10^{-2}$ is given in \cref{tab:supptimeiter10-6} and \cref{tab:supptimeiter10-2}, respectively.

{\setlength{\tabcolsep}{4pt}
	\begin{center}
		\begin{footnotesize}
			\begin{longtable}{|*{14}{c|}}

					\caption{
						CPU times for all four approximations of functions from \cref{T:TPs}. The interpolation data for these functions are noise-free, i.e., $\epsilon=0$. 
						\textit{Function No.} is the number of the function from \cref{T:TPs}. 
						\textit{Sample Type} is the name of the sampling strategy used to pick interpolation points from the domain. 
						For $r_3(x)$, the total time is given as its fit time and multistart time. 
						The number of iterations taken by \cref{A:Polyak} to obtain $r_3(x)$ is also given.  
						$M$ and $SD$ are the average and standard deviation, respectively. 
						Since SG is deterministic in picking points from the domain, a ``-" is placed for \textit{SD} whenever the sampling strategy is SG.
					} 
					\label{tab:supptimeiter0}\\
					\cline{3-14}
					\multicolumn{2}{c|}{}&\multicolumn{12}{|c|}{$\epsilon=0$}
					\\
					\cline{1-14}
					\multirow{2}{*}{\makecell{\textit{Function} \\ \textit{No.}}} 
					&\multirow{2}{*}{\makecell{\textit{Sample}\\ \textit{Type}}} 
					&\multicolumn{2}{|c|}{$r_1(x)$}&\multicolumn{2}{|c|}{$r_2(x)$}
					&\multicolumn{2}{|c|}{\makecell{$r_3(x)$ \\\textit{fit time}}}
					& \multicolumn{2}{|c|}{\makecell{$r_3(x)$ \\\textit{ms time}}}
					& \multicolumn{2}{|c|}{\makecell{$r_3(x)$\\ \textit{iterations}}}
					&\multicolumn{2}{|c|}{$p(x)$}
					\\
					\cline{3-14}
					 &
					&$M$ & $SD$
					&$M$ & $SD$
					&$M$ & $SD$
					&$M$ & $SD$
					&$M$ & $SD$
					&$M$ & $SD$
					\\\hline
					\endfirsthead
					
					\cline{3-14}
					\multicolumn{2}{c|}{}&\multicolumn{12}{|c|}{$\epsilon=0$}
					\\
					\cline{1-14}
					\multirow{2}{*}{\makecell{\textit{Function} \\ \textit{No.}}} 
					&\multirow{2}{*}{\makecell{\textit{Sample}\\ \textit{Type}}} 
					&\multicolumn{2}{|c|}{$r_1(x)$}&\multicolumn{2}{|c|}{$r_2(x)$}
					&\multicolumn{2}{|c|}{\makecell{$r_3(x)$ \\\textit{fit time}}}
					& \multicolumn{2}{|c|}{\makecell{$r_3(x)$ \\\textit{ms time}}}
					& \multicolumn{2}{|c|}{\makecell{$r_3(x)$\\ \textit{iterations}}}
					&\multicolumn{2}{|c|}{$p(x)$}
					\\
					\cline{3-14}
					&
					&$M$ & $SD$
					&$M$ & $SD$
					&$M$ & $SD$
					&$M$ & $SD$
					&$M$ & $SD$
					&$M$ & $SD$
					\\\hline
					
					\endhead
					
					\multirow{3}{*}{\ref{fn:f1}}&SG&2.12&-&2.55&-&0.04&-&15.85&-&1&-&0.41&-
					\\*  \cline{2-14}
					&LHS&3.27&0.12&3.59&0.06&0.05&9.86E-04&15.85&2.81E-04&1&0&0.93&0.06
					\\*  \cline{2-14}
					&d-LHD&3.24&0.11&4.02&0.03&0.05&0.02&15.85&1.19E-04&1&0&0.90&0.02
					\\ \hline
					\multirow{3}{*}{\ref{fn:f2}}&SG&2.25&-&2.55&-&0.05&-&15.85&-&1&-&0.52&-
					\\*  \cline{2-14}
					&LHS&3.54&0.16&3.69&0.13&0.04&1.39E-03&15.85&1.92E-04&1&0&0.90&0.01
					\\*  \cline{2-14}
					&d-LHD&3.36&0.13&4.02&0.02&0.04&7.90E-03&15.85&1.24E-04&1&0&0.93&0.03
					\\ \hline
					\multirow{3}{*}{\ref{fn:f3}}&SG&2.09&-&2.56&-&0.04&-&15.85&-&1&-&0.42&-
					\\*  \cline{2-14}
					&LHS&3.67&0.10&3.98&0.03&0.04&2.21E-03&15.85&2.48E-04&1&0&0.93&0.02
					\\*  \cline{2-14}
					&d-LHD&3.51&0.07&4.04&0.07&0.03&1.95E-03&15.85&2.40E-04&1&0&0.91&0.02
					\\ \hline
					\multirow{3}{*}{\ref{fn:f4}}&SG&2.07&-&2.55&-&0.04&-&15.85&-&1&-&0.42&-
					\\*  \cline{2-14}
					&LHS&3.77&0.05&4.04&0.04&0.06&4.17E-03&15.85&1.03E-04&2&0&0.90&0.04
					\\*  \cline{2-14}
					&d-LHD&3.75&0.08&4.05&0.09&0.03&4.01E-03&15.85&5.91E-05&1&0&0.92&0.02
					\\ \hline
					\multirow{3}{*}{\ref{fn:f5}}&SG&2.07&-&2.55&-&0.04&-&15.85&-&1&-&0.52&-
					\\*  \cline{2-14}
					&LHS&3.80&0.02&4.17&0.15&0.03&1.34E-03&15.85&4.45E-04&1&0&0.86&0.01
					\\*  \cline{2-14}
					&d-LHD&3.80&0.04&4.17&0.09&0.04&6.54E-03&15.85&1.05E-04&1&0&0.89&0.02
					\\ \hline
					\multirow{3}{*}{\ref{fn:f7}}&SG&2.07&-&2.56&-&0.04&-&15.85&-&1&-&0.52&-
					\\*  \cline{2-14}
					&LHS&3.87&0.02&3.73&0.63&0.03&1.02E-03&15.85&2.84E-04&1&0&0.92&0.02
					\\*  \cline{2-14}
					&d-LHD&3.82&0.03&4.12&0.03&0.03&6.44E-04&15.85&3.21E-04&1&0&0.86&0.05
					\\ \hline
					\multirow{3}{*}{\ref{fn:f8}}&SG&2.25&-&2.55&-&0.04&-&15.85&-&1&-&0.42&-
					\\*  \cline{2-14}
					&LHS&2.46&0.14&3.32&0.52&0.38&0.27&15.86&7.98E-03&4&2.28&0.70&0.07
					\\*  \cline{2-14}
					&d-LHD&2.93&0.61&4.03&0.06&0.69&0.37&15.86&9.83E-03&1.80&0.75&0.71&0.05
					\\ \hline
					\multirow{3}{*}{\ref{fn:f9}}&SG&2.24&-&2.55&-&0.04&-&15.85&-&1&-&0.42&-
					\\*  \cline{2-14}
					&LHS&2.36&0.02&3.93&0.06&0.10&0.05&15.85&7.27E-04&2.20&0.75&0.68&8.45E-03
					\\*  \cline{2-14}
					&d-LHD&2.40&0.02&3.87&0.37&0.10&0.06&15.85&6.48E-05&1&0&0.67&0.04
					\\ \hline
					\multirow{3}{*}{\ref{fn:f10}}&SG&7.79&-&7.42&-&3.00&-&3.73E+02&-&1&-&0.56&-
					\\*  \cline{2-14}
					&LHS&3.09&0.03&4.79&0.43&4.68&2.43&3.73E+02&0.04&1.60&1.20&0.67&0.05
					\\*  \cline{2-14}
					&d-LHD&3.16&0.02&4.91&0.01&9.32&8.35&3.73E+02&0.01&1.20&0.40&0.71&5.57E-03
					\\ \hline
					\multirow{3}{*}{\ref{fn:f12}}&SG&2.26&-&2.55&-&0.06&-&15.85&-&1&-&0.52&-
					\\*  \cline{2-14}
					&LHS&2.30&0.01&2.31&8.92E-03&0.02&2.64E-03&15.85&1.03E-04&1&0&0.64&8.47E-03
					\\*  \cline{2-14}
					&d-LHD&2.36&0.04&3.01&0.79&0.03&3.07E-03&15.85&5.89E-03&1&0&0.63&0.05
					\\ \hline
					\multirow{3}{*}{\ref{fn:f13}}&SG&2.07&-&2.56&-&0.07&-&15.85&-&1&-&0.52&-
					\\*  \cline{2-14}
					&LHS&2.38&0.29&2.86&0.46&0.02&4.45E-03&15.85&1.37E-04&1&0&0.61&0.04
					\\*  \cline{2-14}
					&d-LHD&2.33&0.01&2.26&0.04&0.03&5.63E-03&15.85&3.18E-04&1&0&0.62&0.02
					\\ \hline
					\multirow{3}{*}{\ref{fn:f14}}&SG&2.07&-&2.55&-&0.04&-&15.85&-&1&-&0.52&-
					\\*  \cline{2-14}
					&LHS&2.25&6.90E-03&3.22&0.07&0.03&7.98E-03&15.85&3.05E-04&1.20&0.40&0.61&0.01
					\\*  \cline{2-14}
					&d-LHD&2.24&0.01&2.64&0.71&0.04&1.89E-03&15.85&1.53E-04&1&0&0.58&0.05
					\\ \hline
					\multirow{3}{*}{\ref{fn:f15}}&SG&2.08&-&2.56&-&0.04&-&15.85&-&1&-&0.52&-
					\\*  \cline{2-14}
					&LHS&2.23&7.86E-03&3.53&0.21&0.03&3.38E-03&15.85&1.22E-04&1&0&0.59&0.01
					\\*  \cline{2-14}
					&d-LHD&2.35&0.25&4.00&0.02&0.04&1.52E-03&15.85&1.86E-04&1&0&0.52&0.04
					\\ \hline
					\multirow{3}{*}{\ref{fn:f16}}&SG&2.26&-&2.56&-&0.04&-&15.85&-&1&-&0.41&-
					\\*  \cline{2-14}
					&LHS&2.53&0.30&3.54&0.02&0.03&0.01&15.85&5.00E-04&1.20&0.40&0.58&3.20E-03
					\\*  \cline{2-14}
					&d-LHD&2.50&0.19&4.02&0.02&0.04&0.01&15.85&3.38E-04&1&0&0.48&0.05
					\\ \hline
					\multirow{3}{*}{\ref{fn:f17}}&SG&2.50&-&2.69&-&0.39&-&44.94&-&1&-&0.42&-
					\\*  \cline{2-14}
					&LHS&3.06&0.12&3.97&0.13&1.26&0.68&45.08&0.15&5.80&3.54&0.53&0.06
					\\*  \cline{2-14}
					&d-LHD&2.86&0.04&4.10&0.03&0.56&0.21&45.52&0.68&2.40&1.02&0.45&0.04
					\\ \hline
					\multirow{3}{*}{\ref{fn:f18}}&SG&7.89&-&9.87&-&2.03E+02&-&3.73E+02&-&1&-&0.57&-
					\\*  \cline{2-14}
					&LHS&4.01&0.06&4.72&0.04&52.04&6.40&3.73E+02&0.10&27.60&3.77&0.50&0.02
					\\*  \cline{2-14}
					&d-LHD&3.84&0.03&4.94&0.18&21.81&3.65&3.73E+02&0.27&9.20&1.17&0.49&0.06
					\\ \hline
					\multirow{3}{*}{\ref{fn:f19}}&SG&7.76&-&9.11&-&4.01&-&3.73E+02&-&1&-&0.45&-
					\\*  \cline{2-14}
					&LHS&4.18&0.03&4.79&0.29&2.71&0.36&3.73E+02&6.96E-04&1&0&0.48&9.14E-03
					\\*  \cline{2-14}
					&d-LHD&4.00&0.13&5.07&0.08&2.43&0.20&3.73E+02&2.79E-03&1&0&0.49&0.05
					\\ \hline
					\multirow{3}{*}{\ref{fn:f20}}&SG&7.86&-&9.69&-&7.28&-&3.73E+02&-&1&-&0.46&-
					\\*  \cline{2-14}
					&LHS&4.33&0.16&4.50&0.58&45.16&12.52&3.74E+02&0.49&26.40&8.01&0.48&7.19E-03
					\\*  \cline{2-14}
					&d-LHD&4.00&0.35&4.82&0.43&1.36E+02&28.80&3.78E+02&2.29&80.40&30.55&0.46&0.02
					\\ \hline
					\multirow{3}{*}{\ref{fn:f21}}&SG&1.79&-&2.55&-&0.04&-&15.85&-&1&-&0.51&-
					\\*  \cline{2-14}
					&LHS&3.50&0.03&3.58&0.22&0.04&0.01&15.86&0.02&1.20&0.40&0.55&0.03
					\\*  \cline{2-14}
					&d-LHD&3.43&0.03&4.00&0.03&0.03&1.33E-03&15.85&6.93E-04&1&0&0.42&0.01
					\\ \hline
					\multirow{3}{*}{\ref{fn:f22}}&SG&1.99&-&3.64&-&0.05&-&15.85&-&1&-&0.50&-
					\\*  \cline{2-14}
					&LHS&3.62&0.04&3.60&0.10&0.04&0.01&15.85&2.13E-03&1.20&0.40&0.64&0.04
					\\*  \cline{2-14}
					&d-LHD&3.52&0.04&4.00&0.05&0.16&0.12&15.85&2.67E-04&1&0&0.50&0.03
					\\ \hline					
					
			\end{longtable}
		\end{footnotesize}
	\end{center}
}

{\setlength{\tabcolsep}{4pt}
	\begin{center}
		\begin{footnotesize}
			\begin{longtable}{|*{14}{c|}}

				\caption{
					CPU times for all four approximations of functions from \cref{T:TPs}. The relative noise level in the interpolation data  is $\epsilon=10^{-6}$.
					\textit{Function No.} is the number of the function from \cref{T:TPs}. 
					\textit{Sample Type} is the name of the sampling strategy used to pick interpolation points from the domain. 
					For $r_3(x)$, the total time is given as its fit time and multistart time. 
					The number of iterations taken by \cref{A:Polyak} to obtain $r_3(x)$ is also given.  
					$M$ and $SD$ are the average and standard deviation, respectively. 
					Since SG is deterministic in picking points from the domain, a ``-" is placed for \textit{SD} whenever the sampling strategy is SG.
				} 
				\label{tab:supptimeiter10-6}\\
				\cline{3-14}
				\multicolumn{2}{c|}{}&\multicolumn{12}{|c|}{$\epsilon=10^{-6}$}
				\\
				\cline{1-14}
				\multirow{2}{*}{\makecell{\textit{Function} \\ \textit{No.}}} 
				&\multirow{2}{*}{\makecell{\textit{Sample}\\ \textit{Type}}} 
				&\multicolumn{2}{|c|}{$r_1(x)$}&\multicolumn{2}{|c|}{$r_2(x)$}
				&\multicolumn{2}{|c|}{\makecell{$r_3(x)$ \\\textit{fit time}}}
				& \multicolumn{2}{|c|}{\makecell{$r_3(x)$ \\\textit{ms time}}}
				& \multicolumn{2}{|c|}{\makecell{$r_3(x)$\\ \textit{iterations}}}
				&\multicolumn{2}{|c|}{$p(x)$}
				\\
				\cline{3-14}
				&
				&$M$ & $SD$
				&$M$ & $SD$
				&$M$ & $SD$
				&$M$ & $SD$
				&$M$ & $SD$
				&$M$ & $SD$
				\\\hline
				\endfirsthead
				
				\cline{3-14}
				\multicolumn{2}{c|}{}&\multicolumn{12}{|c|}{$\epsilon=10^{-6}$}
				\\
				\cline{1-14}
				\multirow{2}{*}{\makecell{\textit{Function} \\ \textit{No.}}} 
				&\multirow{2}{*}{\makecell{\textit{Sample}\\ \textit{Type}}} 
				&\multicolumn{2}{|c|}{$r_1(x)$}&\multicolumn{2}{|c|}{$r_2(x)$}
				&\multicolumn{2}{|c|}{\makecell{$r_3(x)$ \\\textit{fit time}}}
				& \multicolumn{2}{|c|}{\makecell{$r_3(x)$ \\\textit{ms time}}}
				& \multicolumn{2}{|c|}{\makecell{$r_3(x)$\\ \textit{iterations}}}
				&\multicolumn{2}{|c|}{$p(x)$}
				\\
				\cline{3-14}
				&
				&$M$ & $SD$
				&$M$ & $SD$
				&$M$ & $SD$
				&$M$ & $SD$
				&$M$ & $SD$
				&$M$ & $SD$
				\\\hline
				
				\endhead
				
				\multirow{3}{*}{\ref{fn:f1}}&SG&2.25&-&2.07&-&0.03&-&15.85&-&1&-&0.52&-
				\\*  \cline{2-14}
				&LHS&3.28&0.16&3.34&0.13&0.04&2.97E-03&15.85&4.32E-04&1&0&0.91&0.01
				\\*  \cline{2-14}
				&d-LHD&3.35&0.04&3.45&0.12&0.05&0.02&15.85&2.23E-04&1&0&1.19&0.05
				\\ \hline
				\multirow{3}{*}{\ref{fn:f2}}&SG&2.07&-&2.24&-&0.07&-&15.85&-&1&-&0.53&-
				\\*  \cline{2-14}
				&LHS&3.55&0.11&3.59&0.08&0.04&6.98E-04&15.87&0.04&1&0&0.90&0.02
				\\*  \cline{2-14}
				&d-LHD&3.57&0.10&3.67&0.17&0.04&6.35E-03&15.85&8.69E-05&1&0&1.18&0.03
				\\ \hline
				\multirow{3}{*}{\ref{fn:f3}}&SG&2.24&-&2.07&-&0.03&-&15.85&-&1&-&0.52&-
				\\*  \cline{2-14}
				&LHS&3.68&0.03&3.67&0.06&0.04&1.75E-03&15.85&4.57E-04&1&0&0.83&6.95E-03
				\\*  \cline{2-14}
				&d-LHD&3.60&0.04&3.65&0.04&0.04&3.16E-03&15.85&5.49E-04&1&0&1.19&0.04
				\\ \hline
				\multirow{3}{*}{\ref{fn:f4}}&SG&2.07&-&2.08&-&0.04&-&15.85&-&1&-&0.41&-
				\\*  \cline{2-14}
				&LHS&3.75&0.05&3.89&0.11&0.06&0.01&15.85&3.07E-03&2.20&0.40&0.80&0.01
				\\*  \cline{2-14}
				&d-LHD&3.74&0.06&3.74&0.11&0.03&1.22E-03&15.85&1.51E-04&1&0&1.15&0.13
				\\ \hline
				\multirow{3}{*}{\ref{fn:f5}}&SG&2.07&-&2.26&-&0.04&-&15.85&-&1&-&0.41&-
				\\*  \cline{2-14}
				&LHS&3.75&0.01&4.05&0.08&0.04&5.35E-03&15.85&3.36E-04&1&0&0.77&2.08E-03
				\\*  \cline{2-14}
				&d-LHD&3.70&0.02&4.05&0.12&0.05&0.02&15.85&7.20E-04&1&0&1.19&0.02
				\\ \hline
				\multirow{3}{*}{\ref{fn:f7}}&SG&2.25&-&2.08&-&0.03&-&15.85&-&1&-&0.53&-
				\\*  \cline{2-14}
				&LHS&3.79&0.02&3.92&0.06&0.03&1.43E-03&15.85&3.54E-04&1&0&0.73&0.08
				\\*  \cline{2-14}
				&d-LHD&3.76&0.02&4.07&0.10&0.03&2.91E-03&15.85&4.94E-04&1&0&1.05&0.11
				\\ \hline
				\multirow{3}{*}{\ref{fn:f8}}&SG&2.23&-&2.06&-&0.03&-&15.85&-&1&-&0.40&-
				\\*  \cline{2-14}
				&LHS&2.68&0.24&3.81&0.03&0.17&0.16&15.86&0.02&2.40&1.02&0.58&0.05
				\\*  \cline{2-14}
				&d-LHD&2.48&0.06&3.98&0.17&2.22&1.74&15.86&9.00E-03&2.40&1.02&0.94&9.60E-03
				\\ \hline
				\multirow{3}{*}{\ref{fn:f9}}&SG&2.06&-&2.07&-&0.03&-&15.85&-&1&-&0.41&-
				\\*  \cline{2-14}
				&LHS&2.33&0.07&3.26&0.65&0.23&0.36&15.85&2.15E-03&1.80&1.17&0.55&3.25E-03
				\\*  \cline{2-14}
				&d-LHD&2.37&0.06&3.24&0.63&0.06&0.03&15.85&6.15E-05&1&0&0.92&0.05
				\\ \hline
				\multirow{3}{*}{\ref{fn:f10}}&SG&7.30&-&7.55&-&6.28&-&3.73E+02&-&1&-&0.57&-
				\\*  \cline{2-14}
				&LHS&3.10&0.02&3.40&0.13&5.45&1.86&3.73E+02&0.04&1.20&0.40&0.57&6.26E-03
				\\*  \cline{2-14}
				&d-LHD&3.14&0.02&3.50&0.26&45.66&86.27&3.73E+02&8.23E-04&1&0&0.95&0.02
				\\ \hline
				\multirow{3}{*}{\ref{fn:f12}}&SG&2.07&-&2.26&-&0.05&-&15.85&-&1&-&0.53&-
				\\*  \cline{2-14}
				&LHS&2.30&9.61E-03&2.43&0.28&0.03&4.17E-03&15.85&1.38E-04&1&0&0.57&0.07
				\\*  \cline{2-14}
				&d-LHD&2.30&0.05&2.28&0.06&0.03&7.01E-03&15.85&1.19E-04&1&0&0.94&0.01
				\\ \hline
				\multirow{3}{*}{\ref{fn:f13}}&SG&2.22&-&2.08&-&0.08&-&15.85&-&1&-&0.40&-
				\\*  \cline{2-14}
				&LHS&2.26&0.02&2.26&0.01&0.03&3.98E-03&15.87&0.04&1&0&0.55&4.11E-03
				\\*  \cline{2-14}
				&d-LHD&2.26&0.02&2.38&0.34&0.05&0.04&15.85&1.67E-04&1&0&0.87&0.17
				\\ \hline
				\multirow{3}{*}{\ref{fn:f14}}&SG&2.24&-&2.07&-&0.05&-&15.85&-&1&-&0.41&-
				\\*  \cline{2-14}
				&LHS&2.23&0.01&2.31&0.24&0.04&0.01&15.85&2.45E-04&1.20&0.40&0.53&0.05
				\\*  \cline{2-14}
				&d-LHD&2.20&0.06&2.18&0.07&0.03&7.10E-03&15.85&6.37E-05&1&0&0.95&0.04
				\\ \hline
				\multirow{3}{*}{\ref{fn:f15}}&SG&2.07&-&2.07&-&0.07&-&15.85&-&1&-&0.41&-
				\\*  \cline{2-14}
				&LHS&2.19&0.06&3.00&0.06&0.04&4.29E-03&15.85&1.74E-04&1&0&0.54&8.68E-03
				\\*  \cline{2-14}
				&d-LHD&2.36&0.27&2.39&0.28&0.04&5.36E-03&15.85&1.81E-04&1&0&0.94&7.95E-03
				\\ \hline
				\multirow{3}{*}{\ref{fn:f16}}&SG&2.07&-&2.28&-&0.06&-&15.85&-&1&-&0.41&-
				\\*  \cline{2-14}
				&LHS&2.30&0.24&2.97&0.28&0.04&0.01&15.85&5.76E-04&1.20&0.40&0.54&8.47E-03
				\\*  \cline{2-14}
				&d-LHD&2.42&0.40&2.98&0.06&0.04&0.01&15.85&6.09E-05&1&0&0.93&9.30E-03
				\\ \hline
				\multirow{3}{*}{\ref{fn:f17}}&SG&2.49&-&2.64&-&0.90&-&44.94&-&1&-&0.41&-
				\\*  \cline{2-14}
				&LHS&3.08&0.17&3.11&0.32&1.40&1.09&45.04&0.09&6.80&5.67&0.55&0.03
				\\*  \cline{2-14}
				&d-LHD&3.04&0.03&3.14&0.08&0.49&0.23&44.98&0.05&2.20&0.98&1.12&0.15
				\\ \hline
				\multirow{3}{*}{\ref{fn:f18}}&SG&7.44&-&6.81&-&1.56E+02&-&3.73E+02&-&1&-&0.45&-
				\\*  \cline{2-14}
				&LHS&4.05&0.08&4.08&0.67&52.50&9.45&3.73E+02&0.20&30&5.25&0.54&0.02
				\\*  \cline{2-14}
				&d-LHD&4.02&0.06&3.94&0.03&18.10&2.65&3.73E+02&0.12&9.80&1.60&1.21&0.03
				\\ \hline
				\multirow{3}{*}{\ref{fn:f19}}&SG&7.46&-&7.59&-&3.28&-&3.73E+02&-&1&-&0.44&-
				\\*  \cline{2-14}
				&LHS&4.17&0.03&4.59&0.20&2.72&0.37&3.73E+02&1.43E-03&1&0&0.54&4.34E-03
				\\*  \cline{2-14}
				&d-LHD&4.08&0.05&4.55&0.37&2.53&0.19&3.73E+02&1.35E-03&1&0&1.23&0.05
				\\ \hline
				\multirow{3}{*}{\ref{fn:f20}}&SG&6.26&-&8.46&-&2.83E+02&-&3.73E+02&-&1&-&0.40&-
				\\*  \cline{2-14}
				&LHS&4.15&0.08&4.63&0.11&41.62&12.70&3.74E+02&0.53&26&7.62&0.59&0.12
				\\*  \cline{2-14}
				&d-LHD&4.11&0.08&4.53&0.07&1.00E+02&39.12&3.78E+02&2.39&79.40&30.93&1.22&0.04
				\\ \hline
				\multirow{3}{*}{\ref{fn:f21}}&SG&1.99&-&2.04&-&0.05&-&15.85&-&1&-&0.50&-
				\\*  \cline{2-14}
				&LHS&3.40&0.02&3.53&0.10&0.04&0.01&15.85&5.81E-03&1.20&0.40&0.76&3.22E-03
				\\*  \cline{2-14}
				&d-LHD&3.34&0.03&3.66&0.22&0.03&3.07E-03&15.85&4.13E-04&1&0&1.18&0.05
				\\ \hline
				\multirow{3}{*}{\ref{fn:f22}}&SG&1.74&-&1.81&-&0.06&-&15.85&-&1&-&0.34&-
				\\*  \cline{2-14}
				&LHS&3.49&0.04&3.64&0.02&0.04&0.01&15.85&8.86E-03&1.20&0.40&0.77&4.33E-03
				\\*  \cline{2-14}
				&d-LHD&3.47&0.03&3.44&0.02&0.08&0.05&15.85&2.33E-04&1&0&1.18&6.01E-03
				\\ \hline

			\end{longtable}
		\end{footnotesize}
	\end{center}
}

{\setlength{\tabcolsep}{4pt}
	\begin{center}
		\begin{footnotesize}
			\begin{longtable}{|*{14}{c|}}

				\caption{
					CPU times for all four approximations of functions from \cref{T:TPs}. The relative noise level in the interpolation data  is $\epsilon=10^{-2}$.
					\textit{Function No.} is the number of the function from \cref{T:TPs}. 
					\textit{Sample Type} is the name of the sampling strategy used to pick interpolation points from the domain. 
					For $r_3(x)$, the total time is given as its fit time and multistart time. 
					The number of iterations taken by \cref{A:Polyak} to obtain $r_3(x)$ is also given.  
					$M$ and $SD$ are the average and standard deviation, respectively. 
					Since SG is deterministic in picking points from the domain, a ``-" is placed for \textit{SD} whenever the sampling strategy is SG.
				} 
				\label{tab:supptimeiter10-2}\\
				\cline{3-14}
				\multicolumn{2}{c|}{}&\multicolumn{12}{|c|}{$\epsilon=10^{-2}$}
				\\
				\cline{1-14}
				\multirow{2}{*}{\makecell{\textit{Function} \\ \textit{No.}}} 
				&\multirow{2}{*}{\makecell{\textit{Sample}\\ \textit{Type}}} 
				&\multicolumn{2}{|c|}{$r_1(x)$}&\multicolumn{2}{|c|}{$r_2(x)$}
				&\multicolumn{2}{|c|}{\makecell{$r_3(x)$ \\\textit{fit time}}}
				& \multicolumn{2}{|c|}{\makecell{$r_3(x)$ \\\textit{ms time}}}
				& \multicolumn{2}{|c|}{\makecell{$r_3(x)$\\ \textit{iterations}}}
				&\multicolumn{2}{|c|}{$p(x)$}
				\\
				\cline{3-14}
				&
				&$M$ & $SD$
				&$M$ & $SD$
				&$M$ & $SD$
				&$M$ & $SD$
				&$M$ & $SD$
				&$M$ & $SD$
				\\\hline
				\endfirsthead
				
				\cline{3-14}
				\multicolumn{2}{c|}{}&\multicolumn{12}{|c|}{$\epsilon=10^{-2}$}
				\\
				\cline{1-14}
				\multirow{2}{*}{\makecell{\textit{Function} \\ \textit{No.}}} 
				&\multirow{2}{*}{\makecell{\textit{Sample}\\ \textit{Type}}} 
				&\multicolumn{2}{|c|}{$r_1(x)$}&\multicolumn{2}{|c|}{$r_2(x)$}
				&\multicolumn{2}{|c|}{\makecell{$r_3(x)$ \\\textit{fit time}}}
				& \multicolumn{2}{|c|}{\makecell{$r_3(x)$ \\\textit{ms time}}}
				& \multicolumn{2}{|c|}{\makecell{$r_3(x)$\\ \textit{iterations}}}
				&\multicolumn{2}{|c|}{$p(x)$}
				\\
				\cline{3-14}
				&
				&$M$ & $SD$
				&$M$ & $SD$
				&$M$ & $SD$
				&$M$ & $SD$
				&$M$ & $SD$
				&$M$ & $SD$
				\\\hline
				
				\endhead
				
				\multirow{3}{*}{\ref{fn:f1}}&SG&2.22&-&1.62&-&0.03&-&15.85&-&1&-&0.51&-
				\\*  \cline{2-14}
				&LHS&3.35&0.09&3.76&0.03&0.07&0.03&15.85&3.72E-03&2&0.89&0.85&1.89E-03
				\\*  \cline{2-14}
				&d-LHD&3.23&0.25&3.99&3.57E-03&0.05&0.02&15.85&2.65E-04&1&0&0.91&6.43E-03
				\\ \hline
				\multirow{3}{*}{\ref{fn:f2}}&SG&2.05&-&2.24&-&0.04&-&15.85&-&1&-&0.51&-
				\\*  \cline{2-14}
				&LHS&3.64&0.08&3.86&0.05&0.04&2.23E-03&15.85&4.81E-04&1&0&0.84&4.30E-03
				\\*  \cline{2-14}
				&d-LHD&3.65&0.12&4.00&0.04&0.04&8.95E-03&15.85&4.46E-04&1&0&0.93&0.04
				\\ \hline
				\multirow{3}{*}{\ref{fn:f3}}&SG&2.06&-&2.25&-&0.03&-&15.85&-&1&-&0.51&-
				\\*  \cline{2-14}
				&LHS&3.65&0.07&3.99&0.11&0.09&0.01&15.86&4.93E-03&2.60&0.49&0.85&8.79E-03
				\\*  \cline{2-14}
				&d-LHD&3.72&0.04&4.00&0.01&0.04&9.88E-03&15.85&7.30E-04&1&0&0.95&0.06
				\\ \hline
				\multirow{3}{*}{\ref{fn:f4}}&SG&2.06&-&2.25&-&0.04&-&15.85&-&1&-&0.51&-
				\\*  \cline{2-14}
				&LHS&3.79&0.04&4.05&0.04&0.03&2.62E-03&15.85&7.38E-05&1&0&0.85&0.02
				\\*  \cline{2-14}
				&d-LHD&3.78&0.06&3.96&0.07&0.04&2.03E-03&15.86&0.02&1&0&0.91&0.08
				\\ \hline
				\multirow{3}{*}{\ref{fn:f5}}&SG&2.07&-&2.24&-&0.03&-&15.85&-&1&-&0.51&-
				\\*  \cline{2-14}
				&LHS&3.78&0.02&4.05&0.05&0.04&2.31E-03&15.85&4.10E-04&1&0&0.83&0.01
				\\*  \cline{2-14}
				&d-LHD&3.74&0.03&4.09&0.16&0.06&0.02&15.85&4.24E-04&1&0&0.84&0.02
				\\ \hline
				\multirow{3}{*}{\ref{fn:f7}}&SG&2.05&-&2.25&-&0.04&-&15.85&-&1&-&0.51&-
				\\*  \cline{2-14}
				&LHS&3.84&0.01&4.04&0.08&0.04&2.68E-03&15.85&2.23E-04&1&0&0.79&0.08
				\\*  \cline{2-14}
				&d-LHD&3.77&5.32E-03&3.93&0.10&0.06&0.03&15.85&1.98E-04&1&0&0.83&0.10
				\\ \hline
				\multirow{3}{*}{\ref{fn:f8}}&SG&2.23&-&2.11&-&0.04&-&15.85&-&1&-&0.51&-
				\\*  \cline{2-14}
				&LHS&2.70&0.30&3.98&0.05&0.08&0.01&15.86&7.47E-03&2.20&0.40&0.67&0.06
				\\*  \cline{2-14}
				&d-LHD&2.66&0.30&3.25&0.65&0.05&0.02&15.85&0.01&1.20&0.40&0.72&9.13E-03
				\\ \hline
				\multirow{3}{*}{\ref{fn:f9}}&SG&2.24&-&2.24&-&0.04&-&15.85&-&1&-&0.51&-
				\\*  \cline{2-14}
				&LHS&2.36&7.28E-03&4.02&0.02&0.14&0.04&15.86&9.77E-03&3.80&1.17&0.68&0.04
				\\*  \cline{2-14}
				&d-LHD&2.34&0.06&3.34&0.57&0.05&6.60E-03&15.86&0.02&1&0&0.63&0.05
				\\ \hline
				\multirow{3}{*}{\ref{fn:f10}}&SG&7.19&-&7.89&-&4.10&-&3.73E+02&-&1&-&0.56&-
				\\*  \cline{2-14}
				&LHS&3.09&0.01&5.08&0.04&8.71&3.56&3.73E+02&0.17&3.40&1.02&0.68&0.05
				\\*  \cline{2-14}
				&d-LHD&3.10&8.19E-03&4.38&0.50&1.84&0.25&3.73E+02&1.98E-03&1&0&0.65&0.04
				\\ \hline
				\multirow{3}{*}{\ref{fn:f12}}&SG&2.06&-&2.24&-&0.06&-&15.85&-&1&-&0.51&-
				\\*  \cline{2-14}
				&LHS&2.28&0.01&2.66&0.53&0.03&3.30E-03&15.85&3.44E-04&1&0&0.63&0.05
				\\*  \cline{2-14}
				&d-LHD&2.33&8.95E-03&2.43&0.32&0.04&7.34E-03&15.85&3.09E-04&1&0&0.61&0.02
				\\ \hline
				\multirow{3}{*}{\ref{fn:f13}}&SG&2.06&-&2.24&-&0.08&-&15.85&-&1&-&0.51&-
				\\*  \cline{2-14}
				&LHS&2.25&0.02&2.25&8.32E-03&0.04&3.71E-03&15.85&5.09E-04&1&0&0.61&0.05
				\\*  \cline{2-14}
				&d-LHD&2.30&0.02&2.09&8.19E-03&0.07&0.03&15.85&3.56E-04&1&0&0.54&0.03
				\\ \hline
				\multirow{3}{*}{\ref{fn:f14}}&SG&2.06&-&2.23&-&0.05&-&15.85&-&1&-&0.51&-
				\\*  \cline{2-14}
				&LHS&2.23&7.00E-03&3.25&0.02&0.04&1.74E-03&15.85&1.53E-04&1&0&0.63&0.05
				\\*  \cline{2-14}
				&d-LHD&2.23&8.13E-03&2.43&0.71&0.06&0.03&15.85&8.47E-05&1&0&0.57&0.05
				\\ \hline
				\multirow{3}{*}{\ref{fn:f15}}&SG&2.07&-&2.23&-&0.05&-&15.85&-&1&-&0.51&-
				\\*  \cline{2-14}
				&LHS&2.23&6.88E-03&3.15&0.39&0.03&2.41E-03&15.85&1.83E-04&1&0&0.69&0.06
				\\*  \cline{2-14}
				&d-LHD&2.22&8.53E-03&3.86&0.02&0.07&0.04&15.85&1.79E-04&1&0&0.54&1.97E-03
				\\ \hline
				\multirow{3}{*}{\ref{fn:f16}}&SG&2.06&-&2.24&-&0.05&-&15.85&-&1&-&0.51&-
				\\*  \cline{2-14}
				&LHS&2.38&0.29&3.42&0.02&0.04&6.43E-03&15.85&2.58E-04&1&0&0.64&0.05
				\\*  \cline{2-14}
				&d-LHD&2.23&0.13&3.85&6.24E-03&0.03&7.63E-03&15.85&2.71E-04&1&0&0.54&0.01
				\\ \hline
				\multirow{3}{*}{\ref{fn:f17}}&SG&2.45&-&2.66&-&0.76&-&44.94&-&1&-&0.51&-
				\\*  \cline{2-14}
				&LHS&3.02&0.07&3.19&0.69&1.02&0.54&45.07&0.12&4&2.10&0.85&4.82E-03
				\\*  \cline{2-14}
				&d-LHD&3.06&0.05&3.89&0.16&0.33&0.10&44.94&3.64E-03&1.20&0.40&0.49&0.04
				\\ \hline
				\multirow{3}{*}{\ref{fn:f18}}&SG&6.97&-&8.20&-&2.90E+02&-&3.73E+02&-&1&-&0.56&-
				\\*  \cline{2-14}
				&LHS&4.05&0.06&4.70&0.18&63.03&13.12&3.74E+02&0.57&34&7.51&0.89&0.02
				\\*  \cline{2-14}
				&d-LHD&4.12&0.04&4.56&0.64&17.07&2.42&3.73E+02&0.06&9.20&1.33&0.51&0.05
				\\ \hline
				\multirow{3}{*}{\ref{fn:f19}}&SG&7.15&-&8.21&-&2.62&-&3.73E+02&-&1&-&0.56&-
				\\*  \cline{2-14}
				&LHS&4.11&0.07&4.74&0.10&4.25&1.43&3.73E+02&0.30&1.60&0.49&0.88&5.68E-03
				\\*  \cline{2-14}
				&d-LHD&4.18&0.03&4.82&0.12&2.61&0.18&3.73E+02&8.42E-04&1&0&0.50&0.08
				\\ \hline
				\multirow{3}{*}{\ref{fn:f20}}&SG&6.19&-&8.09&-&4.77&-&3.73E+02&-&1&-&0.55&-
				\\*  \cline{2-14}
				&LHS&4.21&0.05&4.80&0.03&54.06&30.88&3.74E+02&0.74&25&10.73&0.89&0.02
				\\*  \cline{2-14}
				&d-LHD&4.20&0.07&4.47&0.40&1.07E+02&33.93&3.76E+02&2.14&77.80&26.98&0.45&3.29E-03
				\\ \hline
				\multirow{3}{*}{\ref{fn:f21}}&SG&1.99&-&2.03&-&0.04&-&15.85&-&1&-&0.50&-
				\\*  \cline{2-14}
				&LHS&3.42&0.04&3.95&0.06&0.04&0.02&15.85&2.88E-03&1.20&0.40&0.78&0.09
				\\*  \cline{2-14}
				&d-LHD&3.21&0.41&3.65&0.02&0.05&0.01&15.85&6.27E-04&1&0&0.43&0.02
				\\ \hline
				\multirow{3}{*}{\ref{fn:f22}}&SG&1.98&-&2.02&-&0.04&-&15.85&-&1&-&0.50&-
				\\*  \cline{2-14}
				&LHS&3.55&0.03&3.78&0.25&0.04&4.15E-03&15.85&3.33E-04&1&0&0.83&0.01
				\\*  \cline{2-14}
				&d-LHD&3.50&0.03&3.76&0.10&0.05&0.02&15.85&4.52E-04&1&0&0.51&0.03
				\\ \hline

			\end{longtable}
		\end{footnotesize}
	\end{center}
}

%% file: RationalApprox.bbl
\begin{thebibliography}{10}
\expandafter\ifx\csname url\endcsname\relax
  \def\url#1{\texttt{#1}}\fi
\expandafter\ifx\csname urlprefix\endcsname\relax\def\urlprefix{URL }\fi
\expandafter\ifx\csname href\endcsname\relax
  \def\href#1#2{#2} \def\path#1{#1}\fi

\bibitem{Booker1999}
A.~Booker, J.~{Dennis Jr}, P.~Frank, D.~Serafini, V.~Torczon, M.~Trosset, A
  rigorous framework for optimization of expensive functions by surrogates,
  Structural Multidisciplinary Optimization 17 (1999) 1--13.

\bibitem{Matheron1963}
G.~Matheron, {P}rinciples of geostatistics, Economic Geology 58 (1963)
  1246--1266.

\bibitem{Powell1992}
M.~Powell, Advances in Numerical Analysis, vol. 2: wavelets, subdivision
  algorithms and radial basis functions. Oxford University Press, Oxford, pp.
  105-210, Oxford University Press, London, 1992, Ch. The Theory of Radial
  Basis Function Approximation in 1990.

\bibitem{Friedman1991}
J.~Friedman, Multivariate adaptive regression splines, The Annals of Statistics
  19 (1991) 1--141.

\bibitem{Myers1995}
R.~Myers, D.~Montgomery, Response Surface Methodology, Process and Product
  Optimization using Designed Experiments, Wiley-Interscience Publication,
  1995.

\bibitem{Devore1986}
R.~Devore, X.-M. Yu, Multivariate rational approximation, Transactions of the
  American Mathematical Society 293~(1) (1986) 161--169.

\bibitem{New1964}
D.~J. Newman, Rational approximation to $|x|$, Michigan Math. J. 11 (1964)
  11--14.

\bibitem{GPT11}
P.~Gonnet, R.~Pach{\'o}n, L.~N. Trefethen, Robust rational interpolation and
  least-squares, Electron. Trans. Numer. Anal. 38 (2011) 146--167.

\bibitem{PGV12}
R.~Pach{\'o}n, P.~Gonnet, J.~Van~Deun, Fast and stable rational interpolation
  in roots of unity and {C}hebyshev points, SIAM J. Numer. Anal. 50~(3) (2012)
  1713--1734.
\newblock \href {http://dx.doi.org/10.1137/100797291}
  {\path{doi:10.1137/100797291}}.

\bibitem{Bra1986}
D.~Braess, Nonlinear Approximation Theory, Springer, Berlin, 1986.

\bibitem{Che1966}
E.~W. Cheney, Introduction to Approximation Theory, McGraw-Hill, New York,
  1966.

\bibitem{Tre2013}
L.~N. Trefethen, Approximation Theory and Approximation Practice, SIAM,
  Philadelphia, 2013.

\bibitem{Cuy1983}
A.~A.~M. Cuyt, Multivariate {P}ad{\' e}-approximants, J. Math. Anal. Appl. 96
  (1983) 283--293.
\newblock \href {http://dx.doi.org/10.1016/0022-247x(83)90041-0}
  {\path{doi:10.1016/0022-247x(83)90041-0}}.

\bibitem{CV1985}
A.~A.~M. Cuyt, B.~M. Verdonk, Multivariate rational interpolation, Computing 34
  (1985) 41--61.
\newblock \href {http://dx.doi.org/10.1007/bf02242172}
  {\path{doi:10.1007/bf02242172}}.

\bibitem{Cuyt2010}
A.~Cuyt, X.~Yang, A practical error formula for multivariate rational
  interpolation and approximation, Numerical Algorithms 55 (2010) 233--243.

\bibitem{SCGPS2013}
P.~Seshadri, P.~Constantine, P.~Gonnet, G.~Parks, S.~Shahpar, Stable
  multivariate rational interpolation for parameter-dependent aerospace models,
  in: 54th AIAA/ASME/ASCE/AHS/ASC Structures, Structural Dynamics, and
  Materials Conference, 2013.
\newblock \href {http://dx.doi.org/10.2514/6.2013-1680}
  {\path{doi:10.2514/6.2013-1680}}.

\bibitem{GP1973}
G.~H. Golub, V.~Pereyra, The differentiation of pseudo-inverses and nonlinear
  least squares problems whose variables separate, SIAM J. Numer. Anal. 10~(2)
  (1973) 413--432.

\bibitem{Bor2009}
C.~F. Borges, A full-{N}ewton approach to separable nonlinear least squares
  problems and its application to discrete least squares rational
  approximation, Electron. Trans. Numer. Anal. 35 (2009) 57--68.

\bibitem{NST2018}
Y.~Nakatsukasa, O.~S{\` e}te, L.~N. Trefethen, The {AAA} algorithm for rational
  approximation, SIAM J. Sci. Comput. 40 (2018) A1494--A1522.
\newblock \href {http://dx.doi.org/10.1137/16m1106122}
  {\path{doi:10.1137/16m1106122}}.

\bibitem{FNTB2018}
S.-I. Filip, Y.~Nakatsukasa, L.~N. Trefethen, B.~Beckermann, Rational minimax
  approximation via adaptive barycentric representations, SIAM J. Sci. Comput.
  40 (2018) A2427--A2455.
\newblock \href {http://dx.doi.org/10.1137/17m1132409}
  {\path{doi:10.1137/17m1132409}}.

\bibitem{Salazar2007}
O.~{Salazar Celis}, A.~Cuyt, B.~Verdonk, Rational approximation of vertical
  segments, Numerical Algorithms 45 (2007) 375--388.

\bibitem{Buckley:2011ms}
A.~Buckley, et~al., {General-purpose event generators for LHC physics}, Phys.
  Rept. 504 (2011) 145--233.
\newblock \href {http://arxiv.org/abs/1101.2599} {\path{arXiv:1101.2599}},
  \href {http://dx.doi.org/10.1016/j.physrep.2011.03.005}
  {\path{doi:10.1016/j.physrep.2011.03.005}}.

\bibitem{Alves:2017she}
J.~Albrecht, et~al., {A roadmap for HEP software and computing R\&D for the
  2020s}, Comput. Softw. Big Sci. 3~(1) (2019) 7.
\newblock \href {http://arxiv.org/abs/1712.06982} {\path{arXiv:1712.06982}},
  \href {http://dx.doi.org/10.1007/s41781-018-0018-8}
  {\path{doi:10.1007/s41781-018-0018-8}}.

\bibitem{GGT13}
P.~Gonnet, S.~G{\"u}ttel, L.~N. Trefethen, Robust {P}ad{\'e} approximation via
  {SVD}, SIAM Rev. 55~(1) (2013) 101--117.
\newblock \href {http://dx.doi.org/10.1137/110853236}
  {\path{doi:10.1137/110853236}}.

\bibitem{carter2000lebesgue}
M.~Carter, B.~van Brunt,
  \href{https://books.google.com/books?id=qgiLag9dPpMC}{The Lebesgue-Stieltjes
  Integral: A Practical Introduction}, Undergraduate Texts in Mathematics,
  Springer New York, 2000.
\newline\urlprefix\url{https://books.google.com/books?id=qgiLag9dPpMC}

\bibitem{HL2002}
M.~Huhtanen, R.~M. Larsen, On generating discrete orthgonal bivariate
  polynomials, BIT 42 (2002) 393--407.
\newblock \href {http://dx.doi.org/10.1023/A:1021907210628}
  {\path{doi:10.1023/A:1021907210628}}.

\bibitem{Zac2014}
M.~Zaccaron, Discrete orthogonal polynomials and hyperinterpolation over planar
  regions, {M.Sc.}\ thesis, University of Padova (2014).

\bibitem{CLS2010}
D.~Cox, J.~Little, D.~O'Shea, Ideals, Varieties, and Algorithms: An
  Introduction to Computational Algebraic Geometry and Commutative Algebra, 3rd
  Edition, Springer, New York, 2010.

\bibitem{GLRE2005}
L.~Giraud, J.~Langou, M.~Rozlo{\v z}n{\' i}k, J.~van~den Eshof, Rounding error
  analysis of the classical {G}ram-{S}chmidt orthogonalization process, Numer.
  Math. 101 (2005) 87--100.

\bibitem{LBG2013}
S.~J. Leon, {\AA}.~Bj{\" o}rck, W.~Gander, {G}ram-{S}chmidt orthogonalization:
  100 years and more, Numer. Linear Algebra Appl. 20 (2013) 492--532.
\newblock \href {http://dx.doi.org/10.1002/nla.1839}
  {\path{doi:10.1002/nla.1839}}.

\bibitem{Par1998}
B.~N. Parlett, The Symmetric Eigenvalue Problem, SIAM, Philadelphia, 1998.

\bibitem{stein2013bi}
O.~Stein, Bi-Level Strategies in Semi-Infinite Programming, Vol.~71, Springer
  Science \& Business Media, 2013.

\bibitem{HettichKortanek:1993}
R.~Hettich, K.~O. Kortanek, Semi-infinite programming: Theory, methods, and
  applications, SIAM Review 35~(3) (1993) 380--429.
\newblock \href {http://dx.doi.org/10.1137/1035089}
  {\path{doi:10.1137/1035089}}.

\bibitem{ben2009robust}
A.~Ben-Tal, L.~{El Ghaoui}, A.~Nemirovski, Robust Optimization, Princeton
  University Press, 2009.

\bibitem{LEVITIN19661}
E.~Levitin, B.~Polyak, Constrained minimization methods, USSR Computational
  Mathematics and Mathematical Physics 6~(5) (1966) 1 -- 50.

\bibitem{sherali.adams:98}
H.~Sherali, W.~Adams, A Reformulation-Linearization Technique for Solving
  Discrete and Continuous Nonconvex Problems, Kluwer, Dordrecht, 1998.

\bibitem{PhysRev.49.519}
G.~Breit, E.~Wigner,
  \href{https://link.aps.org/doi/10.1103/PhysRev.49.519}{Capture of slow
  neutrons}, Phys. Rev. 49 (1936) 519--531.
\newblock \href {http://dx.doi.org/10.1103/PhysRev.49.519}
  {\path{doi:10.1103/PhysRev.49.519}}.
\newline\urlprefix\url{https://link.aps.org/doi/10.1103/PhysRev.49.519}

\bibitem{Bohm:2004zi}
A.~R. Bohm, Y.~Sato, {Relativistic resonances: Their masses, widths, lifetimes,
  superposition, and causal evolution}, Phys. Rev. D71 (2005) 085018.
\newblock \href {http://arxiv.org/abs/hep-ph/0412106}
  {\path{arXiv:hep-ph/0412106}}, \href
  {http://dx.doi.org/10.1103/PhysRevD.71.085018}
  {\path{doi:10.1103/PhysRevD.71.085018}}.

\bibitem{Barthelmann2000}
V.~Barthelmann, E.~Novak, K.~Ritter, High dimensional polynomial interpolation
  on sparse grids, Advances in Computational Mathematics 12~(4) (2000)
  273--288.

\bibitem{MckayThreeMethods}
M.~D. McKay, R.~J. Beckman, W.~J. Conover, Comparison of three methods for
  selecting values of input variables in the analysis of output from a computer
  code, Technometrics 21~(2) (1979) 239--245.

\bibitem{laug}
E.~Anderson, Z.~Bai, C.~Bischof, S.~Blackford, J.~Demmel, J.~Dongarra,
  J.~Du~Croz, A.~Greenbaum, S.~Hammarling, A.~McKenney, D.~Sorensen, {LAPACK}
  Users' Guide, 3rd Edition, Society for Industrial and Applied Mathematics,
  Philadelphia, PA, 1999.

\bibitem{LEATHERMAN2017346}
E.~R. Leatherman, A.~M. Dean, T.~J. Santner,
  \href{http://www.sciencedirect.com/science/article/pii/S0167947316301773}{Designing
  combined physical and computer experiments to maximize prediction accuracy},
  Computational Statistics and Data Analysis 113 (2017) 346 -- 362.
\newblock \href {http://dx.doi.org/https://doi.org/10.1016/j.csda.2016.07.013}
  {\path{doi:https://doi.org/10.1016/j.csda.2016.07.013}}.
\newline\urlprefix\url{http://www.sciencedirect.com/science/article/pii/S0167947316301773}

\bibitem{Barlow:1993dm}
R.~J. Barlow, C.~Beeston, {Fitting using finite Monte Carlo samples}, Comput.
  Phys. Commun. 77 (1993) 219--228.
\newblock \href {http://dx.doi.org/10.1016/0010-4655(93)90005-W}
  {\path{doi:10.1016/0010-4655(93)90005-W}}.

\bibitem{PhysRevD.98.030001}
M.~Tanabashi, K.~Hagiwara, K.~Hikasa, K.~Nakamura, Y.~Sumino, F.~Takahashi,
  J.~Tanaka, K.~Agashe, G.~Aielli, C.~Amsler, M.~Antonelli, D.~M. Asner,
  H.~Baer, S.~Banerjee, R.~M. Barnett, T.~Basaglia, C.~W. Bauer, J.~J. Beatty,
  V.~I. Belousov, J.~Beringer, S.~Bethke, A.~Bettini, H.~Bichsel, O.~Biebel,
  K.~M. Black, E.~Blucher, O.~Buchmuller, V.~Burkert, M.~A. Bychkov, R.~N.
  Cahn, M.~Carena, A.~Ceccucci, A.~Cerri, D.~Chakraborty, M.-C. Chen, R.~S.
  Chivukula, G.~Cowan, O.~Dahl, G.~D'Ambrosio, T.~Damour, D.~de~Florian,
  A.~de~Gouv\^ea, T.~DeGrand, P.~de~Jong, G.~Dissertori, B.~A. Dobrescu,
  M.~D'Onofrio, M.~Doser, M.~Drees, H.~K. Dreiner, D.~A. Dwyer, P.~Eerola,
  S.~Eidelman, J.~Ellis, J.~Erler, V.~V. Ezhela, W.~Fetscher, B.~D. Fields,
  R.~Firestone, B.~Foster, A.~Freitas, H.~Gallagher, L.~Garren, H.-J. Gerber,
  G.~Gerbier, T.~Gershon, Y.~Gershtein, T.~Gherghetta, A.~A. Godizov,
  M.~Goodman, C.~Grab, A.~V. Gritsan, C.~Grojean, D.~E. Groom, M.~Gr\"unewald,
  A.~Gurtu, T.~Gutsche, H.~E. Haber, C.~Hanhart, S.~Hashimoto, Y.~Hayato, K.~G.
  Hayes, A.~Hebecker, S.~Heinemeyer, B.~Heltsley, J.~J. Hern\'andez-Rey,
  J.~Hisano, A.~H\"ocker, J.~Holder, A.~Holtkamp, T.~Hyodo, K.~D. Irwin, K.~F.
  Johnson, M.~Kado, M.~Karliner, U.~F. Katz, S.~R. Klein, E.~Klempt, R.~V.
  Kowalewski, F.~Krauss, M.~Kreps, B.~Krusche, Y.~V. Kuyanov, Y.~Kwon,
  O.~Lahav, J.~Laiho, J.~Lesgourgues, A.~Liddle, Z.~Ligeti, C.-J. Lin,
  C.~Lippmann, T.~M. Liss, L.~Littenberg, K.~S. Lugovsky, S.~B. Lugovsky,
  A.~Lusiani, Y.~Makida, F.~Maltoni, T.~Mannel, A.~V. Manohar, W.~J. Marciano,
  A.~D. Martin, A.~Masoni, J.~Matthews, U.-G. Mei\ss{}ner, D.~Milstead, R.~E.
  Mitchell, K.~M\"onig, P.~Molaro, F.~Moortgat, M.~Moskovic, H.~Murayama,
  M.~Narain, P.~Nason, S.~Navas, M.~Neubert, P.~Nevski, Y.~Nir, K.~A. Olive,
  S.~Pagan~Griso, J.~Parsons, C.~Patrignani, J.~A. Peacock, M.~Pennington,
  S.~T. Petcov, V.~A. Petrov, E.~Pianori, A.~Piepke, A.~Pomarol, A.~Quadt,
  J.~Rademacker, G.~Raffelt, B.~N. Ratcliff, P.~Richardson, A.~Ringwald,
  S.~Roesler, S.~Rolli, A.~Romaniouk, L.~J. Rosenberg, J.~L. Rosner, G.~Rybka,
  R.~A. Ryutin, C.~T. Sachrajda, Y.~Sakai, G.~P. Salam, S.~Sarkar, F.~Sauli,
  O.~Schneider, K.~Scholberg, A.~J. Schwartz, D.~Scott, V.~Sharma, S.~R.
  Sharpe, T.~Shutt, M.~Silari, T.~Sj\"ostrand, P.~Skands, T.~Skwarnicki, J.~G.
  Smith, G.~F. Smoot, S.~Spanier, H.~Spieler, C.~Spiering, A.~Stahl, S.~L.
  Stone, T.~Sumiyoshi, M.~J. Syphers, K.~Terashi, J.~Terning, U.~Thoma, R.~S.
  Thorne, L.~Tiator, M.~Titov, N.~P. Tkachenko, N.~A. T\"ornqvist, D.~R. Tovey,
  G.~Valencia, R.~Van~de Water, N.~Varelas, G.~Venanzoni, L.~Verde, M.~G.
  Vincter, P.~Vogel, A.~Vogt, S.~P. Wakely, W.~Walkowiak, C.~W. Walter,
  D.~Wands, D.~R. Ward, M.~O. Wascko, G.~Weiglein, D.~H. Weinberg, E.~J.
  Weinberg, M.~White, L.~R. Wiencke, S.~Willocq, C.~G. Wohl, J.~Womersley,
  C.~L. Woody, R.~L. Workman, W.-M. Yao, G.~P. Zeller, O.~V. Zenin, R.-Y. Zhu,
  S.-L. Zhu, F.~Zimmermann, P.~A. Zyla, J.~Anderson, L.~Fuller, V.~S. Lugovsky,
  P.~Schaffner,
  \href{https://link.aps.org/doi/10.1103/PhysRevD.98.030001}{Review of particle
  physics}, Phys. Rev. D 98 (2018) 030001.
\newblock \href {http://dx.doi.org/10.1103/PhysRevD.98.030001}
  {\path{doi:10.1103/PhysRevD.98.030001}}.
\newline\urlprefix\url{https://link.aps.org/doi/10.1103/PhysRevD.98.030001}

\bibitem{Hoferichter:2016nvd}
M.~Hoferichter, P.~Klos, J.~Menéndez, A.~Schwenk, {Analysis strategies for
  general spin-independent WIMP-nucleus scattering}, Phys. Rev. D94~(6) (2016)
  063505.
\newblock \href {http://arxiv.org/abs/1605.08043} {\path{arXiv:1605.08043}},
  \href {http://dx.doi.org/10.1103/PhysRevD.94.063505}
  {\path{doi:10.1103/PhysRevD.94.063505}}.

\bibitem{Cerdeno:2018bty}
D.~G. Cerdeño, A.~Cheek, E.~Reid, H.~Schulz, {Surrogate Models for Direct Dark
  Matter Detection}, JCAP 1808~(08) (2018) 011.
\newblock \href {http://arxiv.org/abs/1802.03174} {\path{arXiv:1802.03174}},
  \href {http://dx.doi.org/10.1088/1475-7516/2018/08/011}
  {\path{doi:10.1088/1475-7516/2018/08/011}}.

\bibitem{Feroz:2008xx}
F.~Feroz, M.~P. Hobson, M.~Bridges, {MultiNest: an efficient and robust
  Bayesian inference tool for cosmology and particle physics}, Mon. Not. Roy.
  Astron. Soc. 398 (2009) 1601--1614.
\newblock \href {http://arxiv.org/abs/0809.3437} {\path{arXiv:0809.3437}},
  \href {http://dx.doi.org/10.1111/j.1365-2966.2009.14548.x}
  {\path{doi:10.1111/j.1365-2966.2009.14548.x}}.

\bibitem{Feroz:2013hea}
F.~Feroz, M.~P. Hobson, E.~Cameron, A.~N. Pettitt, {Importance Nested Sampling
  and the MultiNest Algorithm}, Instrumentation and Methods for
  Astrophysics\href {http://arxiv.org/abs/1306.2144} {\path{arXiv:1306.2144}}.

\bibitem{pymultinest}
J.~{Buchner}, A.~{Georgakakis}, K.~{Nandra}, L.~{Hsu}, C.~{Rangel},
  M.~{Brightman}, A.~{Merloni}, M.~{Salvato}, J.~{Donley}, D.~{Kocevski},
  {X-ray spectral modelling of the AGN obscuring region in the CDFS: Bayesian
  model selection and catalogue}, aap 564 (2014) A125.
\newblock \href {http://arxiv.org/abs/1402.0004} {\path{arXiv:1402.0004}},
  \href {http://dx.doi.org/10.1051/0004-6361/201322971}
  {\path{doi:10.1051/0004-6361/201322971}}.

\bibitem{Fowlie:2016hew}
A.~Fowlie, M.~H. Bardsley, {Superplot: a graphical interface for plotting and
  analysing MultiNest output}, Eur. Phys. J. Plus 131~(11) (2016) 391.
\newblock \href {http://arxiv.org/abs/1603.00555} {\path{arXiv:1603.00555}},
  \href {http://dx.doi.org/10.1140/epjp/i2016-16391-0}
  {\path{doi:10.1140/epjp/i2016-16391-0}}.

\bibitem{ts:05}
M.~Tawarmalani, N.~V. Sahinidis, {A polyhedral branch-and-cut approach to
  global optimization}, Mathematical Programming 103 (2005) 225--249.

\bibitem{sahinidis:baron:17.8.9}
N.~V. Sahinidis, {BARON 17.8.9: Global Optimization of Mixed-Integer Nonlinear
  Programs, {\em User's Manual}} (2017).

\end{thebibliography}
